\documentclass[preprint,12pt]{elsarticle}

\usepackage{hyperref}
\usepackage{threeparttable}

\usepackage{amssymb}
\usepackage{amsmath}
\usepackage[numbers]{natbib}
\usepackage{subcaption}
\usepackage{graphics} 
\usepackage{times} 
\usepackage{amsmath} 
\usepackage{amssymb}  
\usepackage{amsthm}
\usepackage{amsfonts}
\usepackage{mathtools}
\usepackage{multirow} 
\usepackage{rotating} 
\usepackage{algorithm}
\usepackage{algorithmic}
\usepackage{tikz}
\usetikzlibrary{calc}
\usetikzlibrary{trees}
\usetikzlibrary{arrows.meta, positioning}
\usepackage{pdflscape}
\usepackage{graphicx}
\usepackage{epsfig} 
\usepackage{booktabs}
\usepackage{float}
\usepackage{flushend}
\usepackage[table,xcdraw]{xcolor}
\usepackage{csquotes}
\usepackage{bm}
\usepackage{relsize}
\usepackage{colortbl}
\usepackage{lipsum}
\usepackage[mathscr]{eucal}
\usepackage{fancyhdr}
\usepackage{siunitx}
\newcommand{\best}[1]{\textbf{#1}}

\newcommand{\norm}[1]{\left\lVert#1\right\rVert}

\usepackage{tikz}

\usepackage[pagewise]{lineno} 

\journal{Transportation Research Part C: Emerging Technologies}

\begin{document}

\begin{frontmatter}



\title{Integrated Prediction and Distributionally Robust Optimization for Air Traffic Management}


\author[inst1]{
Haochen Wu
}
\author[inst2]{
Xinting Zhu
} 
\author[inst2]{
Lishuai Li
} 
\author[inst1]{
Max Z. Li
} 
\affiliation[inst1]{organization={University of Michigan, Aerospace Engineering},
            addressline={500 S State St}, 
            city={Ann Arbor},
            postcode={48109}, 
            state={Michigan},
            country={United States of America}}

\affiliation[inst2]{organization={City University of Hong Kong\\Kowloon Tong, Kowloon TongKowloon Tong\\},
            addressline={83 Tat Chee Ave}, 
            city={Kowloon Tong},
            postcode={523808}, 
            state={Hong Kong}
}

\begin{highlights}
\item Propose an integrated prediction and distributionally robust optimization framework for Air Traffic Management. 
\item Show that leveraging upstream predictions significantly reduce system-wide delay costs. 
\item Demonstrate that distributionally robust optimization models improve performance under distribution shifts relative to upstream predictions. 
\end{highlights}

\begin{abstract}
Strategic Traffic Management Initiatives (TMIs) such as Ground Delay Programs (GDPs) play a crucial role in mitigating operational costs associated with air traffic demand-capacity imbalances. However, GDPs can only be planned (e.g., duration, delay assignments) with confidence if the future capacities at constrained resources (i.e., airports) are predictable. In reality, such future capacities are uncertain, and predictive models may provide predictions that are vulnerable to errors and distribution shifts. Motivated by the goal of planning optimal GDPs that are \emph{distributionally robust} against airport capacity prediction errors, we study a fully integrated learning-driven optimization framework. We design a deep learning-based prediction model capable of forecasting arrival and departure capacity distributions across a network of airports. We incorporate the predictions into a distributionally robust formulation of the multi-airport ground holding program (\textsc{dr-MAGHP}). Our results demonstrate that \textsc{dr-MAGHP} can achieve up to a $15.6\%$ improvement over the stochastic programming formulation (\textsc{sp-MAGHP}) under airport capacity distribution shifts. We conclude by outlining future research directions aimed at enhancing both the learning and optimization components of the framework.
\end{abstract}
\end{frontmatter}
\section{Introduction}
\label{sec:intro}
Congestion in air transportation systems results from demand-capacity imbalances, often due to airport capacity reductions. Within the US National Airspace System (NAS), the strategic implementation of Traffic Management Initiatives (TMIs) seeks to reduce the operational costs of such imbalances. A prominent example of TMIs is Ground Delay Programs (GDPs), which aims to  delay flights on the ground at the origin airport and mitigate costly airborne delays in response to arrival capacity reductions at the destination airport. The optimal implementation of GDPs is the purview of airport ground holding optimization problems, or GHPs. 
Airport capacities at future time periods play a significant role in GDP implementations. If such capacities are known, the optimal delay allocation decisions can be found by solving GHPs \cite{GHP1}. However, in practice, it is extremely difficult for traffic management decision-makers to ascertain future airport capacities (i.e., Airport Arrival Rates and Airport Departure Rates) due to myriad uncertainties. Such uncertainties stem from environmental (e.g., convective weather predictions \cite{kicinger2016airport}) and operational (e.g., runway availability and traffic volume \cite{tittle2013airport}) sources. Although a large variety of prior work focuses on, e.g., airport runway configuration prediction \cite{avery2015predicting}, weather impact predictions \cite{kicinger2016airport}, and TMI implementation predictions \cite{vlachou2019simultaneous}, all predictions result in a probability distribution over potential outcomes. If such predictions were to be incorrect (e.g., due to distribution shifts or misspecification), the resultant GDP implementation may be sub-optimal. 

\subsection{Motivation and research problem}
\label{subsec:motivation}
In this work, we are motivated by the dual objective of harnessing advances in machine learning (ML) for airport capacity prediction while maintaining a \emph{cautiously optimistic} stance in prescribing ground delay policies. On the one hand, ML models enable data-driven forecasting, which aligns with many future concepts of operations for air traffic management such as the Federal Aviation Administration’s Information-Centric NAS vision \cite{Info_centric_NAS} and SESAR’s European ATM Master Plan 2025 Edition \cite{SESAR_ATM_MASTER_PLAN}. On the other hand, robust decision-making remains essential, as such predictions are inherently imperfect.

ML has seen rapid progress and growing applications in aviation, particularly in airport and airspace operations. By leveraging historical and real-time data, ML models have been used to anticipate and mitigate delays, predict traffic flows, estimate fuel consumption and airport throughput, and optimize flight paths \cite{zhu2021flight,zhu2022short,Wang2022}, thereby enhancing overall airspace efficiency. Programs such as SESAR in Europe and NextGen in the United States are actively working to integrate ML into advanced air traffic management (ATM) systems \cite{Info_centric_NAS,SESAR_ATM_MASTER_PLAN}.

However, the practical use of ML predictions in downstream optimization can be problematic when decision models are highly sensitive to forecast errors. A central challenge is \emph{distribution shift}, i.e., when the underlying data distribution differs between training and deployment \cite{tamang2025handling}. Such shifts often degrade decision-making performance, as illustrated by autonomous vehicles trained in sunny conditions performing poorly in adverse weather \cite{filos2020can}, recommendation systems facing seasonal changes in user preferences \cite{yang2023generic}, or fraud and spam detection models being circumvented through adversarial concept drift \cite{wang2019drifted}. In the context of airport capacity prediction, distribution shifts may lead to inaccurate predictions and, in turn, suboptimal or even unsafe ground holding policies.

To address these challenges, GDPs formulations must explicitly account for uncertainty in upstream predictions. Traditional optimization approaches include the deterministic multi-airport ground holding problem (MAGHP), which assumes fixed capacities, and the stochastic MAGHP, which relaxes this assumption by modeling capacity as a random variable. In this work, we advance this line of research by developing a \emph{distributionally robust} formulation, termed the distributionally robust multi-airport ground holding problem (\textsc{dr-MAGHP}). Distributionally robust optimization (DRO) seeks solutions that perform well under the worst-case distribution within a predefined \emph{ambiguity set}, which in our case is constructed based on the predicted distribution of airport capacities generated by the ML model. Section~\ref{sec:method} provides a rigorous overview of these concepts.

\subsection{Contributions} \label{sec:contribution}
In this work, our contributions are as follows:
\begin{enumerate}
    \item We develop a deep learning framework to provide practical probabilistic predictions of departure and arrival capacity at each airport. Based on experiments with the 30 major US core airports, our model demonstrates that the framework is effective in predicting capacity distributions across diverse airport environments and operational characteristics.
    \item We develop a distributionally robust multi-airport ground holding problem (\textsc{dr-MAGHP}) that incorporates Wasserstein ambiguity sets constructed from upstream distributional predictions. We derive a tractable reformulation of \textsc{dr-MAGHP} and introduce a scenario reduction technique to mitigate the computational challenges posed by the large-scale scenario trees arising from time-series capacity predictions.  
    \item We perform a comprehensive sensitivity analysis to evaluate the performance of \textsc{dr-MAGHP} under varying degrees of airport capacity reductions, meant to mimic distributional uncertainty and mispredictions. This analysis highlights how the choice of the \emph{size} of the ambiguity sets influences robustness and cost efficiency, thereby providing insights into the trade-offs involved in robust policy design. 
\end{enumerate}
We note that a preliminary version of this work has appeared in the 11\textsuperscript{th} International Conference on Research in Air Transportation\cite{wu2024distributionally}.

\section{Literature review} \label{subsec:prior_works}

\subsection{Distributional prediction for airport capacity}

Airport capacity can be defined as the maximum sustainable throughput for arriving (the Airport Arrival Rate, or AAR) and departing (the Airport Departure Rate, or ADR) flights \cite{kicinger2012airport}. In contrast to declared capacities obtained from theoretical analyses or statistical approaches \cite{gilbo1993airport}, real-time capacity is dynamic and challenging to predict in advance. It is influenced by several interconnected operational and environmental factors \cite{choi2021artificial}. As the prediction time horizon increases, forecast uncertainty grows as well, rendering accurate long-term predictions difficult \cite{tien2018using}.

Traditional prediction models in aviation have predominantly relied on deterministic approaches. A sampling of previous work includes analytical approaches such as the Integrated Airport Capacity Model (IACM) \cite{kicinger2012airport}, and more recently, data-driven approaches such as the AAR Distribution Prediction Model (ADPM) \cite{cox2016probabilistic}. However, these deterministic methods often fail to capture the inherent uncertainties in air traffic systems, leading to suboptimal decision-making under uncertain conditions.

The emergence of probabilistic and distributional prediction methods has addressed these limitations by explicitly quantifying prediction uncertainties. 
A non-exhaustive set of recent approaches utilizes probabilistic methods to account for uncertainties from various sources in aircraft trajectory prediction \cite{pang2021data}, delay predictions \cite{zoutendijk2021probabilistic}, and traffic flow prediction \cite{masalonis2004using,zhang2023air}, thereby enabling more robust air traffic management decisions. 
Quantile regression techniques have been developed to capture the full delay distribution rather than simple point estimates, with models demonstrating superior performance compared to statistical baselines by learning various quantiles of departure and arrival delay distributions using features available several days before operations during the pretactical phase~\cite{dalmau2024probabilistic}. 
Ensemble methods, particularly Variational LSTM models incorporating Monte Carlo Dropout, achieve median absolute errors in delay predictions of 5.8 minutes per day across multiple airports while providing well-calibrated prediction intervals crucial for risk management applications~\cite{vandal2018prediction}. 
Bayesian approaches, including Bayesian structural time series models for air traffic demand forecasting~\cite{rodriguez2024air}, demonstrate efficacy in incorporating external factors and handling non-linear demand patterns.
Additionally, probabilistic aircraft trajectory prediction approaches using Bayesian Neural Networks \cite{pang2021data,pang2020probabilistic} can predict aircraft flight paths while quantifying the uncertainties in those predictions, which is particularly valuable for managing weather-related disruptions. 
Separately, Monte Carlo simulation techniques \cite{lecchini2006monte} enable evaluation of post-resolution conflict probabilities, supporting optimization under uncertainty for air traffic conflict resolution.

Few works have been conducted on predicting real-time airport capacities \cite{wang2021prediction}, and flight delay distribution predictions~\cite{Wang2022}. The field continues to evolve with increasing integration of probabilistic methods, advanced deep learning architectures, and multi-source data fusion approaches.

\subsection{Mathematical Programs in Air Traffic Management}

Mathematical programming, particularly integer programming (IP), has been extensively applied in Air Traffic Management (ATM) to address challenges arising from capacity-constrained airspace and airports. One central problem is the Air Traffic Flow Management Problem (ATFMP), which seeks to optimize aircraft flows and reduce congestion while respecting capacity limits. Odoni \cite{odoni1987flow} was among the first to formalize the Flow Management Problem (FMP), modeling it as a network-based formulation to manage congestion across airports, airways, waypoints, and sectors. Building on this, Bertsimas and Patterson \cite{bertsimas2008air} developed an IP model for ATFMP that minimizes total delay costs by optimizing flight arrival and departure times within airspace and airport capacity constraints. Addressing dynamic routing under uncertain weather conditions, Bertsimas and Sim \cite{bertsimas2000traffic} proposed a dynamic multicommodity network flow model, which they solved using a Lagrangian generation algorithm. In a related approach, Sun et al.\ \cite{sun2008multicommodity} introduced a cell-transmission model that recasts air traffic flow as a linear time-invariant dynamical system, solved using a relaxed integer program to manage sector-level aircraft counts over time.

A key subclass of ATFMP is the Ground Holding Problem (GHP), which specifically addresses delays caused by limited capacity at departure and arrival airports. Unlike ATFMP, which accounts for system-wide capacity constraints, GDP models focus on airport-level restrictions. The Multi-Airport Ground-Holding Problem (MAGHP) and its variations—both deterministic and stochastic—have been widely studied \cite{GHP1}. Stochastic formulations incorporate two-stage stochastic programming and chance-constrained programming to model airport capacities as random variables \cite{GHP2,GHP4}. More recent advances include data-driven control frameworks that minimize delay costs while redistributing delays spatially to improve operational resilience \cite{GHP5}. Additional research introduces fairness- and passenger-oriented metrics to ensure equitable treatment of flights across airlines and time slots \cite{GHP6,GHP7}.

Another critical ATM problem is airport slot allocation, especially at highly congested hubs where demand exceeds available runway capacity. This involves assigning takeoff and landing slots to airlines in a manner that adheres to regulatory and operational constraints \cite{IATA_WASG_Ed3_2024} while balancing efficiency and fairness. Pellegrini et al.\ introduced the Simultaneous Slot Allocation Problem (SSAP) and proposed metaheuristic algorithms that outperform traditional ILP methods on large-scale instances by efficiently handling inter-airport coordination \cite{pellegrini2012metaheuristic}. They later extended this work through the SOSTA (Simultaneous Optimisation of the airport Slot Allocation) model, which incorporates aircraft turnaround constraints and conforms to European regulatory standards, demonstrating high computational efficiency (i.e., shorter solution times) on real operational data \cite{pellegrini2017sosta}. From a market design perspective, Ball et al.\ \cite{ball2020quantity} developed a quantity-contingent auction framework in which slot values vary with allocation volume, incorporating constraints to mitigate the risk of market power abuse. Zografos and Jiang \cite{zografos2019bi} introduced a bi-objective optimization model that jointly optimizes schedule efficiency and fairness, providing trade-off analyses under different regulatory regimes. Addressing computational scalability, Ribeiro et al.\ \cite{ribeiro2019large} proposed a large-scale neighborhood search algorithm that delivers near-optimal slot assignments for busy airports with significant runtime improvements over conventional ILP methods.

\subsection{Distributionally Robust Optimization}

Distributionally robust optimization (DRO) provides a robust framework for decision-making under uncertainty by optimizing decisions against the worst-case distribution within a specified ambiguity set. Delage and Ye \cite{DR1} introduced a foundational DRO framework based on moment information, where ambiguity sets are defined using means and covariances, and tractable conic reformulations are derived with statistical guarantees. Esfahani and Kuhn \cite{DR2} extended this approach by defining ambiguity sets using the Wasserstein distance, offering finite-sample performance guarantees and tractable convex reformulations grounded in duality. Kim \cite{DR3} focused on solving two-stage distributionally robust mixed-integer programs under Wasserstein ambiguity using a dual decomposition approach, demonstrating its scalability through discretization and Lagrangian relaxation. Cheramin et al.\ \cite{cheramin2022computationally} proposed computationally efficient approximations for DRO models that combine moment-based and Wasserstein-based ambiguity, significantly improving tractability for large-scale problems. Jiang et al.\ \cite{jiang2019data} applied Wasserstein DRO to appointment scheduling under service time and no-show uncertainty, using copositive and linear programming reformulations to achieve strong out-of-sample performance. Hanasusanto and Kuhn \cite{hanasusanto2018conic} developed conic programming reformulations for two-stage DRO problems over Wasserstein balls, enabling tractable multistage decision-making. Lastly, Wiesemann, Kuhn, and Sim \cite{wiesemann2014distributionally} presented a unifying DRO framework for convex optimization problems, covering a wide range of ambiguity sets and establishing general conditions for tractable reformulations.

\section{Methodology} \label{sec:method}
\begin{figure*}[htbp]
    \centering
    \includegraphics[width=\textwidth]{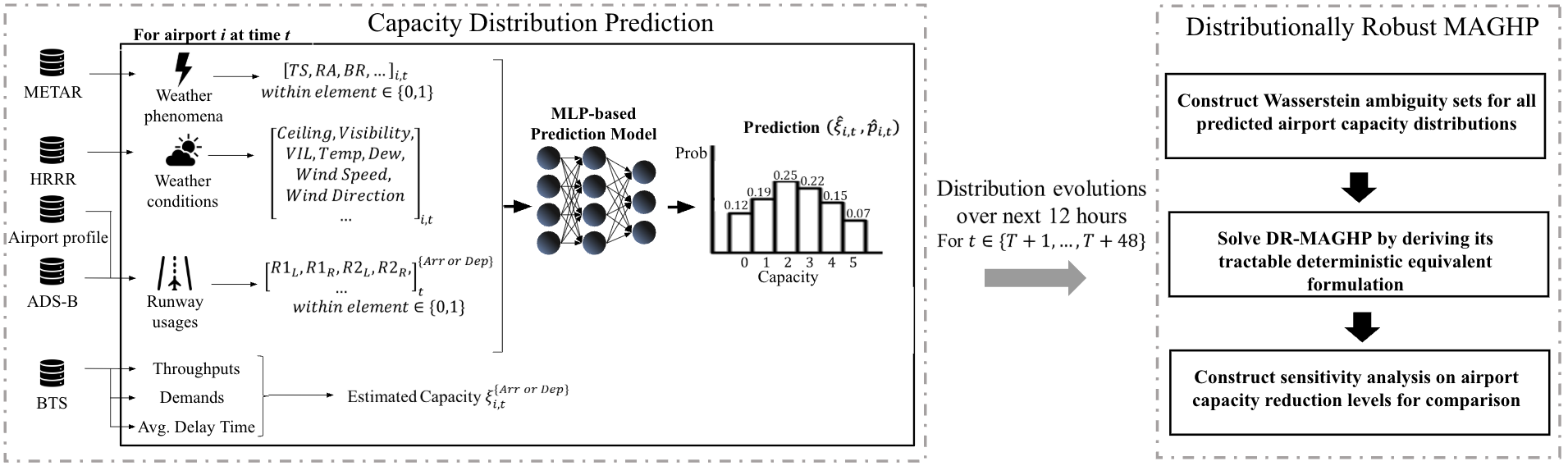}
    \caption{Learning-driven airport capacity distribution prediction and distributionally robust GDP optimization framework. 
    }
    \label{fig:capacity_distribution_prediction_framework} 
\end{figure*}

Figure~\ref{fig:capacity_distribution_prediction_framework} illustrates the proposed learning-driven optimization framework. At each decision epoch $t$, the capacity prediction model generates a probability distribution, where the support corresponds to the set of possible capacity outcomes $\widehat{\xi}_{t}$ and the associated probabilities are denoted by $\widehat{p}_{t}$. This predicted distribution then serves as input parameters for the downstream \textsc{dr-MAGHP}, which leverages them to derive robust ground delay policies. In what follows, we first describe the development of the capacity prediction model, including the derivation of actual capacities (Section~\ref{sssec:derive_actual_capacities}) and the architecture of the probabilistic prediction model (Section~\ref{sssec:prob_cap_pred}). We then explain how \textsc{dr-MAGHP} incorporates these predictive distributions into the decision-making process, focusing on the construction of Wasserstein ambiguity sets (Section~\ref{ssec:dr-MAGHP}) and the derivation of a tractable reformulation (Section~\ref{ssec:scenario_reduce}).







\subsection{Airport capacity distribution prediction}   
\label{ssec:airport_cap_dis_pred} 
We develop a deep learning model for airport capacity distribution prediction. The model predicts airport arrival or departure capacity distributions across a 12-hour prediction horizon, discretized into 15-minute intervals. 
The 12-hour horizon is selected based on operational requirements in air traffic management, where medium-term capacity predictions are critical for strategic flow management and ground delay program decisions \cite{ball2010total, mukherjee2012dynamic}. 
Additionally, this horizon aligns with typical airport coordination timescales and provides sufficient lead time for proactive traffic management while maintaining reasonable prediction accuracy \cite{wang2017airport}. 
The 15-minute interval granularity balances computational efficiency with operational relevance, as it corresponds to standard air traffic management coordination intervals \cite{FAA2016FSM,FAA2017AADC} and captures meaningful variations in airport capacity without excessive noise from minute-by-minute fluctuations \cite{gilbo1993airport, sridhar2011modeling}.

To account for airport-specific operational characteristics and infrastructure differences \cite{liu2019using}, we build separate prediction models for each airport. Furthermore, for each airport, we develop independent models for arrival and departure capacity predictions. We acknowledge that this is a simplifying assumption, as arrival and departure operations are inherently interdependent due to shared runway usage, taxiway conflicts, and gate constraints \cite{gilbo1997airport, simaiakis2014balancing}. 
In practice, high arrival volumes can reduce departure capacity and vice versa, particularly at airports with intersecting or closely spaced parallel runways. However, modeling these complex interdependencies would significantly increase model complexity and data requirements. This independence assumption allows for more tractable model development while still capturing the primary capacity dynamics for each operation type. Future work could explore joint arrival-departure capacity prediction models or post-processing techniques to ensure consistency between arrival and departure predictions.

\subsubsection{Deriving actual capacities from throughput} \label{sssec:derive_actual_capacities}



Training and validating the prediction model requires observations of actual airport capacity values. However, a fundamental challenge in airport capacity prediction is that true operational capacity cannot be directly observed from historical data \cite{ramanujam2009estimation}. 
It is crucial to distinguish between \enquote{throughput} and \enquote{capacity} in aviation operations: \emph{throughput} represents the actual number of aircraft operations (arrivals or departures) processed during a given time period, while \emph{capacity} represents the maximum number of operations the airport system can theoretically handle under specific conditions \cite{neufville2013airport}.
Historical throughput data reflects realized demand constrained by both capacity limits and actual demand levels—when demand is low, observed throughput will be below the airport's true capacity. Conversely, capacity represents the system's theoretical maximum processing rate, which may not be fully utilized due to demand variations, schedule spacing, or operational inefficiencies.

To address this fundamental measurement challenge, we employ a systematic rule-based methodology to estimate capacity by identifying \emph{capacity-saturated} periods where observed throughput accurately reflects true operational capacity. During these time periods $t$, we assume that $\widehat{\text{capacity}_t} = \text{throughput}_t$, meaning the airport is operating at or very close to its maximum capacity limit and demand exceeds the available capacity.

Our approach leverages multiple operational indicators to comprehensively identify capacity-saturated conditions. We define three complementary criteria based on key airport performance metrics: throughput levels, demand-supply relationships, and delay characteristics. For each airport $i$ at time $t$, we apply the following decision rule separately for departure and arrival operations:

\begin{equation}
\label{eq:capacity_estimation}
    \widehat{\text{capacity}}_{i,t} = \text{throughput}_{i,t} \iff \text{(Criterion 1)} \,\cup\, \text{(Criterion 2)} \,\cup\, \text{(Criterion 3)},
\end{equation}
where the three criteria are defined as:

\begin{enumerate}
    \item \emph{Criterion 1}---{Throughput-based criterion:}
    \begin{equation}
        \text{Throughput Ratio} = \frac{\text{actual throughput}_{i,t}}{\text{threshold}_{i}} \geq 1,
    \end{equation}
    where $\text{threshold}_{i}$ refers to the value in the 90\textsuperscript{th} percentile of the throughput distribution for airport $i$.
    
    \item \emph{Criterion 2}---{Demand-based criterion:}
    \begin{equation}
        \text{Throughput-Demand Ratio} = \frac{\text{actual throughput}_{i,t}}{\text{scheduled demand}_{i,t}} \leq \alpha,
    \end{equation}
    where $\alpha$ is a threshold parameter (we use $\alpha = 0.8$ in our implementation)  that identifies periods when the airport cannot fully satisfy scheduled demand. This indicates that meaningful demand pressure exists where scheduled flight demand exceeds the airport's processing capacity, pushing the system toward its operational limits. We note that we did not conduct a comprehensive sensitivity analysis on this parameter, which represents a limitation of our current approach that could be addressed in future work.
    
    \item \emph{Criterion 3}---{Delay-based criterion:}
    \begin{equation}
        \left(\text{Average delay}_{i,t} \geq 15 \text{ mins}\right) \cap \left(\text{Delayed flights}_{i,t} \geq 2\right),
    \end{equation}
    where delayed flights are defined as those experiencing more than 5 minutes of delay relative to the schedule. The 15-minute average delay threshold is defined to reflect significant operational stress indicative of capacity saturation. The requirement for at least 2 delayed flights ensures statistical significance and avoids classification based on isolated incidents.
    We acknowledge that these threshold parameters (15 minutes for average delay, 2 flights for minimum count, and 5 minutes for individual flight delay classification) were selected based on operational judgment rather than comprehensive sensitivity analysis.

\end{enumerate}

This multi-criteria framework ensures robust identification of capacity-saturated periods by requiring when any one of three criteria is satisfied as shown in \eqref{eq:capacity_estimation}. The throughput criterion captures high activity levels, the demand criterion ensures meaningful operational pressure, and the delay criterion confirms that the airport is experiencing capacity-related stress. Only data from time periods satisfying all criteria are used for model training and validation, maximizing the chances that our capacity estimates reflect true operational limits rather than demand-constrained throughput values.

While this approach provides a practical solution for capacity estimation, we acknowledge that it may underestimate true capacity during periods of severe weather or other operational disruptions when capacity limits change but our historical thresholds may not capture these variations. However, alternative approaches face similar challenges\textemdash FAA-published AARs and ADRs are conservative "called rates" set for traffic management rather than reflecting absolute limits. AARs and ADRs are frequently adjusted based on anticipated conditions and risk tolerance, potentially underestimating the airport's true maximum throughput capability under optimal conditions.

Therefore, in the absence of a standardized ground truth for airport capacities, we consider our approach to be a reasonable benchmark for identifying periods when airports are operating near their practical limits under prevailing conditions.



\subsubsection{Distributional capacity prediction} \label{sssec:prob_cap_pred}

Our model includes input feature engineering, a multilayer perceptron (MLP)-based prediction model, and output reformulation. The MLP architecture consists of three layers with 17, 32, and ${MAX}_z + 1$ neurons respectively, where ${MAX}_z$ represents the historical maximum capacity of airport $z$. The network employs ReLU activation functions in the hidden layers, produces probabilistic outputs through a softmax activation function in the final layer, and is trained using cross-entropy loss to optimize the capacity distribution predictions.

We denote the predicted capacity distributions for arrivals and departures at airport $z$ and time $t$ by
\[
\left(\widehat{\xi}^{(z,a)}_{t},\, \widehat{p}^{(z,a)}_{t}\right)
\quad\text{and}\quad
\left(\widehat{\xi}^{(z,g)}_{t},\, \widehat{p}^{(z,g)}_{t}\right),
\]
respectively. We learn a mapping $\Psi^{z,\bullet}$ (with $\bullet\in\{a,g\}$) from inputs $X_t^z$ to a predictive distribution:
\[
\left(\widehat{\xi}^{(z,\bullet)}_{t},\, \widehat{p}^{(z,\bullet)}_{t}\right)
= \Psi^{z,\bullet}\!\left(X_t^z\right).
\]

Rather than predicting a single point estimate of capacity, our model outputs a complete probability distribution over all possible capacity values. This distributional approach is essential for the subsequent distributionally robust optimization (DRO) framework, which requires the full uncertainty characterization of capacity predictions rather than point estimates. Given the inherent aleatoric uncertainty in airport operations due to factors such as individual air traffic manager decisions, aircraft performance variations, and dynamic operational constraints, capacity prediction naturally exhibits distributional characteristics that cannot be adequately captured by point predictions.

To enable this distributional prediction, we reformulate the capacity prediction as a multi-class classification problem. We use a categorical encoding where $\widehat{p}^{(z, \bullet)}_{t}$ is represented as a probability vector over all possible capacity values. This process converts scalar capacity values into a probability distribution across the discrete capacity space. For example, as shown in Figure \ref{fig:capacity_distribution_prediction_framework}, instead of predicting a single capacity value of 2, the model outputs a probability distribution such as ${[0.05, 0.10, 0.70, 0.10, 0.03, 0.02]}_t$, where each element represents the probability that the corresponding capacity value occurs. The length of this vector corresponds to the range of capacity values observed at this airport historically (${MAX}_z + 1$ possible values from 0 to ${MAX}_z$).
During training, the ground truth capacity observations (derived from the capacity-saturated periods as described in Section \ref{sssec:derive_actual_capacities}) are encoded using one-hot vectors. For instance, an observed capacity of 2 would be represented as the sparse vector ${[ 0, 0, 1, 0, 0, 0 ]}_t$. However, the model's output during inference is a complete probability distribution over all capacity values, providing the distributional information required for \textsc{dr-MAGHP}.

Our model incorporates two categories of input features for $X_t$: runway configurations and meteorological information. The runway configuration features capture the operational setup of the airport, including active runway assignments and operational modes, which directly impact the airport's capacity limits. The meteorological features are converted and vectorized, including seven weather variables: ceiling, visibility, vertically integrated liquid (VIL), temperature, dew points, wind direction, surface wind speeds, and one phenomenon indicator of adverse weather. The selection of these features is based on their studied impacts on airport capacity (e.g., see \cite{allan2001analysis,renhe2014meteorological}). By combining both runway configuration and weather information, the model can capture the primary operational and environmental factors that determine airport capacity under different conditions.

\subsection{Distributionally robust MAGHP (\textsc{dr-MAGHP})} \label{ssec:dr-MAGHP}
\begin{figure}
\centering
\subfloat[Gaussian]{\includegraphics[width=4.5cm]{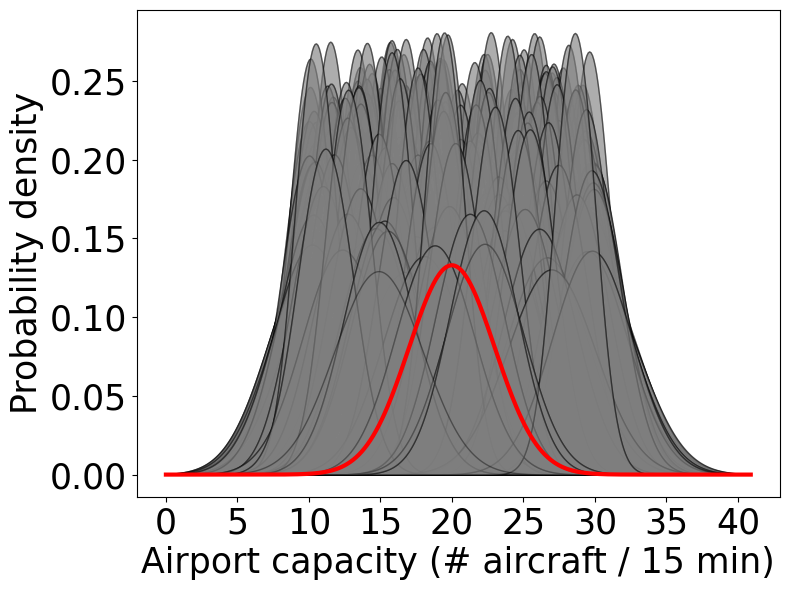}}\hfill
\subfloat[Gamma]{\includegraphics[width=4.5cm]{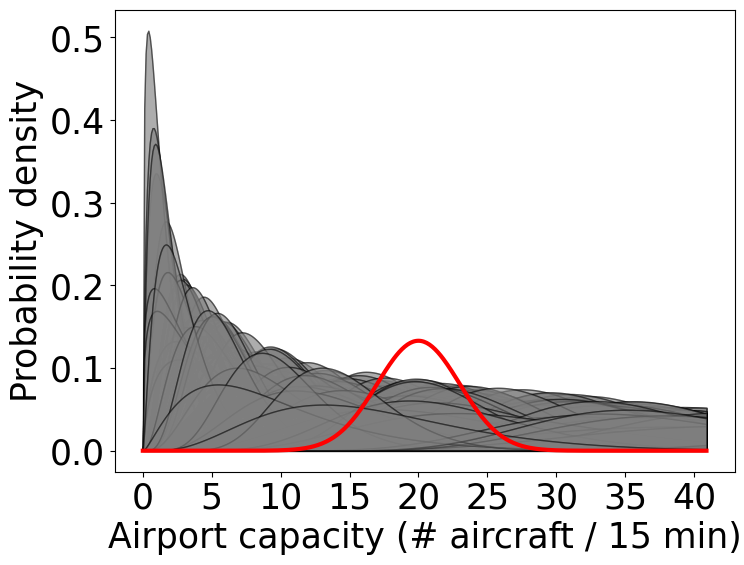}}\hfill
\subfloat[Erlang]{\includegraphics[width=4.5cm]{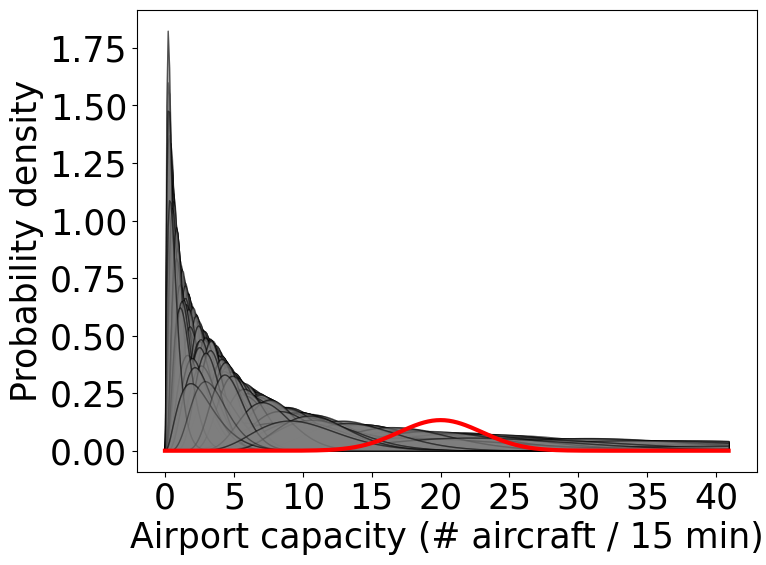}}\vfill
\subfloat[laplace]{\includegraphics[width=4.5cm]{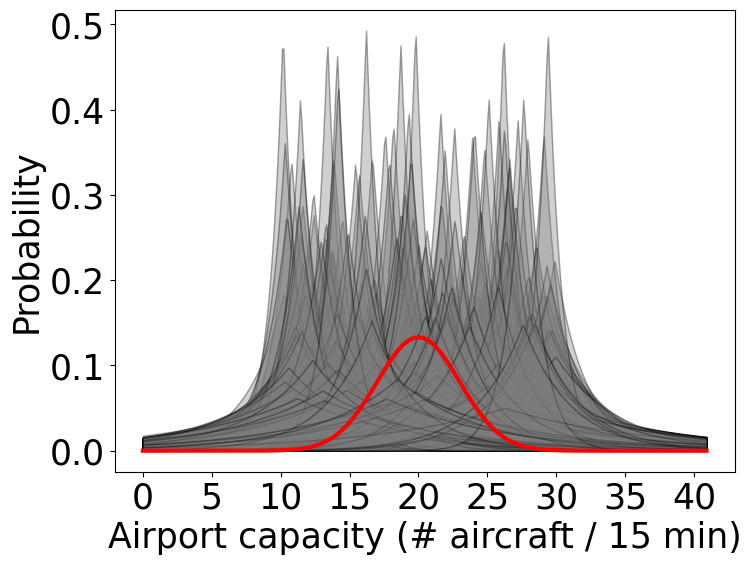}}\hfill
\subfloat[Lognorm]{\includegraphics[width=4.5cm]{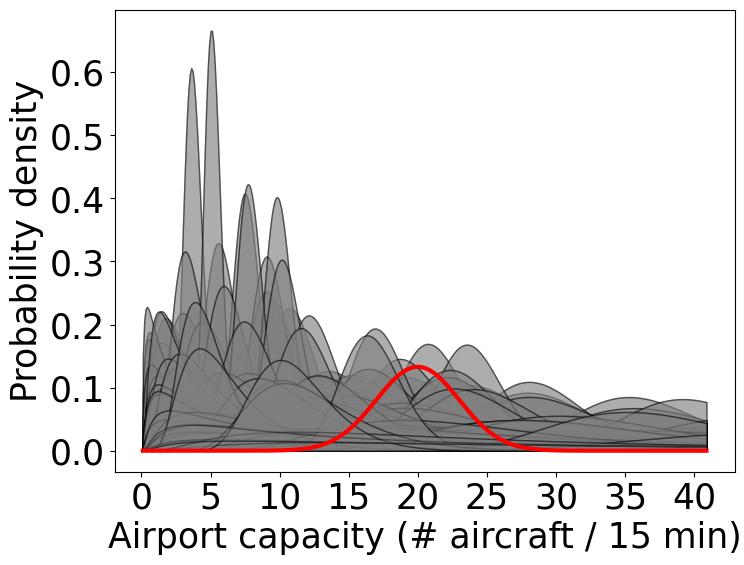}}\hfill
\subfloat[Weibull]{\includegraphics[width=4.5cm]{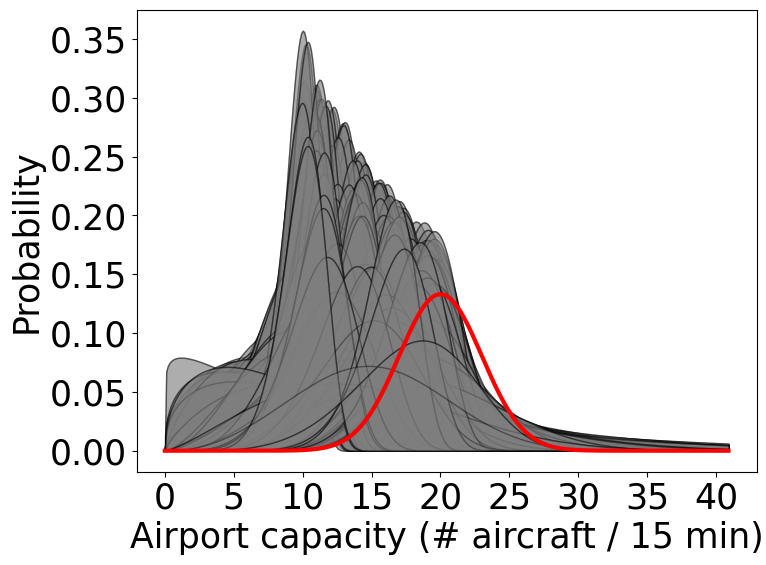}}
\caption{Example of distributions included in the Wasserstein ambiguity set for the empirical Gaussian distribution $\mathcal{N}(20,3)$ with $\epsilon=0.005$.}\label{fig:wsdistributions}
\end{figure}

Recall that we seek to make optimal ground delay allocation decisions that are \emph{robust} with respect to the predicted airport capacity distribution (output of Section \ref{ssec:airport_cap_dis_pred}). To do so, we will solve the MAGHP under a \enquote{worst-case} capacity distribution, located within the Wasserstein ambiguity set centered at the predicted airport capacity distribution. 
Figure~\ref{fig:wsdistributions} presents a numerical example of various distributions—and distribution families—that lie within an ambiguity set of radius $\epsilon = 0.005$, centered at a Gaussian distribution $\mathcal{N}(20,3)$. The parameters are chosen purely for illustrative purposes and carry no operational significance.
Formally, let $\mathcal{M}(\Xi)$ be the space of all probability distributions $\mathbb{Q}$ with support $\Xi$. The Wasserstein distance $d_w : \mathcal{M}(\Xi) \times \mathcal{M}(\Xi) \rightarrow \mathbb{R}_{\geq 0}$ is the minimum transportation cost between distributions $\mathbb{Q}_1 \in \mathcal{M}(\Xi)$ and $\mathbb{Q}_2 \in \mathcal{M}(\Xi)$ \cite{DR2}, and is given explicitly as: 

\vspace{-0.3cm}
\begin{equation}
\begin{aligned}
d_{w}\left(\mathbb{Q}_1,\mathbb{Q}_2\right) = \inf_{\Pi \in \mathcal{D}_{\Pi} \left(\xi_1,\xi_2 \right)}\int_{\Xi^2} \norm{\xi_1 - \xi_2}_2 \, \Pi(d\xi_1,d\xi_2),
\end{aligned}
\label{eq:ws_1.0}
\end{equation}

\noindent
where $\Pi$ is a joint distribution of random variables $\xi_1$ and $\xi_2$ with marginals $\mathbb{Q}_1$ and $\mathbb{Q}_2$, respectively. 
We denote $\mathcal{D}_{\Pi}\left(\xi_1,\xi_2\right)$ as the set of all joint distributions on $\xi_1$ and $\xi_2$ with marginals $\mathbb{Q}_1$ and $\mathbb{Q}_2$. 

Let $Z$ be the set of all airports. For airport $z$, we let $M_{z}$ be the maximum capacity of airport $z$ and $\left\{ \widehat{\xi}_1,\widehat{\xi}_2, \hdots ,\widehat{\xi}_{M_{z}} \right\}$ be the set of airport capacities for airport $z$ with the corresponding estimated probabilities of occurrence $\left\{ \widehat{p}_1,\widehat{p}_2, \hdots ,\widehat{p}_{M_{z}} \right\}$ (i.e., this is precisely the predicted capacity distribution's probability mass function, or PMF, which is the output of the prediction step described in Section \ref{ssec:airport_cap_dis_pred}). The ambiguity set centered around a predicted capacity distribution $\widehat{P}$ with radius $\epsilon > 0$, denoted as $\mathcal{P}_{\epsilon}\left( \widehat{P} \right)$, is given by
\begin{equation}
\begin{aligned}
\mathcal{P}_\epsilon\left(\widehat{P}\right) \coloneqq \left\{ \mathbb{Q} \in \mathcal{M}(\Xi) : d_{w}\left(\widehat{P},\mathbb{Q} \right) \leq \epsilon  \right\}.
\end{aligned}
\label{eq:ws_1.1}
\end{equation}
Note that to differentiate between arrival and departure capacities, we will use subscripts and superscripts with $a$ and $g$, respectively, when we write out the full \textsc{dr-MAGHP} model. We refer readers to \cite{drGHP} for additional technical details on specifying the ambiguity set.

We denote the set of flights by $F$, $F^{(g)}(z), F^{(a)}(z)$ as the set of flights departing from and landing at airport $z$, 
$D_z(t)$ and $R_z(t)$ as the departure and arrival capacity of airport $z \in Z$ at time $t$, respectively, and $\mathcal{C} = F \times F$ as the set of all flight pairs connected by the same aircraft (or tail), where $(f, f') \in \mathcal{C}$ denotes the preceding flight $f$ and successive flight $f'$. 
The decision variables $u_{f,t}$ and $v_{f,t}$ are binary variables where $u_{f,t}$ equals one if flight $f$ will departure at time $t$; similarly, $v_{f,t}$ equals one if flight $f$ will land at time $t$. $z^{g}_{f}$ and $z^{a}_{f}$ represent flights scheduled to depart or land at airport $z$, respectively. $d_f$ is the scheduled departure time of $f$ and $r_f$ is the scheduled arrival time of $f$. $T^{f}_{d}$ is the set of available time periods for $f$ to take off and $T^{f}_{a}$ is the set of available time periods for $f$ to land. $g_{f},a_{f}$ are ground holding delay and airborne delay, respectively, where $g_{f} = \sum_{t\in T^{f}_{d}}tu_{f,t} - d_f$ and $a_{f} = \sum_{t \in T^{f}_{a}}tv_{f,t} - r_{f} - g_{f}$. These definitions follow the standard formulation of the MAGHP \cite{GHP1} and we denote the deterministic MAGHP as \textsc{det-MAGHP}.
The formulation is given explicitly below:
\vspace{-0.3cm}
\begin{subequations}
\begin{alignat}{2}
\min_{u,v} \quad \sum_{f=1}^{F}&\left(C^g_{f}g_f + C^a_{f}a_f\right) \label{eq:det_obj}\\
\textrm{s. t.}\quad \sum_{f:z^{g}_{f} = z}u_{f,t} &\leq D_z(t), && \quad \forall z \in Z, t\in T, \label{eq:det_1b} \\
\sum_{f:z^{a}_{f} = z}v_{f,t} &\leq R_z(t), && \quad \forall z \in Z, t\in T, \label{eq:det_1c} \\
\sum_{t \in T_f^{d}}u_{f,t} &= 1, && \quad \forall f \in F, \label{eq:det_1d}\\
\sum_{t \in T_f^{a}}v_{f,t} &= 1, && \quad \forall f \in F, \label{eq:det_1e}\\
g_{f'} + a_{f'} - s_{f'} &\leq g_{f}, && \quad \forall (f, f') \in \mathcal{C}, \label{eq:det_1f}\\
a_{f},g_{f}&\geq 0, && \quad \forall f \in F, \label{eq:det_1g}\\
u_{f,t},v_{f,t} &\in \{0,1\}, && \quad \forall f \in F, t \in T\label{eq:det_1h}.
\end{alignat}
\label{eq:det_MAGHP}
\end{subequations} 
The objective function in Equation \eqref{eq:det_obj} is the sum of ground holding delay and airborne delay across all flights. Constraints \eqref{eq:det_1b} and \eqref{eq:det_1c} are airport departure and arrival capacity constraints, respectively. Constraints \eqref{eq:det_1d} and \eqref{eq:det_1e} ensure only one departure and arrival time slot is assigned to each flight, respectively,
and \eqref{eq:det_1f} enforces minimum ground turnaround times for connecting flights $(f, f') \in \mathcal{C}$.
We incorporate predicted airport departure and arrival capacity distributions from Section \ref{ssec:airport_cap_dis_pred}, along with robustness guarantees, through a two-stage formulation. 
In the two-stage setting, in addition to the first stage decision variables $u_{f,t}$ and $v_{f,t}$, we introduce second stage decision variables $y^{(g)}_{z,t}$ and $y^{(a)}_{z,t}$.
Decision variables \(y^{(g)}_{z,t}\) denote the additional number of flights entering the departure queue at time period \(t\), while \(y^{(a)}_{z,t}\) represents the additional number of flights entering the arrival queue at the same time period. The second stage of the model focuses on rescheduling the actual departure and arrival times of flights once the actual airport capacity distributions are realized. The objective of the second stage model is to minimize the number of flights joining either queues. We first formulate the two-stage stochastic MAGHP (\textsc{sp-MAGHP}) as follows:
\begin{subequations}
\begin{alignat}{2}
\min_{u,v,g,a,s} \quad \left\{ \sum_{f \in F} (C_g g_f +\right.&\left. C_a a_f )
+ \sum_{z\in Z}\mathbb{E}_{\widehat{P}_{z}^{(g)}}\left[Q_{z}^{(g)} \left(u,\xi_{z}^{(g)}\right)\right] \right. \notag \\ \left.
+ \right. & \left.\sum_{z\in Z}\mathbb{E}_{\widehat{P}_{z}^{(a)}}\left[Q_{z}^{(a)}\left(v,\xi_{z}^{(a)} \right)\right]\right\} &&\label{eq:twostage_objfunc}\\
\text{s.t.} \quad 
&\sum_{t \in T}u_{f,t} = 1, &&\qquad\forall f \in F, \label{eq:twostage_1b}\\
&\sum_{t \in T}v_{f,t} = 1, &&\qquad\forall f \in F, \label{eq:twostage_1c}\\
&g_{f'} + a_{f'} - s_{f'} \leq a_{f}, &&\qquad\forall (f, f') \in \mathcal{C}, \label{eq:twostage_1d} \\
&u_{f,t} \in \{0,1\}, &&\qquad\forall f \in F, t \in T, \label{eq:twostage_1e} \\ 
&g_f, a_f, s_f \geq 0, &&\qquad\forall f \in F, \label{eq:twostage_1f}
\end{alignat}
\label{eq:twostage_a}
\end{subequations}
where $Q^{(g)}\left(u,\xi^{(g)}\right)$ is:
\begin{subequations}
\begin{alignat}{2}
\min_{y^{(g)}} \quad & \sum_{t \in T} C_{g}\left(\xi_{z}^{(g)}\right) y^{(g)}_{z,t}\left(\xi_{z}^{(g)}\right) \label{eq:twostage_21a}\\
\textrm{s. t.}\sum_{f:z^{g}_{f} = z} u_{f,t} &\leq \xi_{z,t}^{(g)} + y^{(g)}_{z,t} \left(\xi_{z}^{(g)} \right), \label{eq:twostage_21b}\\ &\qquad\qquad\;\forall z 
\in Z,t \in T,\xi_{z}^{(g)} \in \Xi_{z}^{(g)} , \nonumber \\
y^{(g)}_{z,t}\left(\xi_{z}^{(g)}\right) &\geq 0, \qquad \forall z 
\in Z,t \in T,\xi_{z}^{(g)} \in \Xi_{z}^{(g)}, \label{eq:twostage_21c}
\end{alignat}
\label{eq:twostage_2a}
\end{subequations}
and $Q^{(a)}\left(u,\xi^{(a)}\right)$ shares a similar formulation:
\begin{subequations}
\begin{alignat}{2}
\min_{y^{(a)}} \quad & \sum_{t \in T} C_{a}\left(\xi_{z}^{(a)}\right)y^{(a)}_{z,t}\left(\xi_{z}^{(a)}\right) \label{eq:twostage_22a}\\
\textrm{s. t.}\sum_{f:z^{a}_{f} = z} v_{f,t} &\leq \xi_{z,t}^{(a)} + y^{(a)}_{z,t} \left(\xi_{z}^{(a)} \right), \label{eq:twostage_22b}\\ &\qquad\qquad\;\forall z \in Z, t \in T, \xi_{z}^{(a)} \in \Xi_{z}^{(a)}, \nonumber \\
y^{(a)}_{z,t}\left(\xi_{z}^{(a)}\right) &\geq 0, \qquad\forall z \in Z, t \in T, \xi_{z}^{(a)} \in \Xi_{z}^{(a)}. \label{eq:twostage_22c}
\end{alignat}
\label{eq:twostage_2b}
\end{subequations}
We denote $\widehat{P}_{z}^{(g)}$ and $\widehat{P}_{z}^{(a)}$ as the predicted departure and arrival capacity distributions, respectively, obtained from the upstream model described in Section~\ref{ssec:airport_cap_dis_pred}.  
For each airport $z$, let $\Xi_{z}^{(g)}$ and $\Xi_{z}^{(a)}$ represent the sets of possible supports (scenarios) for the departure and arrival distributions.  
We further denote $\xi_{z}^{(g)}$ and $\xi_{z}^{(a)}$ as specific scenarios within these sets.  
Finally, $y^{(g)}_{z,t}$ and $y^{(a)}_{z,t}$ correspond to the second-stage decision variables under each scenario at time $t$. The objective function \eqref{eq:twostage_objfunc} includes the first-stage costs plus $\sum_{z\in Z}\mathbb{E}_{\widehat{P}_{z}^{(g)}}\left[Q_{z}^{(g)} \left(u,\xi_{z}^{(g)}\right)\right]$ and $\sum_{z\in Z}\mathbb{E}_{\widehat{P}_{z}^{(a)}}\left[Q_{z}^{(a)}\left(v,\xi_{z}^{(a)} \right)\right]$. Constraints \eqref{eq:twostage_1b}-\eqref{eq:twostage_1d} are the first stage assignment and coupling constraints inherited from the standard, deterministic MAGHP. The first constraint in the second stage minimization problem ensures that, even if the airport capacity is reduced, the number of departing (or arriving) flights at time $t$ does not exceed the realized airport capacity, plus the total number of extra flights allowed to depart (or arrive) at time $t$, across all airports. When there is a drop in airport capacity, the two-stage model will optimally adjust delay allocations based on the weighting between airborne and ground delays (typically airborne delays are 1.2 to 3 times more expensive \cite{delay_cost}). We also assume in this formulation that the capacity distributions across airports are mutually independent. Additionally, recall that we assume that the departure and arrival capacity distributions at each airport are considered independent of one another.

With the formulation of \textsc{sp-MAGHP} and the assumptions we have mentioned above, we then develop the distributionally robust MAGHP model (\textsc{dr-MAGHP}). For each airport $z$, we construct the Wasserstein ambiguity sets around its departure and arrival capacities, and the \textsc{dr-MAGHP} can be written as follows:
\begin{subequations}
\begin{alignat}{2}
\min_{u,v} \quad \left\{\sum_{f \in F} \left(C_g g_f + C_a a_f \right) \right.&\left. + \sum_{z\in Z} \max_{p \in \mathcal{P}_{\epsilon}  \left(\widehat{P}_{z}^{\left(g\right)}\right)} \, \mathbb{E}_{p}\left[Q^{(g)} \left(u,\xi_{z}^{(g)}\right)\right] \notag \right. \\ & \left.+ \sum_{z \in Z}\max_{p \in \mathcal{P}_{\epsilon}\left(\widehat{P}_{z}^{\left(a\right)}\right)} \, \mathbb{E}_{p}\left[Q^{(a)}\left(v,\xi_{z}^{(a)} \right)\right] \right\}  \label{eq:dr_objfunc}\\
\textrm{s. t.}\quad\sum_{t \in T^{f}_{d}}u_{f,t} &= 1, \forall f \in F, \label{eq:dr_twostage_1b}\\
\sum_{t \in T^{f}_{a}}v_{f,t} &= 1, \forall f \in F, \label{eq:dr_twostage_1c}\\
g_{f'} + a_{f'} - s_{f'} &\leq a_{f}, \forall (f, f') \in \mathcal{C}.\label{eq:dr_twostage_1d}
\end{alignat}
\label{eq:dr_a}
\end{subequations}
The objective function \eqref{eq:dr_objfunc} includes two inner maximization problems, which seek the worst-case distribution within each ambiguity set that maximizes the expected second stage cost. Constraints \eqref{eq:dr_twostage_1b}-\eqref{eq:dr_twostage_1d} and $Q^{(g)}, Q^{(a)}$ are the same from the formulation of \textsc{sp-MAGHP} in \eqref{eq:twostage_2b}. $\mathcal{P}_{\epsilon} \left(\widehat{P}_{z}^{\left(g\right)}\right), \mathcal{P}_{\epsilon}  \left(\widehat{P}_{z}^{\left(a\right)}\right)$ are Wasserstein ambiguity sets of size $\epsilon > 0$ constructed based on capacity distribution predictions and $p$ is an arbitrary distribution within each ambiguity set. 

\subsection{Scenario reduction and \textsc{dr-MAGHP} reformulation} \label{ssec:scenario_reduce}
Recall that we optimize ground holding decisions across a prediction horizon of 12 hours, subdivided into 48 time periods (with each unit time period of 15 minutes). From the airport capacity distribution prediction models (Section \ref{ssec:airport_cap_dis_pred}), we are given new predicted distributions of the arrival and departure capacities (and hence, associated Wasserstein ambiguity sets) at each time, across all airports. This leads to an exponential increase in the number of scenarios for the optimization model, resulting in severe computational intractability. The root of this challenge lies in the temporal dependency of decisions: The decision at any given time depends on the sequence of prior decisions. Specifically, if the maximum capacity at an airport within each 15-minute interval is $M$, and the planning horizon contains $|T|$ time periods (discretizations), the total number of scenarios becomes $M^{|T|}$. For instance, at Hartsfield-Jackson Atlanta International Airport (ATL) with a maximum hourly capacity of 35, the scenario space expands to $35^{48}$, rendering a direct solution to be computationally intractable. 

To address this, we apply a scenario reduction strategy aimed at reducing the scenario tree complexity and thus improving model tractability. Our strategy is to reduce discretization points used for capacities and time periods. Firstly, to reduce the number of discretizations $|T|$, we compute the pairwise Wasserstein distance between capacity distributions at consecutive time periods. Time periods with small Wasserstein distances—indicating similar capacity behavior—are grouped together. In contrast, a time step is identified as a change point if its Wasserstein distance (i.e., Equation \eqref{eq:ws_1.0} with $\ell_1$ norm) is among the $n$ largest values across the horizon. These change points divide the full time horizon into $n+1$ intervals. For each interval, a representative capacity distribution is computed as the average (i.e., centroid) of the predicted distributions within that interval. The algorithm for this grouping procedure is given in Algorithm \ref{Alg:wsdcluster}.

{\small
\begin{algorithm}[h]
\caption{Clustering of Capacity Time Series via Similarity Measures}
\begin{algorithmic}[1]
\STATE \textbf{Input}: Predicted PMFs $\left\{\widehat{p}_t^{(z,a)}, \widehat{p}_t^{(z,g)} \right\}_{t \in T}$ for each airport $z \in Z$; number of clusters $n$; Wasserstein distance as the similarity measure $d_w(\cdot, \cdot)$
\FOR{each $z \in Z$}
    \STATE Compute Wasserstein distances: $w_t^{(z,a)} = d_w\left(\widehat{p}_t^{(z,a)}, \widehat{p}_{t-1}^{(z,a)}\right)$, and similarly for departure capacity distributions
    \STATE Identify top $n$ indices $C^{(z,\cdot)} := \arg\max_{S \subset \{1,\dots,|T|-1\},\, |S|=n} \sum_{t \in S} w_t^{(z,\cdot)}$, then sort such that $c_1 < \cdots < c_n$, and $c_1,\cdots,c_n \in C^{(z,\cdot)}$
    \STATE Partition time: $T_1 = [0,c_1],\ T_2 = [c_1+1, c_2],\ \dots,\ T_{n+1} = [c_n+1, |T|-1]$
    \FOR{each interval $T_k^{(z,a)}, k \in \{1,2,\dots,n+1\}$}
        \STATE $\bar{p}_k^{(z,a)} = \frac{1}{|T_k^{(z,a)}|} \sum_{t \in T_k^{(z,a)}} \widehat{p}_t^{(z,a)}$
    \ENDFOR
    \FOR{each interval $T_k^{(z,g)} k \in \{1,2,\dots,n+1\}$}
        \STATE $\bar{p}_k^{(z,g)} = \frac{1}{|T_k^{(z,g)}|} \sum_{t \in T_k^{(z,g)}} \widehat{p}_t^{(z,g)}$
    \ENDFOR
\ENDFOR
\STATE \textbf{Output}: Clustered time intervals $T_k^{(z,\cdot)}$ and representative PMFs $\bar{p}_k^{(z,\cdot)}$ for each airport
\end{algorithmic}
\label{Alg:wsdcluster}
\end{algorithm}
}

To further reduce the complexity of the predicted capacity distributions, we apply $K$-means clustering to the support–probability pairs of each probability mass function (PMF) within the time intervals identified by the Wasserstein-based clustering. For each interval, the original distribution is approximated by a reduced PMF with a fixed number of representative support points. Specifically, the $K$-means algorithm groups the support--probability pairs into $K$ clusters. 
For each cluster, the new support value is computed as the probability-weighted average of the original supports, and the total probability mass is obtained by summing the probabilities of all points in the cluster. This compression scheme preserves the overall shape of the original distribution while significantly reducing its dimensionality, thereby facilitating tractable downstream optimization. The details for compressing the PMFs themselves are laid out in Algorithm \ref{Alg:probmasscluster}.

{\small
\begin{algorithm}[h]
\caption{PMF Clustering via $K$-Means}
\begin{algorithmic}[1]
\STATE \textbf{Input}: For each airport $z\in Z$ and time interval $k$: representative PMF $\left(\bar{p}^{(z,\cdot)}_k, \bar{\xi}^{(z,\cdot)}\right)$ with supports $\bar{\xi}_i \in \mathbb{Z}_{\geq 0}$ for each airport and clustered time interval $k$, and target cluster count $K$
\FOR{each $z \in Z$}
    \FOR{each interval $k$}
        \STATE Form dataset $\mathcal{D}_k^{(z,\cdot)} = \{(\bar{\xi}_i, \bar{p}_i)\}_{i=1}^{N_k}$
        \STATE Apply $K$-means on $\mathcal{D}_k^{(z,\cdot)}$ with $K$ clusters \cite{likas2003global}
        \FOR{each cluster $j = 1,\dots,K$}
            \STATE Let $\mathcal{I}_j$ be indices assigned to cluster $j$
            \STATE Compute cluster weight: $\bar{p}_{k,j}^{(z,\cdot)} = \sum_{i \in \mathcal{I}_j} \bar{p}_i$
            \STATE Compute cluster centroid: $\bar{\xi}_{k,j}^{(z,\cdot)} = \mathrm{round}\left( \frac{\sum_{i \in \mathcal{I}_j} \bar{\xi}_i \bar{p}_i}{\sum_{i \in \mathcal{I}_j} \bar{p}_i} \right)$
        \ENDFOR
        \STATE Define clustered PMF: $\widetilde{p}_k^{(z,\cdot)} = \{ (\bar{\xi}_{k,j}^{(z,\cdot)}, \bar{p}_{k,j}^{(z,\cdot)}) \}_{j=1}^K$
    \ENDFOR
\ENDFOR
\STATE \textbf{Output}: Clustered PMFs $\{\widetilde{p}_k^{(z,\cdot)}\}$ for all $z$ and $k$
\end{algorithmic}
\label{Alg:probmasscluster}
\end{algorithm}
}

An illustration of the reduced scenario tree is shown in Figure~\ref{fig:scenario_tree}, where $\bar{\xi}_{k,j}$ denotes the $j^{\text{th}}$ probability mass center (i.e., representative support) of the $k^{\text{th}}$ time cluster. The scenario tree captures the joint distribution of an airport’s arrival or departure capacities across the planning horizon, with each edge representing a possible support of this joint distribution. For example, the leftmost path in Figure~\ref{fig:scenario_tree}, given by the sequence $\left\{\bar{\xi}_{1,1}^{(z,\cdot)}, \bar{\xi}_{2,1}^{(z,\cdot)}, \bar{\xi}_{3,1}^{(z,\cdot)}\right\}$, corresponds to a scenario in which, at each of the three time clusters, the capacity takes the first representative value of the respective PMF. The probability associated with each scenario in the reduced scenario tree is computed as the product of the marginal probabilities of the selected probability mass clusters at each time stage. Specifically, for a scenario represented by the sequence $\boldsymbol{\bar{\xi}}^{(z,\cdot)} = \left\{\bar{\xi}_{1,j_1}^{(z,\cdot)}, \bar{\xi}_{2,j_2}^{(z,\cdot)}, \dots, \bar{\xi}_{K,j_K}^{(z,\cdot)}\right\}$, where $K$ is the number of time clusters and $\bar{p}_{k,j_k}$ denotes the probability assigned to the $j_k^{\text{th}}$ mass cluster at time cluster $k$, the joint scenario probability is given by:
\[
\widehat{\mathbb{P}}^{(z,\cdot)}(\bar{\xi}_{1,j_1}^{(z)}, \dots, \bar{\xi}_{K,j_K}^{(z)}) = \prod_{k=1}^{K} \bar{p}^{(z,\cdot)}_{k,j_k}.
\]
This formulation assumes conditional independence across time clusters in the reduced scenario construction, such that the joint distribution can be expressed as a product of marginal distributions over time. Therefore, given the reduced scenario tree, we update the formulation of \textsc{sp-MAGHP} as 
\begin{subequations}
\begin{alignat}{2}
\min_{u,v,g,a,s} \quad & \left\{ 
\sum_{f \in F} \left(C_g g_f + C_a a_f \right) 
+ \sum_{z\in Z}\mathbb{E}_{\widehat{\mathbb{P}}^{(z,g)}}\left[Q_{z}^{(g)} \left(u,\boldsymbol{\bar{\xi}}^{(z,g)}\right)\right] 
+ \right. \\ & \left. \sum_{z\in Z}\mathbb{E}_{\widehat{\mathbb{P}}^{(z,a)}}\left[Q_{z}^{(a)}\left(v,\boldsymbol{\bar{\xi}}^{(z,a)}\right)\right] 
\right\} \label{eq:twostage_objfunc_b}\\
\text{s.t.} \quad 
& \text{Constraints } \eqref{eq:twostage_1b} \text{-} \eqref{eq:twostage_1f}, 
\end{alignat}
\label{eq:twostage_b}
\end{subequations}
where both the supports and probabilities are obtained from the joint distribution encoded by the reduced scenario tree, rather than from the original marginal capacity distributions. The second stage value function $Q^{(g)}_{z}\left(u,\boldsymbol{\Bar{\xi}}^{(g)}\right)$ is reformulated as
\begin{subequations}
\label{eq:twostage_2c}
\begin{alignat}{2}
\min_{y^{(g)}}\quad 
& \sum_{t\in T} C_{g}\!\left(\bar{\boldsymbol{\xi}}^{(z,g)}\right)\, y^{(g)}_{t}\!\left(\bar{\boldsymbol{\xi}}^{(z,g)}\right) 
\label{eq:twostage_23a}\\[0.5ex]
\text{s.t.}\quad
& \sum_{f:\, z^{g}_{f}=z} u_{f,t} 
  \;\leq\; D_{t}\!\left(\bar{\boldsymbol{\xi}}^{(z,g)}\right) + y^{(g)}_{t}\!\left(\bar{\boldsymbol{\xi}}^{(z,g)}\right),
  \label{eq:twostage_23b}\\
& \qquad \forall\, t\in T,\ \bar{\boldsymbol{\xi}}^{(z,g)}\in \boldsymbol{\Xi}^{(z,g)}, \nonumber \\[0.5ex]
& y^{(g)}_{t}\!\left(\bar{\boldsymbol{\xi}}^{(z,g)}\right) \;\geq\; 0,
  \label{eq:twostage_23c}\\
& \qquad \forall\, t\in T,\ \bar{\boldsymbol{\xi}}^{(z,g)}\in \boldsymbol{\Xi}^{(z,g)}, \nonumber
\end{alignat}
\end{subequations}
where we denote $D_{t}\left(\boldsymbol{\bar{\xi}}^{(z,g)}\right)$ as the capacity value assigned at time $t$ under scenario $\boldsymbol{\bar{\xi}}^{\left(z,g\right)}$. Then,
$D_{t}\left(\boldsymbol{\bar{\xi}}^{(z,g)}\right) = \Bar{\xi}_{k,j_k} \text{ if } t \in T^{(z,g)}_k$,
where $T^{(z,g)}_k$ is the set of time steps associated with time cluster $k$, and $\boldsymbol{\bar{\xi}}^{(z,g)}= \left(\bar{\xi}^{(z,g)}_{1,j_1}, \bar{\xi}^{(z,g)}_{2,j_2}, \dots, \bar{\xi}^{(z,g)}_{K,j_K}\right)$ is a scenario in the reduced scenario tree. We also note that $Q^{(a)}_{z}\left(u,\boldsymbol{\Bar{\xi}}^{(a)}\right)$ is of a similar form, and for brevity we skip writing down its formulation in its entirety here. The formulation of \textsc{dr-MAGHP} \eqref{eq:dr_a} is updated accordingly as follow:

\begin{subequations}
\begin{alignat}{2}
\min_{u,v} \quad & \sum_{f \in F} \left(C_g g_f + C_a a_f \right) \notag \\
&+ \sum_{z \in Z} \max_{p \in \mathcal{P}_{\epsilon}(\widehat{\mathbb{P}}^{(z,g)})} \, \mathbb{E}_{p}\left[Q^{(g)}\left(u,\boldsymbol{\bar{\xi}}^{(z,g)}\right)\right] \notag \\
&+ \sum_{z \in Z} \max_{p \in \mathcal{P}_{\epsilon}(\widehat{\mathbb{P}}^{(z,a)})} \, \mathbb{E}_{p}\left[Q^{(a)}\left(v,\boldsymbol{\bar{\xi}}^{(z,a)}\right)\right] \label{eq:dr_objfunc_b} \\
\text{s.t.} \quad & \text{Constraints } \eqref{eq:dr_twostage_1b} \text{--} \eqref{eq:dr_twostage_1d}.
\end{alignat}
\label{eq:dr_b}
\end{subequations}
The Wasserstein ambiguity sets appearing in the objective function is explicitly expressed as follows
\begin{equation}
\begin{aligned}
\mathcal{P}_{\epsilon} \left(\widehat{\mathbb{P}}^{(z,\cdot)}\right) \coloneqq \left\{ \mathbb{Q} \in M\left(\boldsymbol{\Xi}^{(z,\cdot)}\right) : d_{w}\left(\widehat{\mathbb{P}}^{(z,\cdot)}, \mathbb{Q} \right) \leq \epsilon \right\}.
\end{aligned}
\label{eq:ws_2}
\end{equation}
However, the formulation in~\eqref{eq:dr_b} is generally computationally intractable, primarily due to the integral expression in~\eqref{eq:ws_1.0}, which leads to an infinite-dimensional optimization problem. Furthermore, the inherent min–max structure introduces additional complexity, rendering the problem even more challenging to solve directly. To address these challenges, we reformulate the \textsc{dr-MAGHP} in~\eqref{eq:dr_b} by transforming the inner second-stage maximization problems into equivalent minimization problems, resulting in a semi-infinite program. We then apply discretization techniques to handle the continuous support of this semi-infinite formulation, ultimately yielding a deterministic equivalent representation of the \textsc{dr-MAGHP}. For a comprehensive derivation of the deterministic reformulation, we refer the reader to~\cite{drGHP}. In this paper, we directly present the deterministic equivalent formulation of \textsc{dr-MAGHP} as follows:
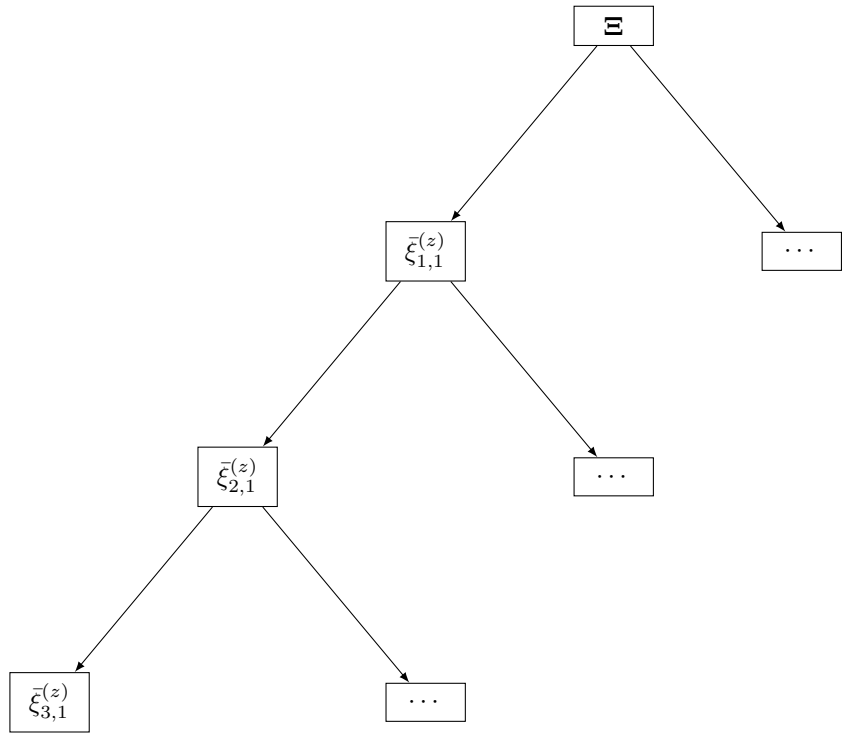
\begin{figure}
\centering
\begin{tikzpicture}[
  scale=2.5,
  grow=down,
  level distance=1.2cm,
  sibling distance=2cm,
  edge from parent/.style={draw, -latex},
  every node/.style={
    draw, 
    rectangle, 
    font=\footnotesize, 
    minimum width=2.5em,
    minimum height=1.2em,
    align=center
  }
]

\node {$\boldsymbol{\Xi}$}
  child {node {$\bar{\xi}_{1,1}^{(z)}$}
    child {node {$\bar{\xi}_{2,1}^{(z)}$}
      child {node {$\bar{\xi}_{3,1}^{(z)}$}}
      child {node {$\cdots$} [draw=none]} 
    }
    child {node {$\cdots$} [draw=none]} 
  }
  child {node {$\cdots$} [draw=none]}; 

\end{tikzpicture}
\caption{Illustration of the reduced scenario tree for capacity distributions. Each node $\Bar{\xi}_{k,j}$ represents the $j^{\text{th}}$ support value at stage $k$.}
\label{fig:scenario_tree}
\end{figure}

 \begin{subequations}
     \begin{alignat}{2}
         \min_{\mathbf{\alpha},\mathbf{\beta},\textbf{y},u,v} \left.\biggl\{\sum_{f \in F} \Bigg(C_gg_f + C_aa_f\Bigg) \right. & \left. + \Phi^{(g)} + \Phi^{(a)} \right.\biggl\} \label{eq:drmaghp_objfunc} \\
         \mathrm{s. t.}\quad\alpha^{(z,g)} \norm{\boldsymbol{\Bar{\xi}}^{(z,g)} - \boldsymbol{\xi}^{(z,g)}}_{2} +\beta^{(z,g)}\left(\boldsymbol{\bar{\xi}}^{(z,g)}\right) &\geq \sum_{t \in T} C_g\left(\boldsymbol{\bar{\xi}}^{(z,g)}\right) y^{(g)}_{t}\left(\boldsymbol{\bar{\xi}}^{(z,g)}\right), \notag \\ \forall \boldsymbol{\bar{\xi}}^{(z,g)} \in \boldsymbol{\Xi}^{(z,g)}\,\forall z \in Z, \, & \boldsymbol{\xi}^{(z,g)} \in \boldsymbol{\Xi}^{(z,g)}, \label{eq:drmaghp_1} \\
         \quad \alpha^{(z,a)} \norm{\Bar{\boldsymbol{\xi}}^{(z,a)} - \boldsymbol{\xi}^{(z,a)}}_{2} +\beta^{(z,a)}\left(\boldsymbol{\Bar{\xi}}^{(z,a)}\right) &\geq \sum_{t \in T} C_a\left(\boldsymbol{\bar{\xi}}^{(z,a)}\right) y^{(a)}_{t}\left(\boldsymbol{\bar{\xi}}^{(z,a)}\right), \notag \\ \forall \boldsymbol{\bar{\xi}}^{(z,a)} \in \boldsymbol{\Xi}^{(z,a)},\,\forall z \in Z, \, & \boldsymbol{\xi}^{(z,a)} \in \boldsymbol{\Xi}^{(z,a)}, \label{eq:drmaghp_2} \\
         \sum_{t \in T} u_{f,t} &= 1, \notag \\ \forall f \in F,\label{eq:drmaghp_3}\\
         \sum_{t \in T} v_{f,t} &= 1, \notag \\ \forall f \in F,\label{eq:drmaghp_4}\\
          D_{t}\left(\boldsymbol{\xi}^{(z,g)}\right) + y^{(g)}_{t}\left(\boldsymbol{\xi}^{(z,g)}\right)& \geq \sum_{f \in F^{(g)}(z)}u_{f,t} , \notag \\ \forall t \in T \,  & \boldsymbol{\xi}^{(z,g)} \in \boldsymbol{\Xi}^{(z,g)}, \label{eq:drmaghp_5}\\
         R_{t}\left(\boldsymbol{\xi}^{(z,a)}\right) + y^{(a)}_{t}\left(\boldsymbol{\xi}^{(z,a)}\right) & \geq \sum_{f \in F^{(a)}(z)}v_{f,t}, \notag \\ \forall t \in T, \,  & \boldsymbol{\xi}^{(z,a)} \in \boldsymbol{\Xi}^{(z,a)}, \label{eq:drmaghp_6}\\
         y^{(\cdot)}_{0}\left(\boldsymbol{\bar{\xi}}^{(z,\cdot)}\right) = 0, \, y^{(\cdot)}_{t}\left(\boldsymbol{\bar{\xi}}^{(z,\cdot)}\right)&\geq 0, \, \alpha^{(z,\cdot)} \geq 0,\notag \\ \forall t \in T, \, \boldsymbol{\xi}^{(z,\cdot)} \in & \boldsymbol{\Xi}^{(z,\cdot)} \, z \in Z. \label{eq:drmaghp_7}
     \end{alignat}
     \label{eq:drmaghp}
 \end{subequations}
 In the objective function \eqref{eq:drmaghp_objfunc}, $\Phi^{(g)}$ and $\Phi^{(a)}$ are substitutions for the dual objective function; explicitly, we have that:
 \begin{equation}
     \begin{aligned}
         \Phi^{(g)} = \epsilon^{(g)} \sum_{z \in Z}\alpha^{(z,g)} \notag + \sum_{z \in Z}\sum_{\boldsymbol{\bar{\xi}}^{(z,g)} \in \boldsymbol{\Xi}^{(z,g)}} \widehat{\mathbb{P}}^{(z,g)}(\boldsymbol{\bar{\xi}}^{(z,g)})\beta^{(z,g)}\left(\boldsymbol{\bar{\xi}}^{(z,g)}\right),\\
         \Phi^{(a)} = \epsilon^{(a)} \sum_{z \in Z}\alpha^{(z,a)} \notag + \sum_{z \in Z}\sum_{\boldsymbol{\bar{\xi}}^{(z,a)} \in \boldsymbol{\Xi}^{(z,a)}} \widehat{\mathbb{P}}^{(z,a)}(\boldsymbol{\bar{\xi}}^{(z,a)})\beta^{(z,a)}\left(\boldsymbol{\bar{\xi}}^{(z,a)}\right).
     \end{aligned}
 \end{equation}
 $\mathbf{\alpha}$ and $\mathbf{\beta}$ are dual variables for the two-stage dr-MAGHP, $\epsilon$ is the radius for the Wasserstein ambiguity set, and $\widehat{\mathbb{P}}^{(z,\cdot)}(\boldsymbol{\bar{\xi}}^{(z,\cdot)})$ is the probability assigned to each scenario within the empirical distribution. Constraints \eqref{eq:drmaghp_1} and \eqref{eq:drmaghp_2} are derived from the Wasserstein distance constraints of \eqref{eq:twostage_b}. \eqref{eq:drmaghp_3} and \eqref{eq:drmaghp_4} are the first stage assignment constraints, and \eqref{eq:drmaghp_5} and \eqref{eq:drmaghp_6} are the capacity constraints for each scenario of each distribution within the ambiguity set. 
 
\section{Numerical Experiments and Discussion} \label{sec:Experiments}

In the following numerical experiments, we evaluate the proposed integrated prediction and DRO framework as a whole. The upstream capacity prediction model generates probabilistic predictions, which are then used as input distributions for the downstream \textsc{sp-MAGHP} and \textsc{dr-MAGHP}. We first examine the quality of the predictive distributions and then assess how these predictions translate into prescriptive actions (i.e., air traffic management strategic rescheduling) through the optimization stage.


\subsection{Data Description}
We obtained airport throughput data from the US Department of Transportation's Bureau of Transportation Statistics (BTS), and weather data from the US National Oceanic and Atmospheric Administration's High-Resolution Rapid Refresh (HRRR) database. BTS provides detailed information for each flight, including the scheduled departure and arrival times, actual departure and arrival times, delay duration.
Using these data points and the procedure described in Section \ref{sssec:derive_actual_capacities}, we estimate the capacity of each FAA Core 30 airport in our study (note that we can scale to a larger network of airports if needed).
HRRR provides weather data on a $3$ km $\times$ $3$ km grid covering all 50 US states, with a forecast horizon of up to 23 hours from the current hour. We collect data for the entire year of 2019. Each day is divided into 96 quarter-hour intervals, resulting in 35,040 time periods in total.
\subsection{Prediction Model Setup and Evaluation}
\subsubsection{Experiment setup}

We use data comprising 60 datasets in total from the 30 US airports, with each airport providing 2 datasets for arrival and departure operations. Each dataset contains approximately 13\%-48\%  of time periods remaining after rule-based capacity estimation presented in \eqref{eq:capacity_estimation}, with an overall average of 25.7\% for arrivals and 26.7\% for departures. 
As a reminder, the rule-based approach identifies capacity-saturated periods when any one of three criteria is satisfied: (1) throughput-based criterion (actual throughput $\geq$ 90th percentile threshold), (2) demand-based criterion (throughput-demand ratio $\leq$ 0.8, indicating unmet demand), and (3) delay-based criterion (average delay $\geq$ 15 minutes with at least 2 delayed flights). The exact number of remaining time periods for each airport across both arrival and departure operations is detailed in~\ref{appendix:rule_estimation_table}.

For the train-validation-test split, we follow temporal order by using the first 10 weeks of each quarter for training, the 11\textsuperscript{th} week for validation, and the final 12\textsuperscript{th} week for testing. We normalize all numerical features using min-max normalization applied to the training set to prevent the model from being dominated by a subset of variables, with the same normalization parameters applied to both validation and test sets to ensure consistent scaling across all data splits. Hyperparameter tuning is performed via grid search, resulting in a learning rate of 0.0001, 300 epochs, and a batch size of 16.


\subsubsection{Evaluation Metrics}
We use four metrics to evaluate predictive model performance: \emph{Root Mean Squared Error (RMSE)}, \emph{Mean Absolute Error (MAE)}, \emph{Prediction Interval Coverage Probability (PICP)}, and \emph{Mean Prediction Interval Width (MPIW)}. These metrics provide a comprehensive assessment of both point prediction accuracy and distributional prediction quality.

\textbf{Point prediction results:} The model predicts the capacity probability distribution for each airport $i$ at time $t$. For point prediction evaluation, we select the capacity value with the highest predicted probability:
\begin{equation}
\hat{c}_{i^*,t} = \underset{i}{\mathrm{argmax}} \, p_{i,t},
\end{equation}
where $p_{i,t}$ represents the predicted probability for capacity value $i$ at time $t$ for a specific airport, and $\hat{c}_{i^*,t}$ is the predicted capacity value with the highest probability. The RMSE measures the prediction accuracy in terms of squared errors:
\begin{equation}
\text{RMSE} = \sqrt{\frac{1}{N_t} \sum_{t=1}^{N_t} (\hat{c}_{i^*,t} - c_t)^2},
\end{equation}
and the MAE provides a linear penalty for prediction errors:
\begin{equation}
\text{MAE} = \frac{1}{N_t} \sum_{t=1}^{N_t} |\hat{c}_{i^*,t} - c_t|,
\end{equation}
where $c_t$ is the actual observed capacity value at time $t$.

\textbf{Distributional prediction quality:}
For evaluating the quality of predicted capacity distributions, we employ prediction interval-based metrics. The tolerance rate $TR_t$ represents the prediction interval at time $t$ under a 90\% confidence level, defined as the set of capacity values whose cumulative probability reaches or exceeds 90\%:
\begin{equation}
\sum_{i \in TR_t} p_{i,t} \geq 0.9.
\end{equation}

The PICP measures how often the actual capacity falls within the predicted tolerance rate:
\begin{equation}
\text{PICP} = \frac{1}{N_t} \sum_{t=1}^{N_t} I(c_t \in TR_t),
\end{equation}
where $I(\cdot)$ is the indicator function that equals 1 if the condition is true, and 0 otherwise. The MPIW quantifies the average width of the prediction intervals:
\begin{equation}
\text{MPIW} = \frac{1}{N_t} \sum_{t=1}^{N_t} |TR_t|,
\end{equation}
where $|TR_t|$ represents the range or size of the tolerance rate.

The prediction model is optimized to minimize RMSE and MAE for point prediction accuracy, while simultaneously maximizing PICP and minimizing MPIW for distributional prediction quality. A well-calibrated model should achieve high PICP (ideally close to 0.9, corresponding to the 90\% confidence level) while maintaining narrow prediction intervals (low MPIW).

\subsubsection{Capacity Distribution Prediction Results}

The capacity distribution prediction performance across the core 30 US airports, as shown in Tables~\ref{tab:prediction_results_arrival} and \ref{tab:prediction_results_departure}, reveals significant variability in model accuracy and calibration quality, reflecting the diverse operational characteristics and complexity of different airport environments. The RMSE and MAE results demonstrate considerable variation in model performance across airports. The model achieves relatively strong prediction accuracy at medium or relatively smaller large hub airports such as MEM (arrival RMSE: 0.81, MAE: 0.45), MDW (arrival RMSE: 1.59, MAE: 1.11), and BWI (arrival RMSE: 1.73, MAE: 1.23). In contrast, major hub airports present substantial prediction challenges, with LAX showing the highest arrival prediction errors (RMSE: 8.55, MAE: 6.84), followed by DEN (RMSE: 7.93, MAE: 6.09) and CLT (RMSE: 6.99, MAE: 5.50). This pattern suggests that airports with complex operational environments, multiple runway configurations, and high traffic variability inherently pose greater difficulties for capacity prediction models.

Additionally, departure capacity prediction exhibits a notably different performance pattern compared to arrivals. Several airports show significantly higher departure prediction errors, with DTW demonstrating particularly poor departure performance (RMSE: 10.04, MAE: 8.00) and IAH showing substantial challenges for both arrivals (RMSE: 5.66) and departures (RMSE: 10.41). This asymmetry between arrival and departure prediction accuracy likely reflects the greater complexity of departure operations, which are more susceptible to ground delays, pushback sequencing constraints, airline scheduling decisions, and air traffic control departure sequencing.

The PICP results reveal systematic calibration issues across the majority of airports, with most facilities achieving coverage probabilities substantially below the target 90\% confidence level. The best-performing airports in terms of distributional prediction include DTW (arrival PICP: 77.89\%, departure PICP: 79.34\%) and IAD (departure PICP: 80.62\%), though even these fall short of the ideal 90\% coverage. Most airports exhibit PICP values in the 45-75\% range, indicating significant underconfidence in the predicted distributions. The consistent underestimation of prediction intervals suggests that the current approach may not adequately capture all sources of capacity variability, including weather-related variations, operational disruptions, and seasonal traffic patterns. 
The MPIW results show that prediction interval width generally correlates with airport complexity and size. Major hub airports like DTW (departure MPIW: 19.87), CLT (departure MPIW: 17.71), and DEN (arrival MPIW: 14.13) exhibit the widest prediction intervals, reflecting the inherent uncertainty in predicting capacity at complex operational environments. Conversely, relatively smaller airports like MDW (arrival MPIW: 3.01) and TPA (arrival MPIW: 3.08) demonstrate narrower intervals, suggesting more predictable capacity patterns. The inherent prediction uncertainties observed across these diverse airport environments motivate the need for distributionally robust optimization procedures to effectively counter these prediction challenges in subsequent capacity allocation decisions.

\begin{table}[htbp]
\centering
\begin{threeparttable}
\resizebox{0.7\textwidth}{!}{
\begin{tabular}{c|c|cccc}
\hline
\textbf{Index} & \textbf{Airport} & \textbf{RMSE} ($\downarrow$) & \textbf{MAE} ($\downarrow$) & \textbf{PICP (\%)} ($\uparrow$) & \textbf{MPIW} ($\downarrow$) \\
\hline
1  & ATL & 6.81 & 5.24 & 60.42 & 10.82 \\
2  & BOS & 2.36 & 1.79 & 65.95 & 4.43 \\
3  & BWI & 1.73 & 1.23 & 67.35 & 3.07 \\
4  & CLT & 6.99 & 5.50 & 45.79 & 8.32 \\
5  & DCA & 2.30 & 1.69 & 58.41 & 3.42 \\
6  & DEN & 7.93 & 6.09 & 71.77 & 14.13 \\
7  & DFW & 5.86 & 4.50 & 54.02 & 7.76 \\
8  & DTW & 6.42 & 5.19 & 77.89 & 14.48 \\
9  & EWR & 2.39 & 1.87 & 66.43 & 4.52 \\
10 & FLL & 2.03 & 1.56 & 62.82 & 3.65 \\
11 & HNL & 1.37 & 0.85 & 76.82 & 3.01 \\
12 & IAD & 3.68 & 2.76 & 70.56 & 7.00 \\
13 & IAH & 5.66 & 4.45 & 53.05 & 8.79 \\
14 & JFK & 2.29 & 1.72 & 64.09 & 4.04 \\
15 & LAS & 2.62 & 2.00 & 50.80 & 3.41 \\
16 & LAX & 8.55 & 6.84 & 48.82 & 11.35 \\
17 & LGA & 2.63 & 2.03 & 58.14 & 3.97 \\
18 & MCO & 2.34 & 1.82 & 63.40 & 4.32 \\
19 & MDW & 1.59 & 1.11 & 66.16 & 3.01 \\
20 & MEM & 0.81 & 0.45 & 92.69 & 3.03 \\
21 & MIA & 2.50 & 1.82 & 59.43 & 3.71 \\
22 & MSP & 3.96 & 3.05 & 49.63 & 5.23 \\
23 & ORD & 6.72 & 5.11 & 58.21 & 10.03 \\
24 & PHL & 3.46 & 2.70 & 66.01 & 6.19 \\
25 & PHX & 3.60 & 2.73 & 45.83 & 4.27 \\
26 & SAN & 1.71 & 1.22 & 68.52 & 3.31 \\
27 & SEA & 3.19 & 2.45 & 64.38 & 5.37 \\
28 & SFO & 2.90 & 2.31 & 61.67 & 4.90 \\
29 & SLC & 3.84 & 3.01 & 46.26 & 4.74 \\
30 & TPA & 1.48 & 1.03 & 72.02 & 3.08 \\
\hline
\end{tabular}
}
\begin{tablenotes}
\footnotesize
\item Note: Arrows indicate preferred direction: ($\downarrow$) lower is better, ($\uparrow$) higher is better. Airports follow the standard IATA codes. A complete list of airport names corresponding to IATA codes is provided in~\ref{appendix:airport_codes}.
\end{tablenotes}
\caption{Arrival capacity distribution prediction performance for the core 30 US airports in the test set.}
\label{tab:prediction_results_arrival}
\end{threeparttable}
\end{table}

\begin{table}[htbp]
\centering
\begin{threeparttable}
\resizebox{0.7\textwidth}{!}{
\begin{tabular}{c|c|cccc}
\hline
\textbf{Index} & \textbf{Airport} & \textbf{RMSE} ($\downarrow$) & \textbf{MAE} ($\downarrow$) & \textbf{PICP (\%)} ($\uparrow$) & \textbf{MPIW} ($\downarrow$) \\
\hline
1  & ATL & 8.07 & 6.48 & 57.74 & 13.04 \\
2  & BOS & 3.20 & 2.47 & 65.57 & 6.02 \\
3  & BWI & 2.19 & 1.64 & 55.10 & 3.18 \\
4  & CLT & 8.97 & 7.12 & 79.14 & 17.71 \\
5  & DCA & 3.38 & 2.66 & 55.04 & 5.04 \\
6  & DEN & 9.19 & 6.73 & 64.77 & 12.93 \\
7  & DFW & 9.31 & 7.11 & 62.95 & 14.47 \\
8  & DTW & 10.04 & 8.00 & 79.34 & 19.87 \\
9  & EWR & 3.24 & 2.56 & 65.47 & 5.80 \\
10 & FLL & 2.00 & 1.57 & 65.54 & 3.89 \\
11 & HNL & 1.12 & 0.61 & 87.57 & 3.00 \\
12 & IAD & 5.62 & 4.38 & 80.62 & 14.67 \\
13 & IAH & 10.41 & 8.07 & 47.73 & 10.81 \\
14 & JFK & 3.38 & 2.56 & 56.20 & 4.95 \\
15 & LAS & 3.18 & 2.56 & 61.40 & 5.36 \\
16 & LAX & 3.66 & 2.79 & 65.80 & 6.37 \\
17 & LGA & 3.64 & 2.86 & 67.32 & 6.72 \\
18 & MCO & 2.70 & 2.05 & 57.13 & 4.23 \\
19 & MDW & 1.95 & 1.43 & 60.28 & 3.26 \\
20 & MEM & 1.05 & 0.54 & 89.16 & 3.00 \\
21 & MIA & 3.21 & 2.26 & 52.13 & 3.68 \\
22 & MSP & 6.17 & 4.74 & 46.04 & 7.16 \\
23 & ORD & 7.62 & 5.91 & 63.66 & 12.65 \\
24 & PHL & 4.54 & 3.62 & 53.89 & 7.13 \\
25 & PHX & 4.41 & 3.47 & 64.66 & 7.46 \\
26 & SAN & 2.01 & 1.44 & 59.75 & 3.11 \\
27 & SEA & 4.50 & 3.63 & 58.03 & 7.36 \\
28 & SFO & 3.49 & 2.81 & 57.14 & 5.49 \\
29 & SLC & 5.02 & 4.01 & 37.61 & 5.22 \\
30 & TPA & 1.69 & 1.17 & 67.98 & 3.16 \\
\hline
\end{tabular}
}
\begin{tablenotes}
\footnotesize
\item Note: Arrows indicate preferred direction: ($\downarrow$) lower is better, ($\uparrow$) higher is better. Airports follow the standard IATA codes. A complete list of airport names corresponding to IATA codes is provided in~\ref{appendix:airport_codes}.
\end{tablenotes}
\caption{Departure capacity distribution prediction performance for the core 30 US airports in the test set.}
\label{tab:prediction_results_departure}
\end{threeparttable}
\end{table}

\begin{figure}[htbp]
\centering
\begin{subfigure}{0.85\textwidth}
    \centering
    \includegraphics[width=\textwidth]{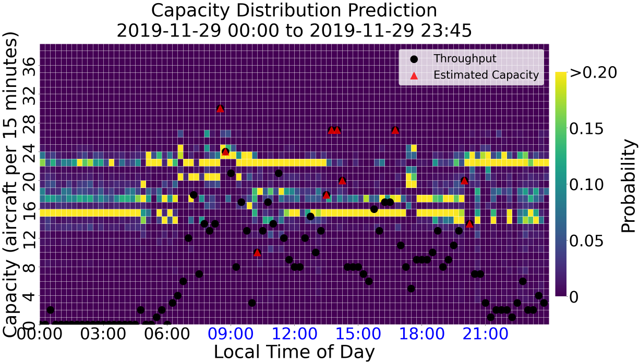}
    \caption{Departure}
    \label{fig:case1}
\end{subfigure}
\hfill
\begin{subfigure}{0.85\textwidth}
    \centering
    \includegraphics[width=\textwidth]{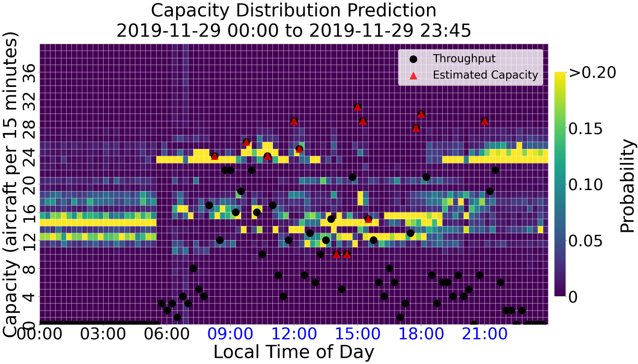}
    \caption{Arrival}
    \label{fig:case2}
\end{subfigure}
\caption{Capacity distribution prediction examples from November 29, 2019, showing predicted probability distributions (heatmap), actual throughput (black dots), and estimated capacity (red triangles) over a full operational day. The blue area from 9:00-21:00 is the 12-hour solution window of the DR-MAGHP.}
\label{fig:case_studies}
\end{figure}

Moreover, Figure~\ref{fig:case_studies} illustrates the model's capacity distribution prediction performance through two case studies from November 29, 2019. The probability heatmaps demonstrate that the model successfully captures temporal patterns of airport capacity variations, with higher probability concentrations during peak operational hours (06:00-22:00) and appropriate uncertainty quantification during transitional periods. The predicted distributions show concentrated probability mass around specific capacity values during stable operations, while becoming more dispersed during uncertain conditions. Comparing predicted distributions with actual throughput (black dots) and estimated capacity values (red triangles), the model shows reasonable alignment during morning and evening periods, but exhibits notable misalignment during midday operations where actual values frequently fall outside the high-probability regions. This visual evidence supports the quantitative findings of suboptimal PICP performance, indicating that while the model captures general capacity trends and temporal dependencies effectively, it requires enhanced uncertainty quantification to better account for operational disruptions and dynamic capacity-affecting factors for improved practical utility.

While the prediction results reveal certain limitations—such as coverage probabilities falling short of the nominal 90\% level and occasional misalignments during disruptive periods—they nevertheless capture meaningful temporal patterns and distributional features of airport capacities. In the next section, we examine whether these imperfect yet informative predictions can still enhance downstream decision-making compared to purely deterministic approaches, and whether the distributionally robust \textsc{dr-MAGHP} can further mitigate the impact of prediction errors by providing more robust ground delay policies.

\newpage
\subsection{Comparative Study of \textsc{det-MAGHP}, \textsc{sp-MAGHP} and \textsc{dr-MAGHP}}
Building on the predicted distributions from Section~\ref{sssec:prob_cap_pred}, we next evaluate how these predictions inform downstream decision-making. In particular, we compare the performance of prediction-driven models (\textsc{sp-MAGHP} and \textsc{dr-MAGHP}) against the deterministic baseline (\textsc{det-MAGHP}) in generating ground delay policies. This comparison allows us to assess both (\emph{i}) the extent to which upstream predictions improve decision quality relative to deterministic assumptions, and (\emph{ii}) the additional robustness provided by \textsc{dr-MAGHP} when distribution shifts occur between predicted and realized capacities (e.g., realized capacities being systematically lower than what was predicted).

\subsubsection{Experiment Setup}
To ensure consistency with the methodology used for predicting capacity distributions, we also construct flight schedules using BTS flight data from the year 2019. Due to the computational complexity of solving the \textsc{det-MAGHP}, \textsc{sp-MAGHP}, and \textsc{dr-MAGHP} models, our experiments focus on a selected set of representative days. To identify these days, we first compute the mean predicted capacity distribution for each airport and day over the full year. 
We then select the three dates exhibiting the greatest overestimation discrepancies, where predicted capacities were higher than the realized capacities. The selected dates used for evaluation are (in YYYY-MM-DD format): \text{2019-05-04}, \text{2019-05-30}, \text{2019-06-16}, \text{2019-11-13}, \text{2019-12-25}, and \text{2019-12-26}.

When computing airborne delays and enforcing connecting constraints, we restrict our attention to flights that depart from within the core 30 US airports. For flights originating from or arriving at out-of-network airports, we assume that those airports have unlimited arrival and departure capacities. As such, these flights are not subject to connecting constraints, and their airborne delays are calculated as the difference between scheduled and actual arrival times. 

The planning horizon for the \textsc{dr-MAGHP} model spans 12 hours (from 09:00 to 21:00, expressed in the local time of each airport), discretized into 48 time periods ($|T|=48$), each representing a 15-minute interval. We adopt a realistic minimum turnaround time of 45 minutes \cite{GHP1} and follow standard MAGHP formulations, including the use of an additional time period to accommodate excess flight volume \cite{GHP1}.
In terms of cost modeling, we assume that the unit cost of airborne holding delay is three times higher than that of ground holding delay, thereby establishing a 3:1 cost ratio. Furthermore, the costs associated with second-stage departure and arrival delays are assumed to be equivalent to those of airborne holding. This assumption reflects two factors: (i) second-stage departure delays typically involve aircraft waiting in taxiways or departure queues with engines running, leading to substantial fuel consumption\cite{delay_cost}; and (ii) such delays disrupt gate assignment schedules, which in turn can generate additional airborne holding for arriving flights.

For the Wasserstein ambiguity sets, we impose a uniform radius $\epsilon$ across all airports. The distance between any two scenarios in the scenario tree is measured using the $\ell_{2}$ norm. This metric is then used as the basis for constructing the Wasserstein ambiguity sets. Because the magnitude of inter-scenario distances is inherently scale-dependent, we normalize these distances prior to constructing the Wasserstein ambiguity sets. This procedure ensures that the ambiguity radius $\epsilon$ captures relative variations in the predicted capacity distributions, rather than being driven by the absolute scale of the underlying capacity values. Consequently, $\epsilon$ can be interpreted as a dimensionless measure of robustness, which facilitates consistent comparisons across airports and experimental settings with differing capacity ranges.

\subsubsection{Evaluation Framework}
We define the outputs from the upstream capacity distribution prediction procedure as predicted capacity distributions. We test the performance of the \textsc{DR-MAGHP} by comparing outcomes when the predicted capacity distributions, derived from the upstream prediction model, differ from the realized capacity distributions, which we adopt as testing distributions. Testing distributions of interest are the cases where realized capacity distributions are shifted to the left, i.e., the realized probabilistic capacities are lower than anticipated. We also solve the \textsc{det-MAGHP} and \textsc{sp-MAGHP} to compare their delay costs with those of the \textsc{dr-MAGHP} when predicted capacity distributions are not accurate. We use the terms \emph{in-sample performance} ($\phi_{IS}$) and \emph{out-of-sample performance} ($\phi_{OS}$) to refer to the costs (i.e., optimal value in expectation) of \textsc{det-MAGHP}, \textsc{sp-MAGHP} and \textsc{dr-MAGHP} when evaluated on predicted distributions and testing distributions respectively. Finally, we emphasize that these comparisons must be made on a day-by-day basis, just as unique GDP policies must be developed for each encountered NAS scenario. 

\begin{algorithm}
\caption{Capacity Resampling Algorithm}
\begin{algorithmic}[1]
    \STATE \textbf{Input}: Predicted PMF $\{\widehat{p}^{(1)},\widehat{p}^{(2)},\dots, \widehat{p}^{(|Z|)}\}$ of each airport's predicted capacity $\{\widehat{\xi}^{(1)},\widehat{\xi}^{(2)},\dots, \widehat{\xi}^{(|Z|)}\}$ with mean values $\{\widehat{\mu}^{(1)},\widehat{\mu}^{(2)},\dots,\widehat{\mu}^{|Z|}\}$, capacity reduction level $r$ for all airports and a maximum variability rate of probability $\delta$.
        \FOR{each airport $z$}
            \STATE Target mean $\mu^{*} = \widehat{\mu}^{(z)}  \left(1-r\right)$
            \STATE Update weights for each support $p_{(z)}$:
            \begin{equation}
            \begin{aligned}
                p^{(z)} = \arg\min_{p} \sum_{i=1}^{|\widehat{\xi}^{(z)}|}p_{i}\widehat{\xi}^{(z)}_{i} & \\
                \textrm{s.t.} \quad \sum_{i=1}^{|\widehat{\xi}^{(z)}|} p_{i}\widehat{\xi}^{(z)}_{i} &\geq \mu^{*},\quad\sum_{i=1}^{|\widehat{\xi}^{(z)}|}p_{i} = 1,\\
                p_{i} - \widehat{p}^{(z)}_{i} \leq \delta\widehat{p}^{(z)}_{i},\quad-\delta\widehat{p}^{(z)}_{i} &\leq p_{i} - \widehat{p}^{(z)}_{i},\quad \forall i \in N.
                \label{eq:lp_capreduc}
            \end{aligned}
            \end{equation}
        \ENDFOR
\STATE Draw i.i.d samples $\Tilde{\xi}^{(z)} \sim p^{(z)}, z = 1,2,\dots,|Z|$ \\
\STATE \textbf{Output}: Reduced testing capacities for all airports $\{\Tilde{\xi}^{(1)},\Tilde{\xi}^{(2)},\dots,\Tilde{\xi}^{(|Z|)}\}$
\end{algorithmic}
\label{Alg:alg1}
\end{algorithm}

\begin{table}[!htbp]
\centering
\resizebox{\textwidth}{!}{%
\begin{tabular}{llrrrccc}
\toprule
\textbf{Day} & \textbf{Reduction} & \textbf{Det.} & \textbf{Stoch.} & \textbf{DR} & \textbf{\%↓ vs Det.} & \textbf{\%↓ vs Stoch.} & $\epsilon^*$ \\
\midrule
\multirow{5}{*}{2019-11-13}
 & 10\% & \num{163693.06} & \best{\num{94115.79}} & \num{94181.56} & 42.5\% & -0.1\% & 0.00 \\
 & 20\% & \num{209085.40} & \num{121744.95} & \best{\num{117647.07}} & 43.7\% & 3.4\% & 0.04 \\
 & 30\% & \num{253888.48} & \num{152774.67} & \best{\num{141632.65}} & 44.2\% & 7.3\% & 0.09 \\
 & 40\% & \num{298865.38} & \num{183891.93} & \best{\num{164472.45}} & 45.0\% & 10.6\% & 0.10 \\
 & 50\% & \num{344322.97} & \num{221223.15} & \best{\num{192907.44}} & 44.0\% & 12.8\% & 0.10 \\
\midrule
\multirow{5}{*}{2019-12-25}
 & 10\% & \num{94752.95} & \num{56742.00} & \best{\num{56696.37}} & 40.2\% & 0.1\% & 0.00 \\
 & 20\% & \num{120845.96} & \num{73642.38} & \best{\num{70129.41}} & 42.0\% & 4.8\% & 0.05 \\
 & 30\% & \num{147568.70} & \num{92025.57} & \best{\num{82585.26}} & 44.0\% & 10.3\% & 0.06 \\
 & 40\% & \num{181714.55} & \num{117025.35} & \best{\num{100484.91}} & 44.7\% & 14.1\% & 0.06 \\
 & 50\% & \num{220189.49} & \num{148585.35} & \best{\num{126541.02}} & 42.5\% & 14.9\% & 0.07 \\
\midrule
\multirow{5}{*}{2019-12-26}
 & 10\% & \num{142287.07} & \num{86047.91} & \best{\num{85394.12}} & 40.0\% & 0.8\% & 0.02 \\
 & 20\% & \num{183530.89} & \num{114501.26} & \best{\num{108540.52}} & 40.9\% & 5.2\% & 0.03 \\
 & 30\% & \num{226631.65} & \num{143543.39} & \best{\num{128234.60}} & 43.4\% & 10.7\% & 0.09 \\
 & 40\% & \num{273066.79} & \num{177167.00} & \best{\num{151937.59}} & 44.4\% & 14.2\% & 0.10 \\
 & 50\% & \num{316751.41} & \num{212267.00} & \best{\num{179178.43}} & 43.5\% & 15.6\% & 0.10 \\
\midrule
\multirow{5}{*}{2019-05-04}
 & 10\% & \num{127984.75} & \num{70429.23} & \best{\num{69890.68}} & 45.4\% & 0.8\% & 0.03 \\
 & 20\% & \num{159001.87} & \num{90412.11} & \best{\num{85840.19}} & 46.0\% & 5.1\% & 0.04 \\
 & 30\% & \num{193583.71} & \num{115190.28} & \best{\num{106207.73}} & 45.1\% & 7.8\% & 0.04 \\
 & 40\% & \num{229113.76} & \num{144020.88} & \best{\num{131114.62}} & 42.8\% & 9.0\% & 0.05 \\
 & 50\% & \num{270145.24} & \num{181063.59} & \best{\num{162968.52}} & 39.7\% & 10.0\% & 0.09 \\
\midrule
\multirow{5}{*}{2019-05-30}
 & 10\% & \num{172395.43} & \num{96131.94} & \best{\num{95623.95}} & 44.5\% & 0.5\% & 0.01 \\
 & 20\% & \num{220632.73} & \num{127279.71} & \best{\num{121618.56}} & 44.9\% & 4.4\% & 0.06 \\
 & 30\% & \num{269383.51} & \num{162797.16} & \best{\num{146544.81}} & 45.6\% & 10.0\% & 0.06 \\
 & 40\% & \num{315159.25} & \num{199288.83} & \best{\num{173199.30}} & 45.0\% & 13.1\% & 0.10 \\
 & 50\% & \num{363924.94} & \num{241438.92} & \best{\num{207220.14}} & 43.1\% & 14.2\% & 0.10 \\
\midrule
\multirow{5}{*}{2019-06-16}
 & 10\% & \num{162692.60} & \best{\num{94799.79}} & \num{94824.74} & 41.7\% & -0.0\% & 0.00 \\
 & 20\% & \num{207200.21} & \num{121499.67} & \best{\num{114618.74}} & 44.7\% & 5.7\% & 0.05 \\
 & 30\% & \num{256369.01} & \num{155088.27} & \best{\num{138248.20}} & 46.1\% & 10.9\% & 0.06 \\
 & 40\% & \num{303713.81} & \num{191647.74} & \best{\num{167155.00}} & 45.0\% & 12.8\% & 0.10 \\
 & 50\% & \num{344636.30} & \num{227784.66} & \best{\num{199335.52}} & 42.2\% & 12.5\% & 0.10 \\
\bottomrule
\end{tabular}}
\caption{Out-of-sample cost comparison across different days and capacity reductions. Columns show deterministic (Det.), stochastic (Stoch.), and distributionally robust (DR) MAGHP. Bold values indicate the lowest (best) cost in each row.}
\label{tab:maghp_cost_comparison}
\end{table}

\subsubsection{Sensitivity Analysis Framework}
A sensitivity analysis is conducted to assess the impact of discrepancies between predicted and realized capacity distributions on delay costs, as well as the ability of \textsc{dr-MAGHP} to mitigate these impacts. This analysis generates testing distributions at various levels of capacity reductions to compare the out-of-sample performance of \textsc{det-MAGHP}, \textsc{sp-MAGHP}, and \textsc{dr-MAGHP}. Ground delay policies generated by \textsc{sp-MAGHP} and \textsc{dr-MAGHP} are constructed using the full predicted capacity distributions, capturing the uncertainty in future capacities. In contrast, \textsc{det-MAGHP} derives its policies from the single best-capacity scenario of the day, reflecting an optimistic assumption that airports will operate at their maximum predicted capacities. 

\subsubsection{Construction of Reduced Capacity Distributions}
To sample from reduced capacity distributions at various reduction levels, we introduce a linear program in Equation \eqref{eq:lp_capreduc} that performs valid adjustments to the PMFs to minimize the deviation between the current mean value and the targeted (reduced) mean value, while maintaining the probabilistic properties of the weights. We also introduce a parameter $\delta$ for the maximum variability rate to ensure a more uniform distribution of probability mass. 

The details of the sampling procedure are presented in Algorithm~\ref{Alg:alg1}. Given a general ground delay policy $\Tilde{x}$, we generate 100 samples from the reduced capacity distributions from Algorithm \ref{Alg:alg1} and compute the out-of-sample performance using the expression $\phi_{OS}(\Tilde{x}) = \sum_{i=1}^{|\Tilde{\xi}|} \widetilde{\phi_{OS}}(\Tilde{x}, \Tilde{\xi}_{i}) / |\Tilde{\xi}|$, where $\widetilde{\phi_{OS}}(\Tilde{x}, \Tilde{\xi}_{i})$ denotes the objective value of the proposed ground holding policy $\Tilde{x}$ under the $i\textsuperscript{th}$ capacity sample $\Tilde{\xi}_{i}$. We emphasize that $\widetilde{\phi_{OS}}$ corresponds to an optimization model structurally equivalent to \textsc{sp-MAGHP}, sharing the same objective function and constraint set. The key distinction lies in the fact that $\widetilde{\phi_{OS}}$ accepts a fixed first-stage ground delay decision $\Tilde{x}$ as input, and evaluates its second-stage performance under a specified realization of the capacity distribution. 
\subsubsection{\emph{In-Sample} \emph{vs.} \emph{Out-of-Sample} Performance}
We investigate the impact of erroneous capacity predictions, specifically cases in which realized capacities are reduced by 10\% to 50\% relative to predicted capacity distributions. To evaluate both in-sample and out-of-sample performance of \textsc{dr-MAGHP} under varying levels of distributional robustness, we conduct a sweep over the Wasserstein radius $\epsilon$, ranging from 0 to 0.1. It is important to note that, since inter-scenario distances in the scenario tree are normalized, even relatively small values of $\epsilon$ can induce substantial changes in the objective function. Intuitively, normalization compresses the scale of distances so that scenarios appear closer together, which makes the Wasserstein ball ``tighter’’ around the predicted distribution. As a result, a modest increase in $\epsilon$ is sufficient to admit a much larger set of probability distributions into the ambiguity set, thereby altering the optimization landscape and, in turn, the optimal objective function value (i.e., ground delay policy cost).

We begin by examining the in-sample performance of \textsc{sp-MAGHP} and \textsc{dr-MAGHP} across all six selected evaluation dates, using the predicted capacity distributions. As illustrated in Figure~\ref{fig:In_sample_costs}, both models yield identical objective values when $\epsilon = 0$. This is expected, since in this case the Wasserstein ambiguity set collapses to a trivial set with one element, which is exactly the predicted distribution. Solving \textsc{dr-MAGHP} over this trivial ambiguity set is equivalent to solving the original \textsc{sp-MAGHP} problem. As $\epsilon$ increases from 0 to 0.1, we observe a consistent rise in the objective values of \textsc{dr-MAGHP} across all test dates. This trend indicates that greater distributional robustness—i.e., larger ambiguity sets—leads to deteriorated in-sample performance when evaluated on the nominal predicted capacities. However, with regard to out-of-sample performance, when the realized capacities deviate only slightly from the predictions (i.e., a 10\% reduction), a small ambiguity set (i.e., small $\epsilon$) is sufficient for \textsc{dr-MAGHP} to moderately outperform \textsc{sp-MAGHP} on most evaluation dates, including on  \text{2019-12-25}, \text{2019-12-26}, \text{2019-05-04}, and \text{2019-05-30}. However, for \text{2019-11-13} and \text{2019-06-16}, the out-of-sample performance of \textsc{sp-MAGHP} is slightly better that of \textsc{dr-MAGHP} across all tested values of $\epsilon$. 
\subsubsection{Effect of the Wasserstein Radius $\epsilon$}
It is important to recognize that the out-of-sample performance of \textsc{sp-MAGHP} and \textsc{dr-MAGHP} may diverge, even when their in-sample objective values are identical under a Wasserstein radius of $\epsilon = 0$. This divergence arises from differences in the delay assignment policies produced by the two models, which, despite being theoretically equivalent at $\epsilon = 0$, can vary due to numerical tolerances and solver precision. As a result, matching in-sample performance does not necessarily imply identical out-of-sample behavior. 

Furthermore, as the gap between realized and predicted capacity distributions widens, larger ambiguity sets (i.e., higher values of $\epsilon$) are typically needed for \textsc{dr-MAGHP} to deliver improved out-of-sample performance. Table~\ref{tab:maghp_cost_comparison} summarizes the out-of-sample delay assignment costs associated with \textsc{det-MAGHP}, \textsc{sp-MAGHP}, and \textsc{dr-MAGHP} across varying levels of capacity reductions, along with the corresponding optimal ambiguity set radius $\epsilon^*$. Notably, Table~\ref{tab:maghp_cost_comparison} reveals a monotonic increase in the optimal $\epsilon^*$ as the severity of the capacity reduction intensifies. Additionally, we observe that the out-of-sample cost of \textsc{dr-MAGHP} tends to increase again once $\epsilon$ exceeds its optimal value $\epsilon^*$. Although this trend may plateau at higher levels due to the experimental cap of $\epsilon = 0.1$, it underscores the trade-off between robustness and conservativeness inherent in the design of distributionally robust policies. 
\subsubsection{Comparative Performance under Capacity Reductions}
Table~\ref{tab:maghp_cost_comparison} further reveals that under substantial capacity reductions—specifically, when realized capacities are 30\%, 40\%, or 50\% lower than the predicted values—the \textsc{dr-MAGHP} model achieves notable cost savings relative to even the \textsc{sp-MAGHP}, ranging from 7.30\% to 10.85\% at 30\% reductions, 8.97\% to 14.21\% at 40\% reductions, and 9.95\% to 15.58\% at 50\% reductions. These reductions represent cases where the predictions are overly optimistic compared to the realized (testing) distributions. Moreover, the improvements are even more pronounced when compared with \textsc{det-MAGHP}, exceeding 40\% in all tested scenarios. These results demonstrate that robustified ground delay policies can yield meaningful cost reductions under adverse operational conditions. 

To better understand this effect, we examine the results in Tables~\ref{tab:maghp_arrival_delay} and~\ref{tab:maghp_departure_delay}. These tables show that policies derived from \textsc{dr-MAGHP} consistently lead to lower second-stage departure and arrival delays, with the effect being more pronounced for departures. Because second-stage delays carry higher unit costs than first-stage delays, the conservative nature of \textsc{dr-MAGHP}—which may increase first-stage delay costs—ultimately reduces overall costs by mitigating more expensive downstream disruptions, particularly in scenarios with high degrees of disruptions across the NAS.

\subsubsection{Key Findings}
In summary, the numerical experiments provide three key insights: (1) leveraging airport capacity predictions enables data-driven models such as \textsc{sp-MAGHP} and \textsc{dr-MAGHP} to reduce operational costs; (2) when capacity predictions are expected to be reliable, \textsc{sp-MAGHP} may be preferable due to its cost efficiency, whereas in more uncertain conditions, \textsc{dr-MAGHP} provides more robust performance; and (3) the choice of the Wasserstein radius $\epsilon$ plays a critical role, balancing conservativeness and robustness in decision making.

\section{Concluding Remarks}
\subsection{Summary}
This paper proposed an integrated framework that combines learning-based prediction and Distributionally Robust Optimization (DRO) for Air Traffic Management (ATM), with a focus on flight rescheduling decision-making within the context of Ground Delay Programs (GDP), given uncertainty in airport arrival and departure capacities. The upstream module utilizes a machine learning model to generate probabilistic airport capacity predictions across the planning horizon. These predictions are then used by the downstream \textsc{DR-MAGHP} model to construct Wasserstein ambiguity sets, enabling the generation of distributionally robust ground delay policies.

Our results show that prediction-based models such as \textsc{sp-MAGHP} and \textsc{dr-MAGHP} substantially reduce operational costs compared to the baseline \textsc{det-MAGHP} model that does not incorporate predictive information. Notably, \textsc{dr-MAGHP} demonstrates superior performance under distribution shifts, highlighting its ability to better manage uncertainty in airport capacities. This work thus bridges data-driven prediction and robust optimization for strategic air traffic management decision-making, offering tractable reformulations and scalable scenario reduction techniques for real-world deployment. It also identifies the conditions under which DRO-based models outperform traditional stochastic optimization approaches.

\subsection{Limitations and Future Work}
While our approach demonstrates strong performance on data from the FAA Core 30 airports, it is built on several simplifying assumptions that may limit its fidelity to real-world operations. First, we assume independence between arrival and departure capacities at each airport, as well as independence across airports. In reality, operational dependencies are common: adverse weather or other disruptions often reduce both arrival and departure capacities simultaneously at a single airport, and severe regional weather phenomena can affect multiple neighboring airports at once. Second, we impose a uniform value of $\epsilon$ across all airports when constructing ambiguity sets. This assumption overlooks the heterogeneity of uncertainty across airports, as different airports face distinct traffic patterns, weather conditions, and operational environments that may warrant different ambiguity set sizes. Finally, due to limited access to historical airport capacity data, model evaluation was conducted using synthetically generated capacity reduction scenarios, which may not fully replicate the complexity of real-world disruptions.

Future research will focus on several extensions to enhance the realism and applicability of the proposed framework. First, we aim to model the correlation between arrival and departure capacities at individual airports. Second, we plan to incorporate inter-airport correlations to better capture the regional impacts of weather and operational disruptions. Third, we will explore assigning heterogeneous levels of uncertainty to different airports, allowing the ambiguity set size to reflect local traffic patterns, weather variability, and operational characteristics, and investigate the trade-off between efficiency and robustness for each airport. Finally, we will evaluate the proposed DRO models on historical disruption days to assess their ability to mitigate the negative impacts of adverse events. Together, these directions are expected to strengthen the practical relevance of the framework for air traffic management.
\appendix
\section{Number of time periods remaining after rule-based capacity estimation for each airport.}\label{appendix:rule_estimation_table}
\begin{table}[H]
\centering
\small
\begin{threeparttable}
\begin{tabular}{c|c|cc|cc}
\hline
\textbf{Index} & \textbf{Airport} & \multicolumn{2}{c|}{\textbf{Arrival}} & \multicolumn{2}{c}{\textbf{Departure}} \\
& & \textbf{Time Periods} & \textbf{\%} & \textbf{Time Periods} & \textbf{\%} \\
\hline
1  & ATL & 12,816 & 36.6\% & 14,120 & 40.3\% \\
2  & BOS & 10,023 & 28.6\% & 9,730 & 27.8\% \\
3  & BWI & 6,758 & 19.3\% & 8,430 & 24.1\% \\
4  & CLT & 9,503 & 27.1\% & 11,578 & 33.0\% \\
5  & DCA & 9,065 & 25.9\% & 8,467 & 24.2\% \\
6  & DEN & 13,004 & 37.1\% & 13,714 & 39.1\% \\
7  & DFW & 12,548 & 35.8\% & 13,958 & 39.8\% \\
8  & DTW & 7,556 & 21.6\% & 7,548 & 21.5\% \\
9  & EWR & 12,279 & 35.0\% & 11,472 & 32.7\% \\
10 & FLL & 6,848 & 19.5\% & 7,616 & 21.7\% \\
11 & HNL & 5,370 & 15.3\% & 4,323 & 12.3\% \\
12 & IAD & 4,981 & 14.2\% & 4,548 & 13.0\% \\
13 & IAH & 8,288 & 23.7\% & 8,154 & 23.3\% \\
14 & JFK & 7,813 & 22.3\% & 7,744 & 22.1\% \\
15 & LAS & 10,060 & 28.7\% & 10,625 & 30.3\% \\
16 & LAX & 14,202 & 40.5\% & 13,421 & 38.3\% \\
17 & LGA & 11,891 & 33.9\% & 11,078 & 31.6\% \\
18 & MCO & 8,673 & 24.8\% & 10,014 & 28.6\% \\
19 & MDW & 6,365 & 18.2\% & 8,360 & 23.9\% \\
20 & MEM & 7,391 & 21.1\% & 7,101 & 20.3\% \\
21 & MIA & 5,394 & 15.4\% & 5,929 & 16.9\% \\
22 & MSP & 7,192 & 20.5\% & 6,722 & 19.2\% \\
23 & ORD & 16,395 & 46.8\% & 16,820 & 48.0\% \\
24 & PHL & 8,173 & 23.3\% & 7,962 & 22.7\% \\
25 & PHX & 9,721 & 27.7\% & 9,714 & 27.7\% \\
26 & SAN & 6,989 & 19.9\% & 7,312 & 20.9\% \\
27 & SEA & 9,784 & 27.9\% & 9,905 & 28.3\% \\
28 & SFO & 12,682 & 36.2\% & 10,931 & 31.2\% \\
29 & SLC & 6,059 & 17.3\% & 5,400 & 15.4\% \\
30 & TPA & 5,737 & 16.4\% & 6,239 & 17.8\% \\
\hline
\textbf{Total} & & \textbf{270,660} & \textbf{25.7\%} & \textbf{280,376} & \textbf{26.7\%} \\
\hline
\end{tabular}
\caption{Number of time periods remaining after rule-based capacity estimation for each airport.}
\end{threeparttable}
\end{table}
\section{IATA airport codes and corresponding airport names of the core 30 US airports}\label{appendix:airport_codes}
\begin{table}[H]
\centering
\small
\begin{tabular}{ccc}
\hline
\textbf{Index} & \textbf{IATA Code} & \textbf{Airport Name} \\
\hline
1 & ATL & Hartsfield-Jackson Atlanta International Airport \\
2 & BOS & Logan International Airport \\
3 & BWI & Baltimore/Washington International Airport \\
4 & CLT & Charlotte Douglas International Airport \\
5 & DCA & Ronald Reagan Washington National Airport \\
6 & DEN & Denver International Airport \\
7 & DFW & Dallas/Fort Worth International Airport \\
8 & DTW & Detroit Metropolitan Wayne County Airport \\
9 & EWR & Newark Liberty International Airport \\
10 & FLL & Fort Lauderdale-Hollywood International Airport \\
11 & HNL & Daniel K. Inouye International Airport \\
12 & IAD & Washington Dulles International Airport \\
13 & IAH & George Bush Intercontinental Airport \\
14 & JFK & John F. Kennedy International Airport \\
15 & LAS & McCarran International Airport \\
16 & LAX & Los Angeles International Airport \\
17 & LGA & LaGuardia Airport \\
18 & MCO & Orlando International Airport \\
19 & MDW & Chicago Midway International Airport \\
20 & MEM & Memphis International Airport \\
21 & MIA & Miami International Airport \\
22 & MSP & Minneapolis-Saint Paul International Airport \\
23 & ORD & Chicago O'Hare International Airport \\
24 & PHL & Philadelphia International Airport \\
25 & PHX & Phoenix Sky Harbor International Airport \\
26 & SAN & San Diego International Airport \\
27 & SEA & Seattle-Tacoma International Airport \\
28 & SFO & San Francisco International Airport \\
29 & SLC & Salt Lake City International Airport \\
30 & TPA & Tampa International Airport \\
\hline
\end{tabular}
\caption{IATA airport codes and corresponding airport names.}
\end{table}

\section{Second stage delay comparison}\label{appendix:second_stage_delay_comparison}
\begin{table}[H]
\centering
\small
\resizebox{\textwidth}{!}{%
\begin{tabular}{llrrrcc}
\toprule
\textbf{Day} & \textbf{Reduction} & \textbf{Det.} & \textbf{Stoch.} & \textbf{DR} & \textbf{\%↓ vs Det.} & \textbf{\%↓ vs Stoch.} \\
\midrule
\multirow{5}{*}{2019-11-13}
 & 10\% & \num{1.87e4} & \best{\num{8.42e3}} & \num{8.46e3} & 54.8\% & -0.5\% \\
 & 20\% & \num{2.41e4} & \num{1.26e4} & \best{\num{1.09e4}} & 54.6\% & 13.1\% \\
 & 30\% & \num{2.89e4} & \num{1.65e4} & \best{\num{1.31e4}} & 54.7\% & 20.7\% \\
 & 40\% & \num{3.37e4} & \num{2.04e4} & \best{\num{1.61e4}} & 52.4\% & 21.3\% \\
 & 50\% & \num{3.89e4} & \num{2.50e4} & \best{\num{2.00e4}} & 48.5\% & 19.8\% \\
\midrule
\multirow{5}{*}{2019-12-25}
 & 10\% & \num{1.05e4} & \num{4.36e3} & \best{\num{4.34e3}} & 58.6\% & 0.5\% \\
 & 20\% & \num{1.40e4} & \num{6.90e3} & \best{\num{5.33e3}} & 61.9\% & 22.7\% \\
 & 30\% & \num{1.72e4} & \num{9.24e3} & \best{\num{6.98e3}} & 59.4\% & 24.5\% \\
 & 40\% & \num{2.12e4} & \num{1.25e4} & \best{\num{9.52e3}} & 55.2\% & 23.5\% \\
 & 50\% & \num{2.53e4} & \num{1.59e4} & \best{\num{1.22e4}} & 51.7\% & 23.0\% \\
\midrule
\multirow{5}{*}{2019-12-26}
 & 10\% & \num{1.42e4} & \num{6.86e3} & \best{\num{6.45e3}} & 54.7\% & 5.9\% \\
 & 20\% & \num{1.92e4} & \num{1.05e4} & \best{\num{9.13e3}} & 52.5\% & 13.4\% \\
 & 30\% & \num{2.38e4} & \num{1.38e4} & \best{\num{1.08e4}} & 54.5\% & 21.4\% \\
 & 40\% & \num{2.88e4} & \num{1.77e4} & \best{\num{1.38e4}} & 52.2\% & 22.3\% \\
 & 50\% & \num{3.37e4} & \num{2.20e4} & \best{\num{1.74e4}} & 48.3\% & 20.8\% \\
\midrule
\multirow{5}{*}{2019-05-04}
 & 10\% & \num{1.34e4} & \num{6.12e3} & \best{\num{5.51e3}} & 58.8\% & 10.0\% \\
 & 20\% & \num{1.71e4} & \num{8.82e3} & \best{\num{7.51e3}} & 56.2\% & 14.9\% \\
 & 30\% & \num{2.19e4} & \num{1.25e4} & \best{\num{1.04e4}} & 52.7\% & 16.9\% \\
 & 40\% & \num{2.67e4} & \num{1.65e4} & \best{\num{1.36e4}} & 49.3\% & 17.8\% \\
 & 50\% & \num{3.15e4} & \num{2.09e4} & \best{\num{1.62e4}} & 48.7\% & 22.7\% \\
\midrule
\multirow{5}{*}{2019-05-30}
 & 10\% & \num{1.79e4} & \num{8.89e3} & \best{\num{8.38e3}} & 53.1\% & 5.7\% \\
 & 20\% & \num{2.40e4} & \num{1.32e4} & \best{\num{1.07e4}} & 55.5\% & 19.4\% \\
 & 30\% & \num{2.92e4} & \num{1.70e4} & \best{\num{1.38e4}} & 52.7\% & 18.5\% \\
 & 40\% & \num{3.48e4} & \num{2.13e4} & \best{\num{1.66e4}} & 52.3\% & 22.3\% \\
 & 50\% & \num{4.00e4} & \num{2.61e4} & \best{\num{2.09e4}} & 47.9\% & 20.0\% \\
\midrule
\multirow{5}{*}{2019-06-16}
 & 10\% & \num{1.89e4} & \num{8.82e3} & \best{\num{8.78e3}} & 53.6\% & 0.5\% \\
 & 20\% & \num{2.43e4} & \num{1.26e4} & \best{\num{1.10e4}} & 54.6\% & 12.4\% \\
 & 30\% & \num{2.97e4} & \num{1.65e4} & \best{\num{1.41e4}} & 52.5\% & 14.5\% \\
 & 40\% & \num{3.52e4} & \num{2.11e4} & \best{\num{1.68e4}} & 52.3\% & 20.2\% \\
 & 50\% & \num{4.00e4} & \num{2.54e4} & \best{\num{2.08e4}} & 47.9\% & 18.0\% \\
\bottomrule
\end{tabular}}
\caption{Second-stage \emph{arrival delay} comparison (in delay units) among deterministic (Det.), stochastic (Stoch.), and distributionally robust (DR) MAGHP across different days and capacity reductions. Bold values denote the lowest delay in each row.}
\label{tab:maghp_arrival_delay}
\end{table}

\begin{table}[H]
\centering
\small
\resizebox{\textwidth}{!}{%
\begin{tabular}{llrrrcc}
\toprule
\textbf{Day} & \textbf{Reduction} & \textbf{Det.} & \textbf{Stoch.} & \textbf{DR} & \textbf{\%↓ vs Det.} & \textbf{\%↓ vs Stoch.} \\
\midrule
\multirow{5}{*}{2019-11-13}
 & 10\% & \num{3.38e4} & \best{\num{5.97e3}} & \num{5.98e3} & 82.3\% & -0.2\% \\
 & 20\% & \num{4.36e4} & \num{1.10e4} & \best{\num{7.87e3}} & 81.9\% & 28.7\% \\
 & 30\% & \num{5.37e4} & \num{1.75e4} & \best{\num{8.07e3}} & 85.0\% & 53.8\% \\
 & 40\% & \num{6.38e4} & \num{2.39e4} & \best{\num{1.21e4}} & 81.1\% & 49.5\% \\
 & 50\% & \num{7.38e4} & \num{3.18e4} & \best{\num{1.76e4}} & 76.2\% & 44.7\% \\
\midrule
\multirow{5}{*}{2019-12-25}
 & 10\% & \num{2.02e4} & \best{\num{4.65e3}} & \num{4.66e3} & 76.9\% & 0.0\% \\
 & 20\% & \num{2.54e4} & \num{7.74e3} & \best{\num{3.79e3}} & 85.1\% & 51.0\% \\
 & 30\% & \num{3.11e4} & \num{1.15e4} & \best{\num{5.40e3}} & 82.6\% & 53.2\% \\
 & 40\% & \num{3.84e4} & \num{1.66e4} & \best{\num{8.82e3}} & 77.0\% & 47.0\% \\
 & 50\% & \num{4.72e4} & \num{2.38e4} & \best{\num{1.37e4}} & 71.0\% & 42.4\% \\
\midrule
\multirow{5}{*}{2019-12-26}
 & 10\% & \num{3.12e4} & \num{6.75e3} & \best{\num{6.07e3}} & 80.6\% & 10.0\% \\
 & 20\% & \num{4.00e4} & \num{1.25e4} & \best{\num{9.90e3}} & 75.3\% & 21.1\% \\
 & 30\% & \num{4.98e4} & \num{1.90e4} & \best{\num{6.99e3}} & 86.0\% & 63.2\% \\
 & 40\% & \num{6.03e4} & \num{2.63e4} & \best{\num{1.13e4}} & 81.2\% & 56.8\% \\
 & 50\% & \num{6.99e4} & \num{3.37e4} & \best{\num{1.68e4}} & 76.0\% & 50.2\% \\
\midrule
\multirow{5}{*}{2019-05-04}
 & 10\% & \num{2.72e4} & \num{4.81e3} & \best{\num{3.59e3}} & 86.8\% & 25.5\% \\
 & 20\% & \num{3.38e4} & \num{8.78e3} & \best{\num{6.13e3}} & 81.9\% & 30.1\% \\
 & 30\% & \num{4.05e4} & \num{1.34e4} & \best{\num{1.01e4}} & 75.2\% & 24.8\% \\
 & 40\% & \num{4.75e4} & \num{1.90e4} & \best{\num{1.31e4}} & 72.5\% & 31.2\% \\
 & 50\% & \num{5.64e4} & \num{2.69e4} & \best{\num{1.76e4}} & 68.8\% & 34.5\% \\
\midrule
\multirow{5}{*}{2019-05-30}
 & 10\% & \num{3.71e4} & \num{6.55e3} & \best{\num{6.42e3}} & 82.7\% & 2.0\% \\
 & 20\% & \num{4.71e4} & \num{1.26e4} & \best{\num{5.76e3}} & 87.8\% & 54.3\% \\
 & 30\% & \num{5.81e4} & \num{2.07e4} & \best{\num{1.09e4}} & 81.3\% & 47.3\% \\
 & 40\% & \num{6.78e4} & \num{2.85e4} & \best{\num{1.47e4}} & 78.3\% & 48.3\% \\
 & 50\% & \num{7.88e4} & \num{3.78e4} & \best{\num{2.18e4}} & 72.4\% & 42.4\% \\
\midrule
\multirow{5}{*}{2019-06-16}
 & 10\% & \num{3.25e4} & \best{\num{6.82e3}} & \num{6.85e3} & 78.9\% & -0.4\% \\
 & 20\% & \num{4.20e4} & \num{1.20e4} & \best{\num{6.31e3}} & 85.0\% & 47.2\% \\
 & 30\% & \num{5.30e4} & \num{1.92e4} & \best{\num{1.03e4}} & 80.6\% & 46.6\% \\
 & 40\% & \num{6.32e4} & \num{2.69e4} & \best{\num{1.35e4}} & 78.6\% & 49.8\% \\
 & 50\% & \num{7.20e4} & \num{3.45e4} & \best{\num{2.02e4}} & 72.0\% & 41.6\% \\
\bottomrule
\end{tabular}}
\caption{Second-stage \emph{departure delay} comparison (in delay units) among deterministic (Det.), stochastic (Stoch.), and distributionally robust (DR) MAGHP across different days and capacity reductions. Bold values denote the lowest delay in each row.}
\label{tab:maghp_departure_delay}
\end{table}

\section{\emph{Out-of-sample} performance of \textsc{sp-MAGHP} and \textsc{dr-MAGHP}}\label{appendix:OSP_results}
\begin{figure}
    \centering
    \begin{subfigure}[b]{0.49\textwidth}
        \includegraphics[width=\textwidth]{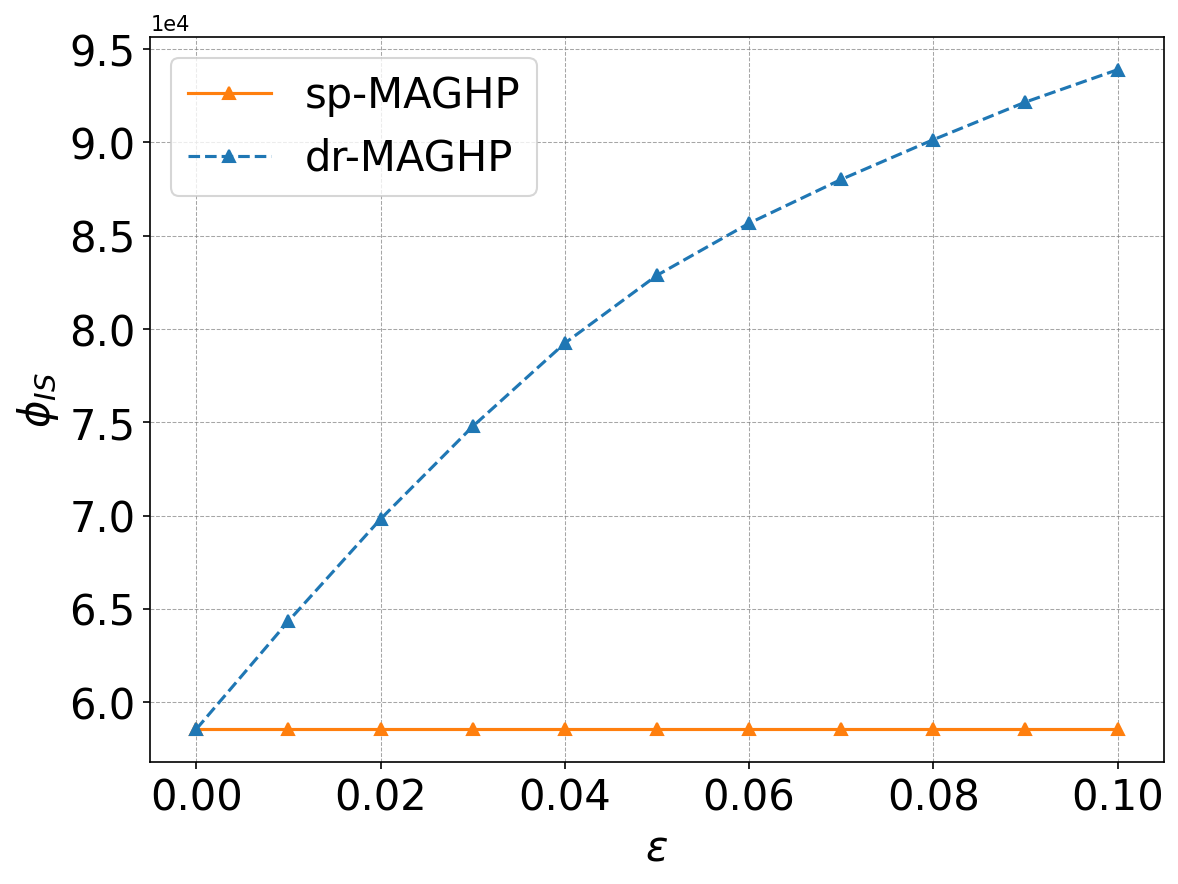}
        \caption{2019-05-04}
    \end{subfigure}
    \hfill
    \begin{subfigure}[b]{0.49\textwidth}
        \includegraphics[width=\textwidth]{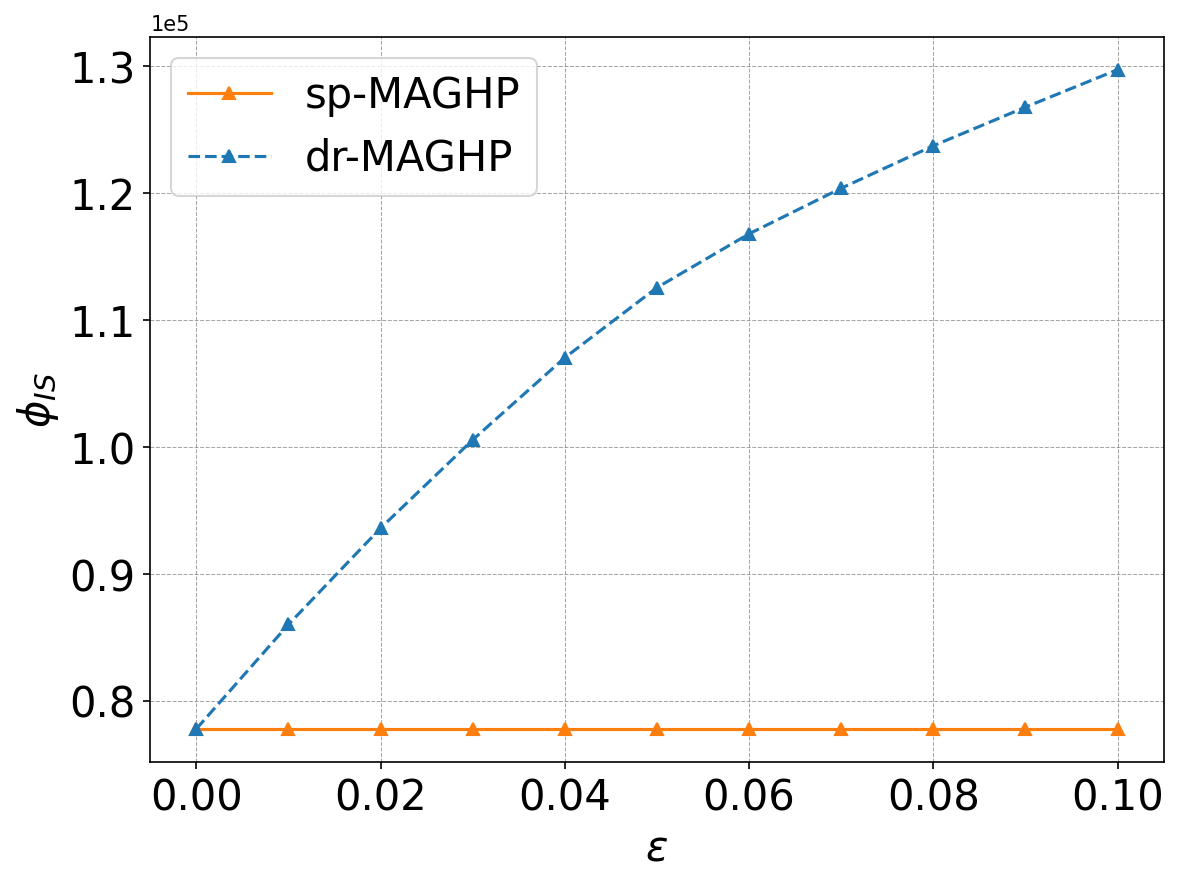}
        \caption{2019-05-30}
    \end{subfigure}

    \vspace{0.7em}

    \begin{subfigure}[b]{0.49\textwidth}
        \includegraphics[width=\textwidth]{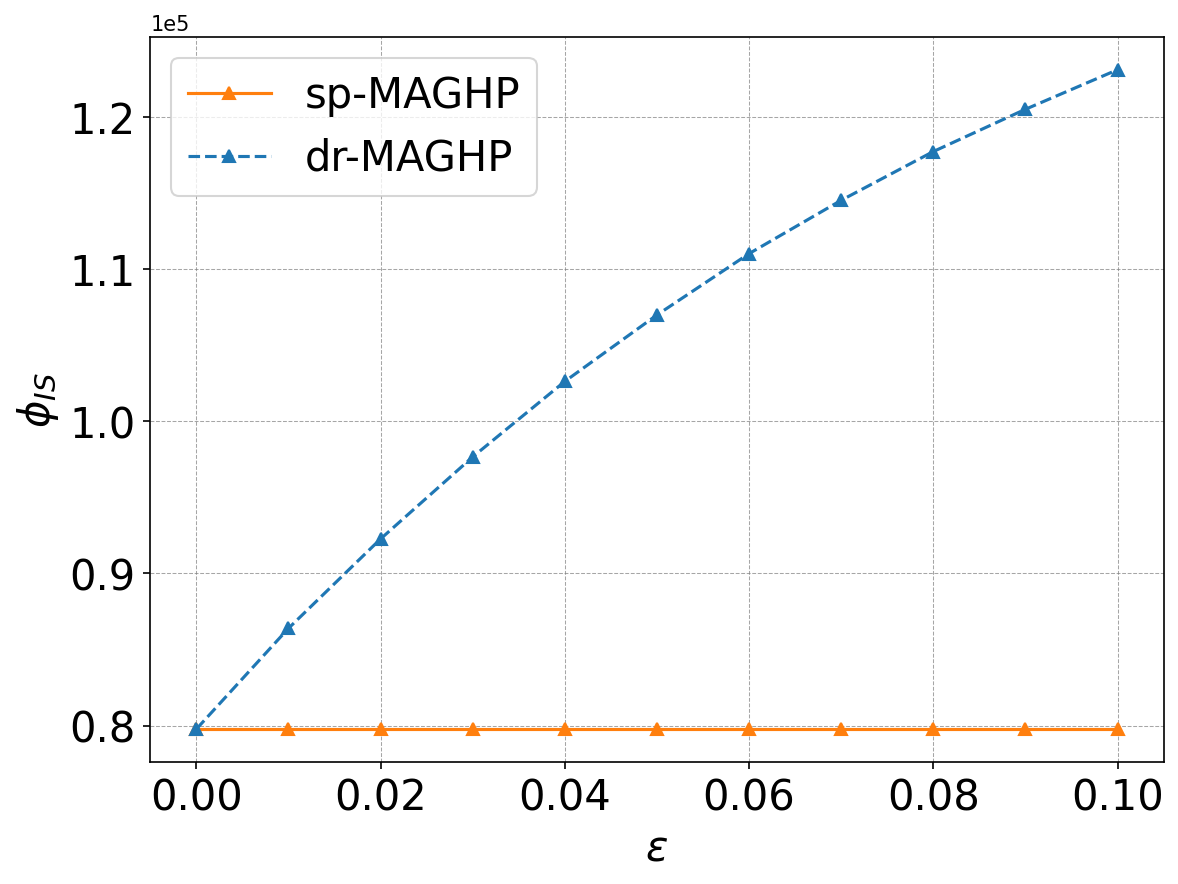}
        \caption{2019-06-16}
    \end{subfigure}
    \hfill
    \begin{subfigure}[b]{0.49\textwidth}
        \includegraphics[width=\textwidth]{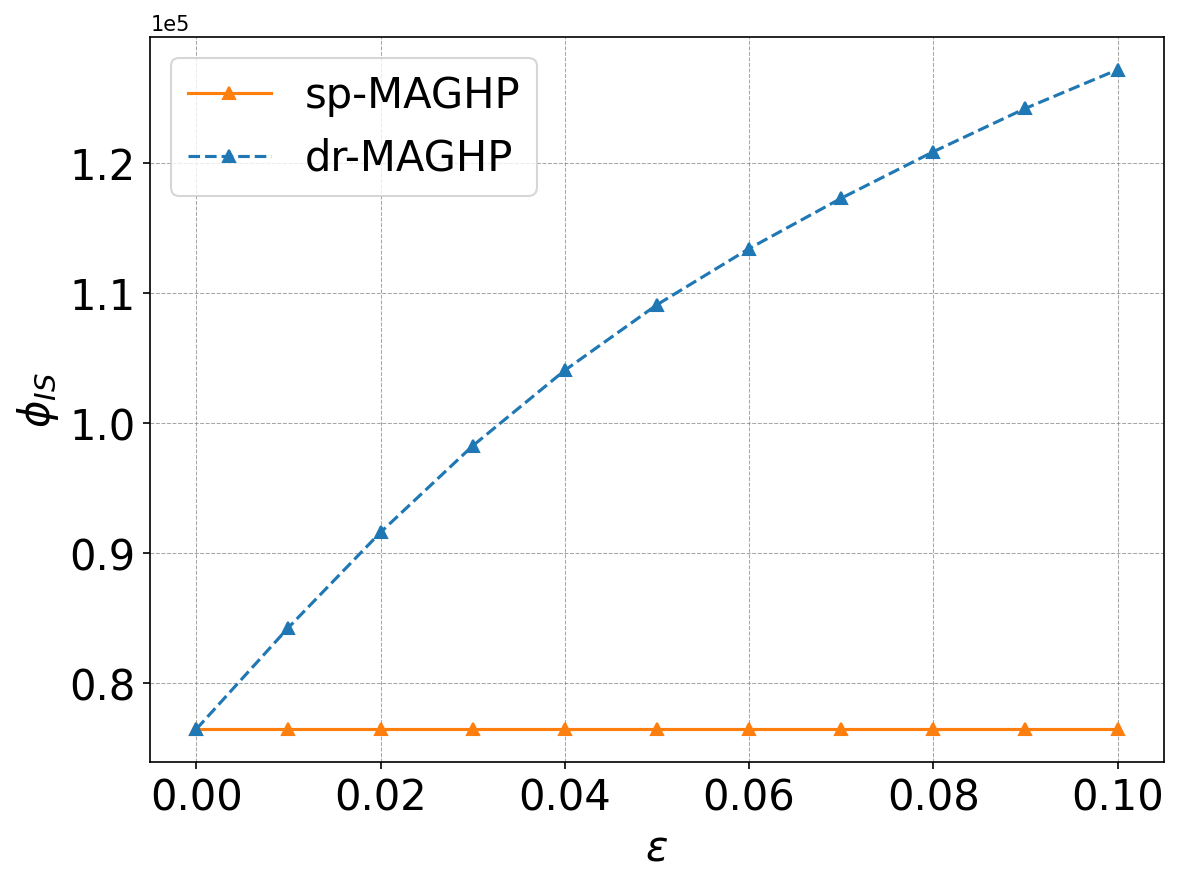}
        \caption{2019-11-13}
    \end{subfigure}

    \vspace{0.7em}

    \begin{subfigure}[b]{0.5\textwidth}
        \centering
        \includegraphics[width=\textwidth]{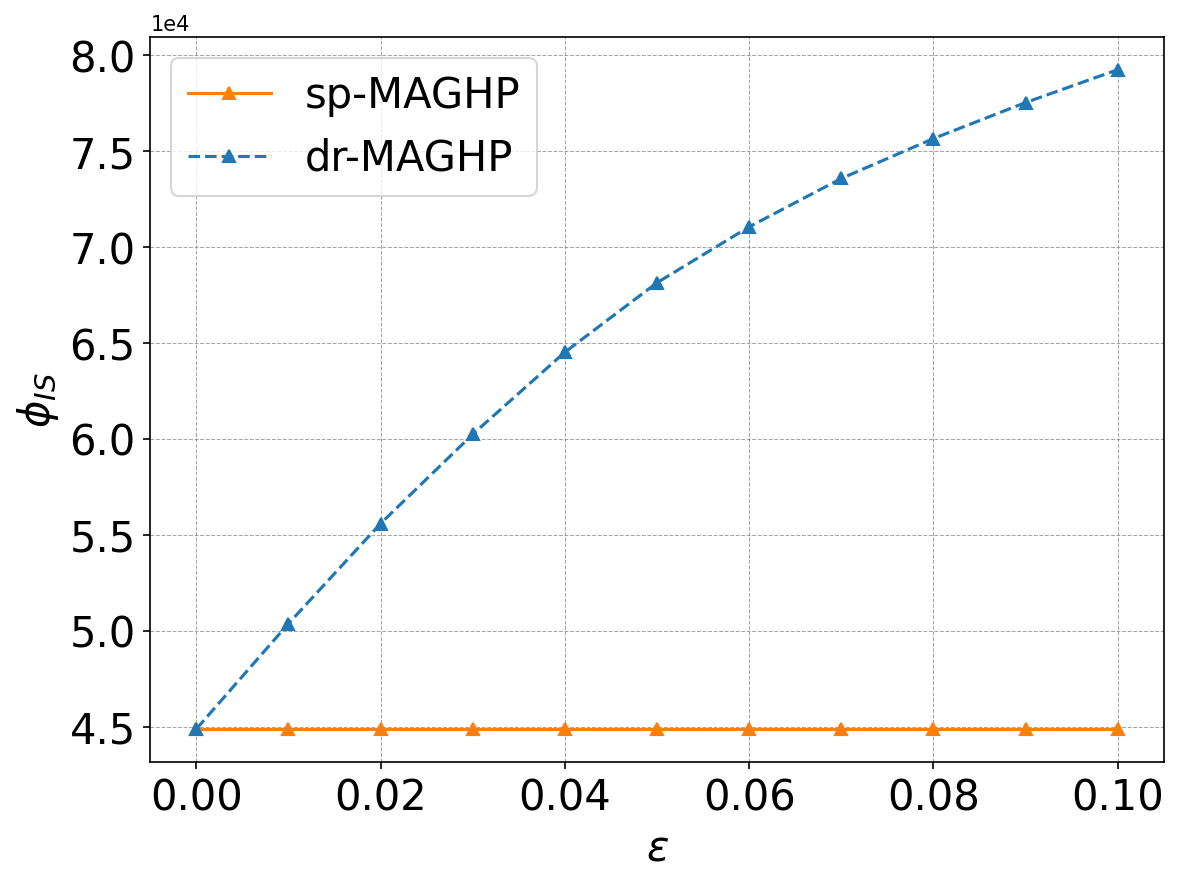}
        \caption{2019-12-25}
    \end{subfigure}
    \hfill
    \begin{subfigure}[b]{0.49\textwidth}
        \includegraphics[width=\textwidth]{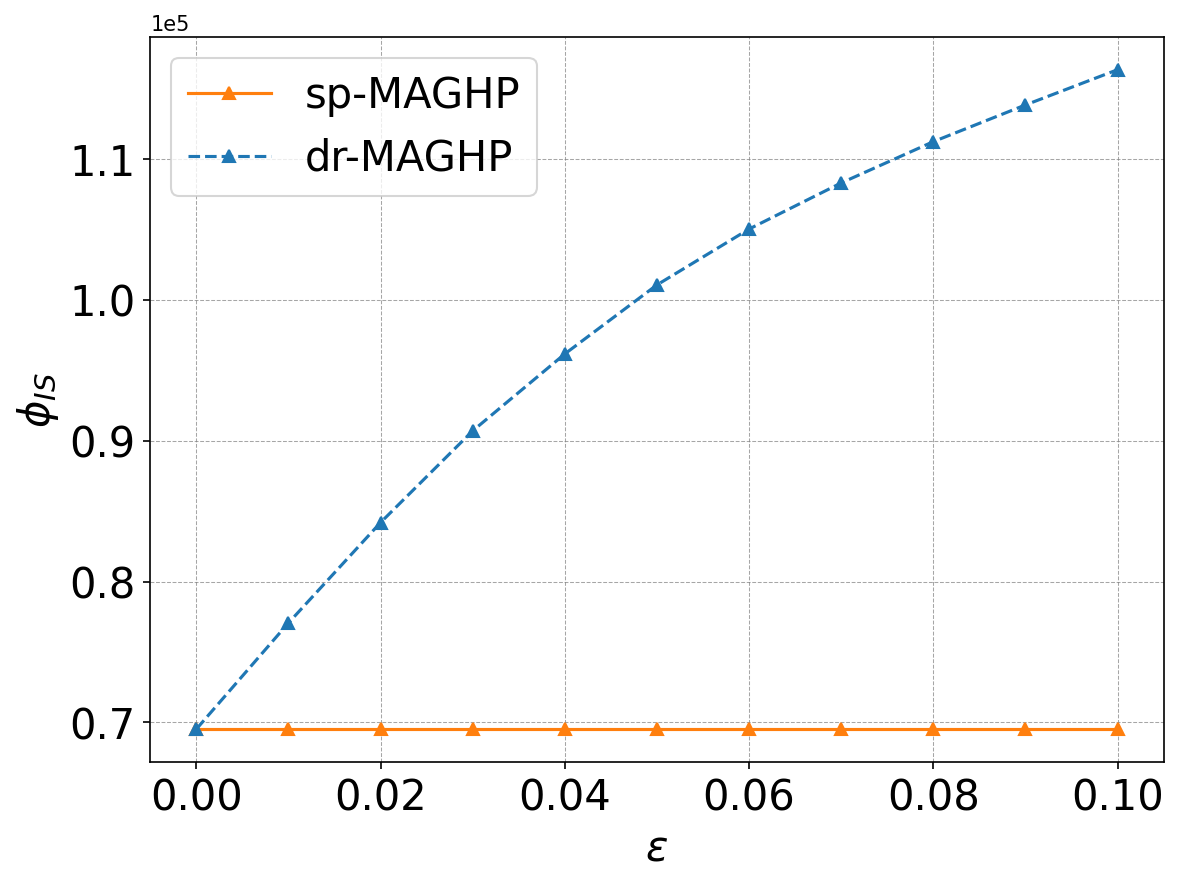}
        \caption{2019-12-26}
    \end{subfigure}

    \caption{\emph{in-sample} performance of \textsc{dr-MAGHP} with different $\epsilon$ values under overestimation scenario on all testing dates.}
    \label{fig:In_sample_costs}
\end{figure}

\begin{figure}
    \centering
    \begin{subfigure}[b]{0.49\textwidth}
        \includegraphics[width=\textwidth]{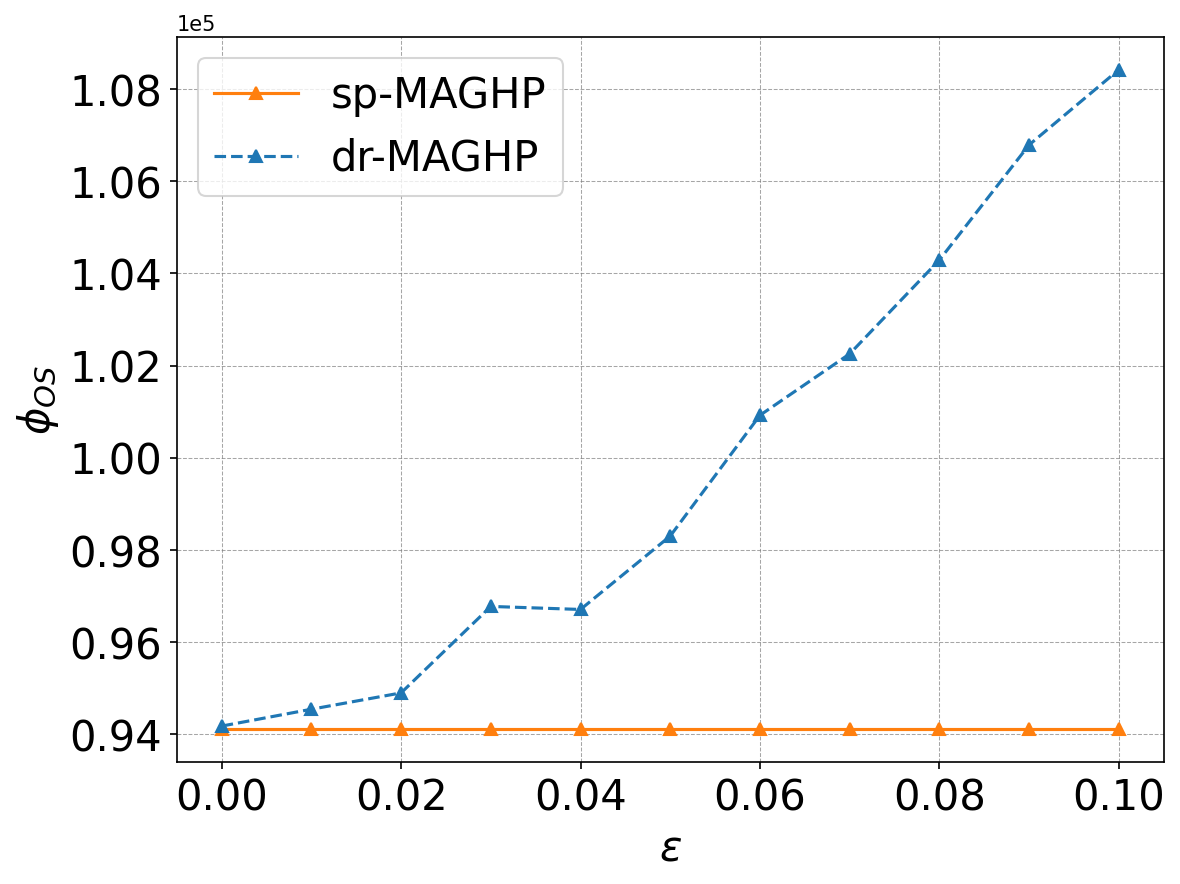}
        \caption{rl = 0.1}
    \end{subfigure}
    \hfill
    \begin{subfigure}[b]{0.49\textwidth}
        \includegraphics[width=\textwidth]{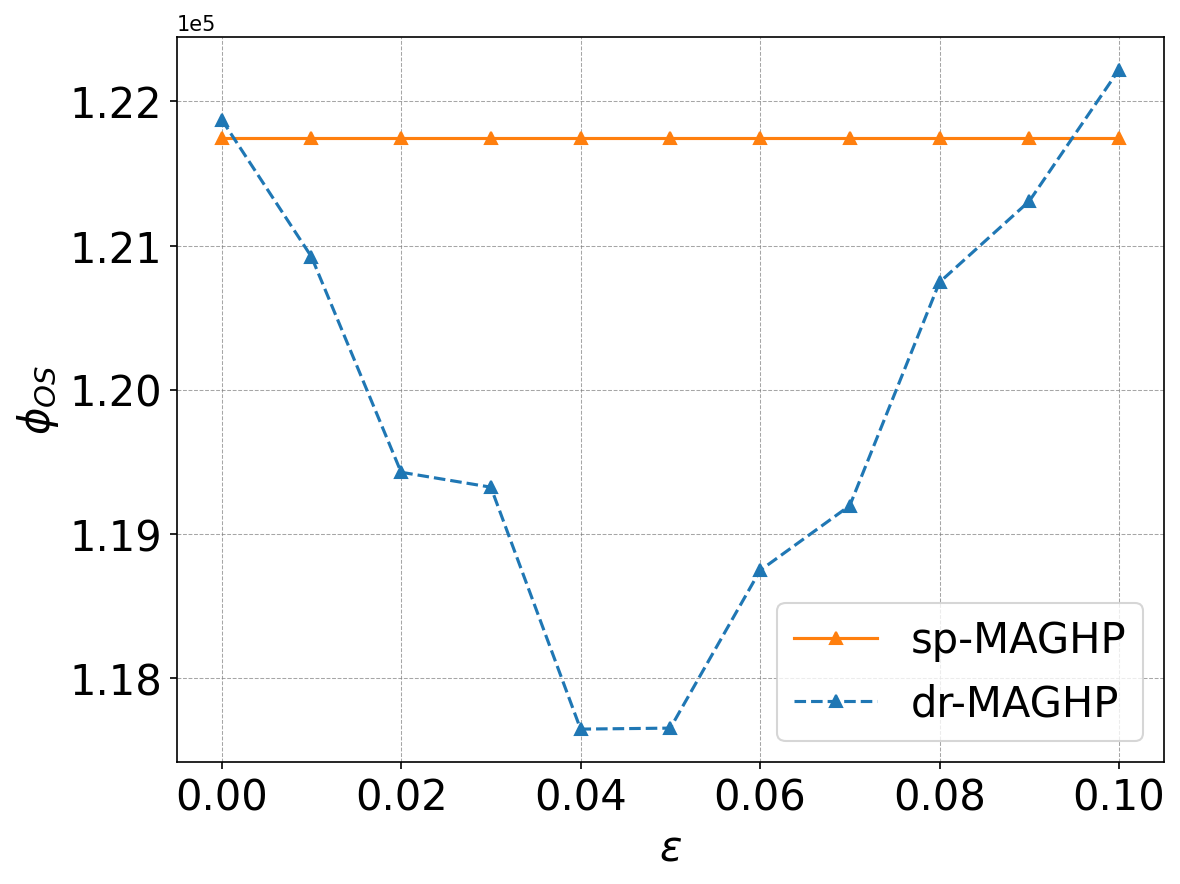}
        \caption{rl = 0.2}
    \end{subfigure}

    \vspace{0.7em}

    \begin{subfigure}[b]{0.49\textwidth}
        \includegraphics[width=\textwidth]{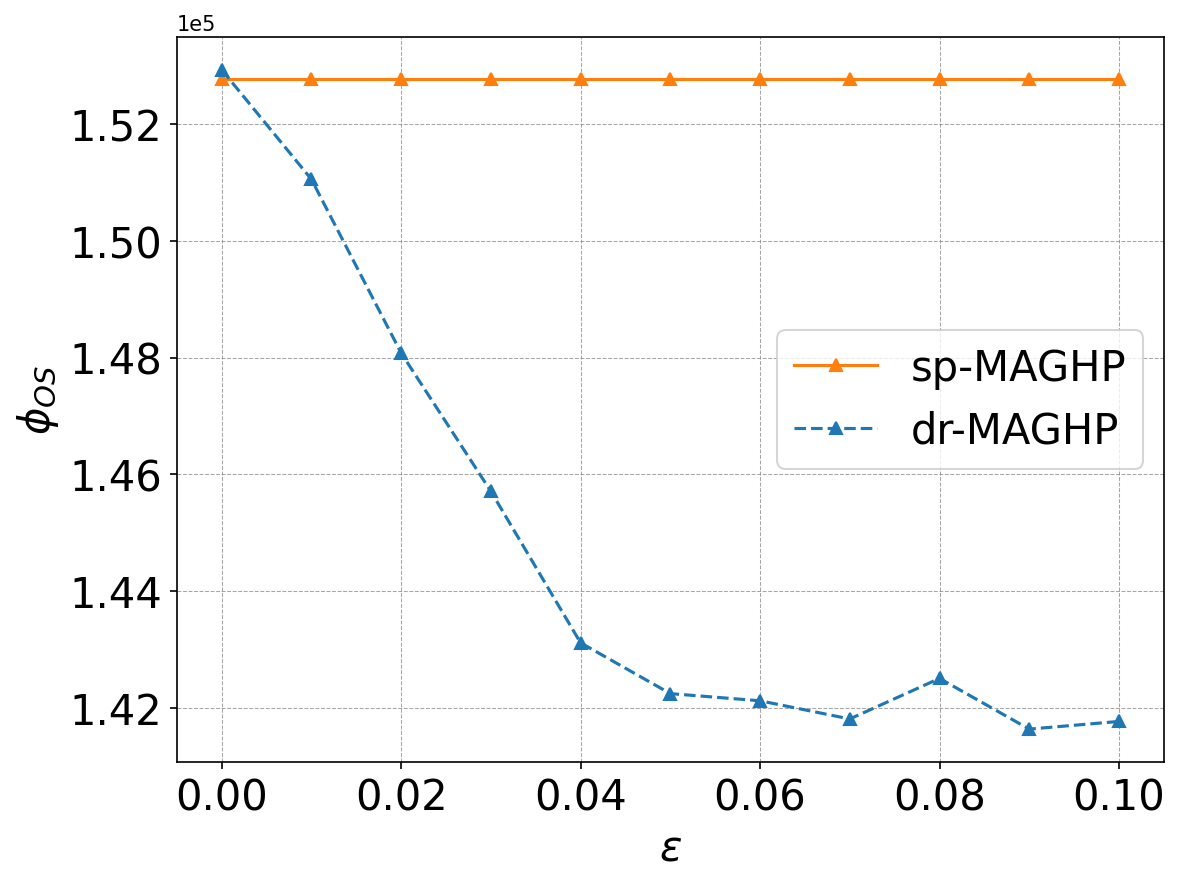}
        \caption{rl = 0.3}
    \end{subfigure}
    \hfill
    \begin{subfigure}[b]{0.49\textwidth}
        \includegraphics[width=\textwidth]{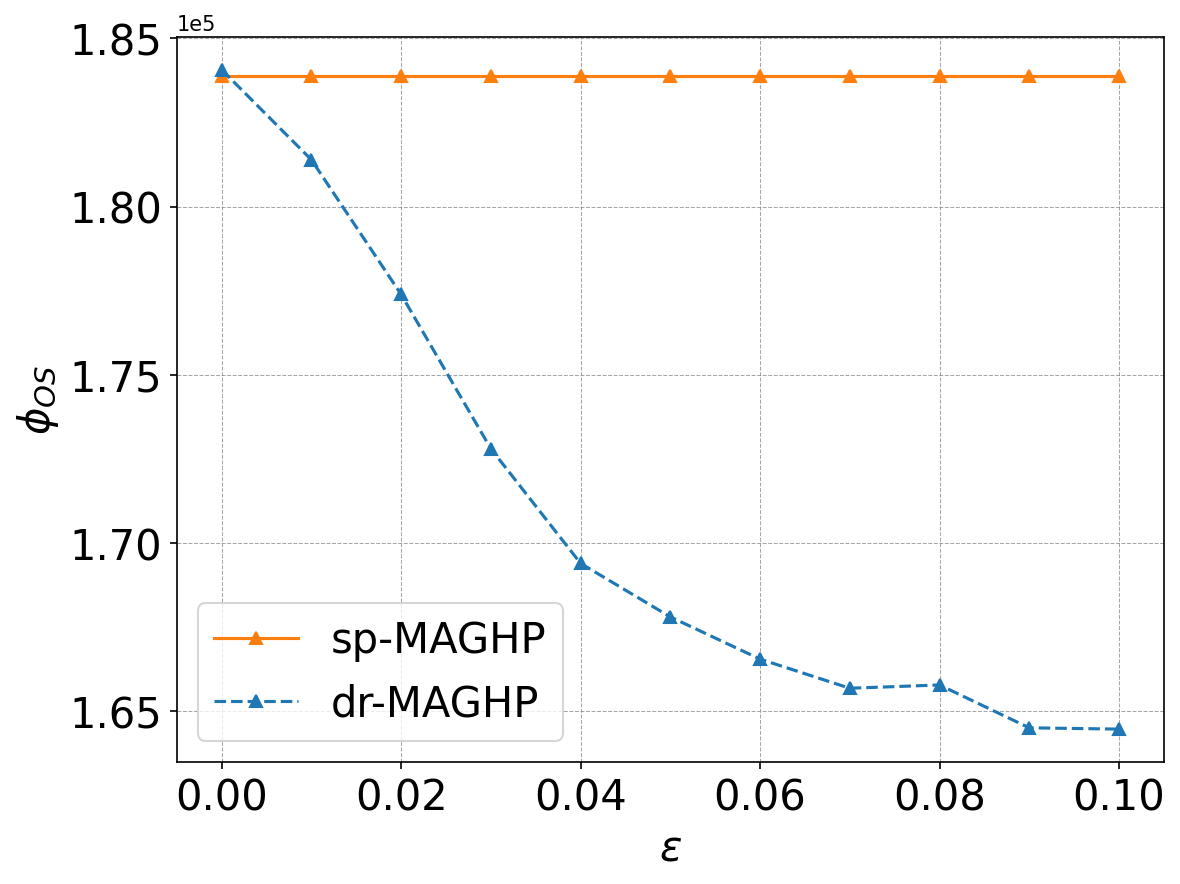}
        \caption{rl = 0.4}
    \end{subfigure}

    \vspace{0.7em}

    \begin{subfigure}[b]{0.5\textwidth}
        \centering
        \includegraphics[width=\textwidth]{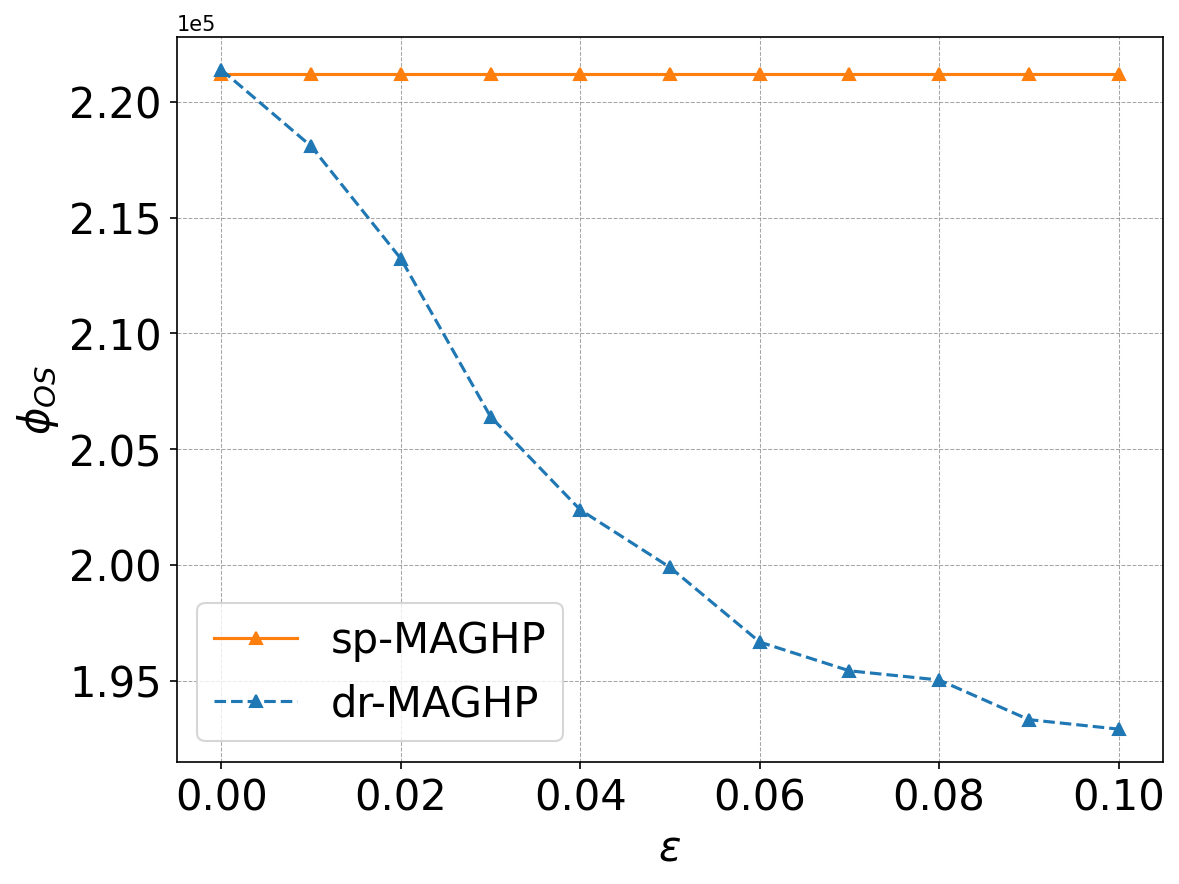}
        \caption{rl = 0.5}
    \end{subfigure}

    \caption{\emph{out-of-sample} performance of \textsc{sp-MAGHP} and \textsc{dr-MAGHP} with different $\epsilon$ values under overestimation scenario on 2019-11-13.}
    \label{fig:dr_cost_eps_overest_20191113}
\end{figure}

\begin{figure}
    \centering
    \begin{subfigure}[b]{0.49\textwidth}
        \includegraphics[width=\textwidth]{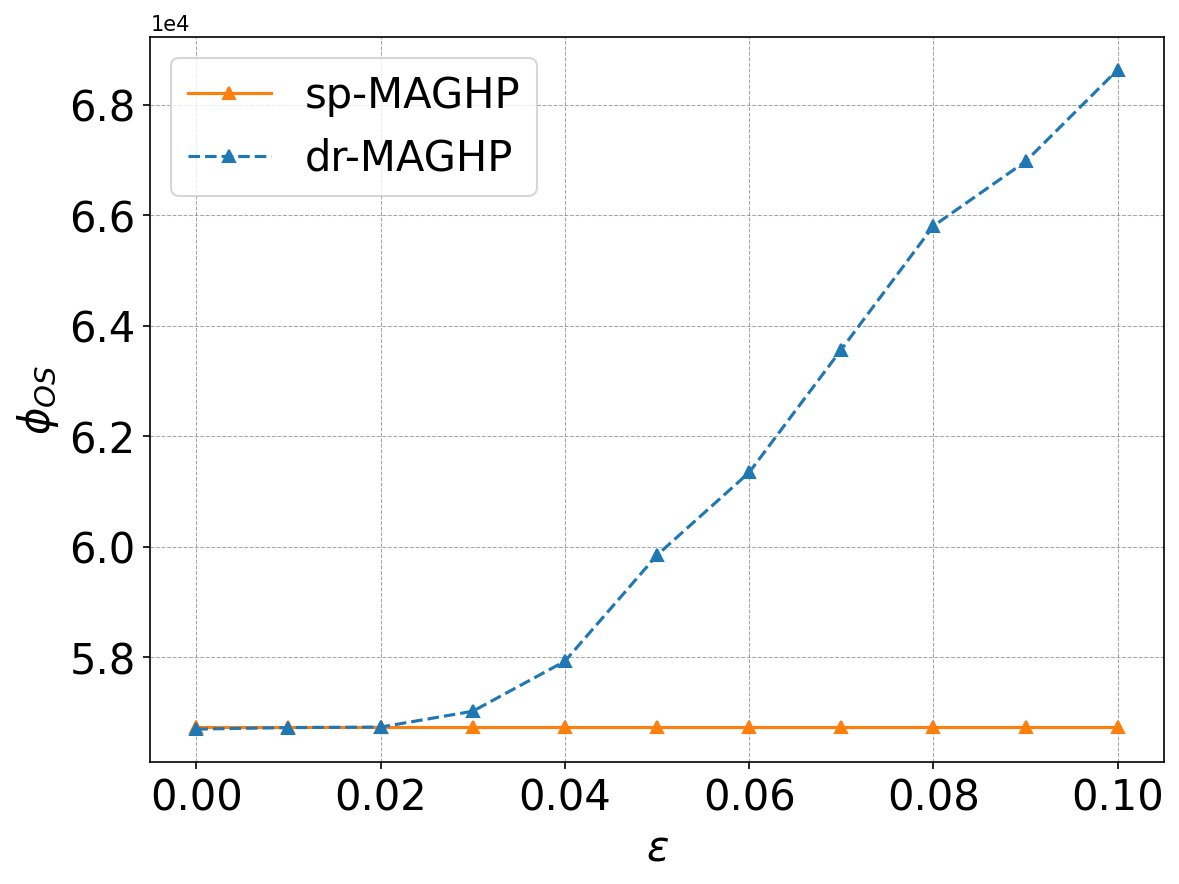}
        \caption{rl = 0.1}
    \end{subfigure}
    \hfill
    \begin{subfigure}[b]{0.49\textwidth}
        \includegraphics[width=\textwidth]{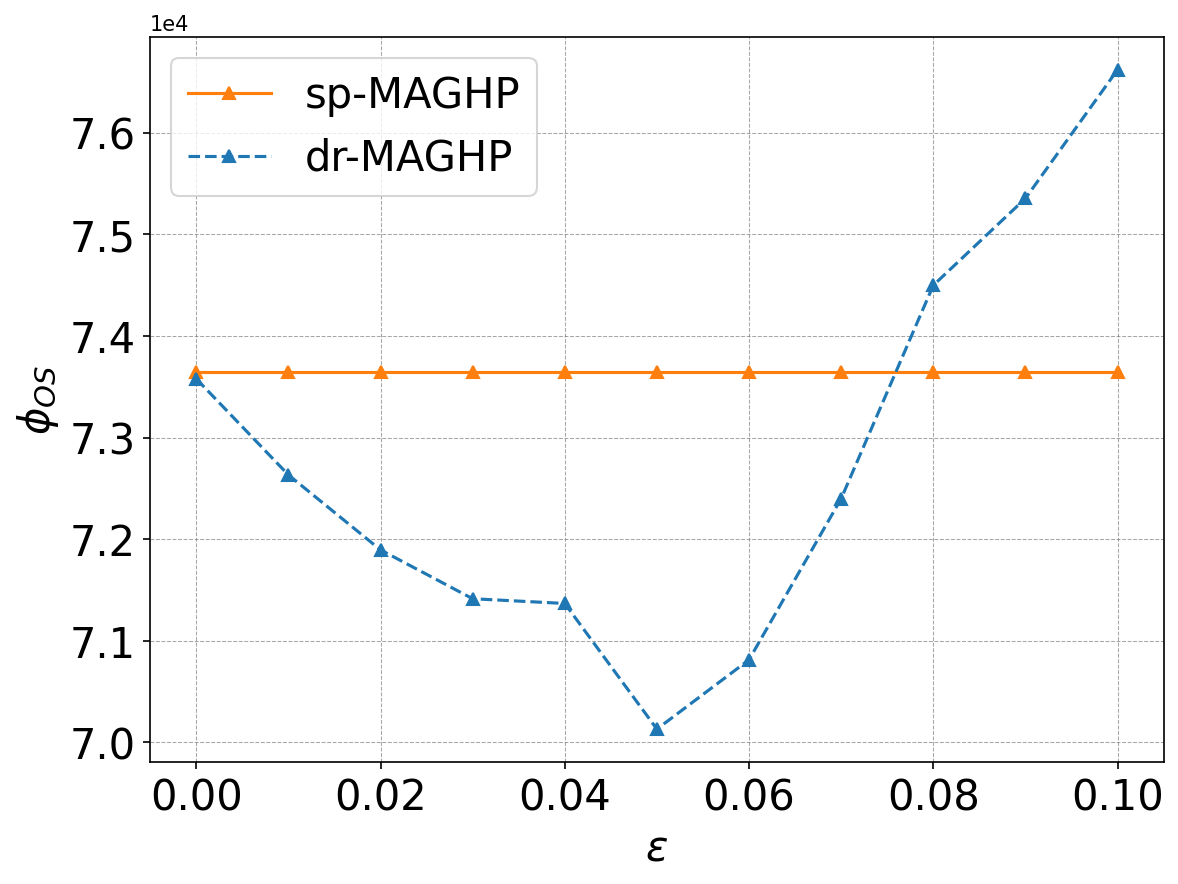}
        \caption{rl = 0.2}
    \end{subfigure}

    \vspace{0.7em}

    \begin{subfigure}[b]{0.49\textwidth}
        \includegraphics[width=\textwidth]{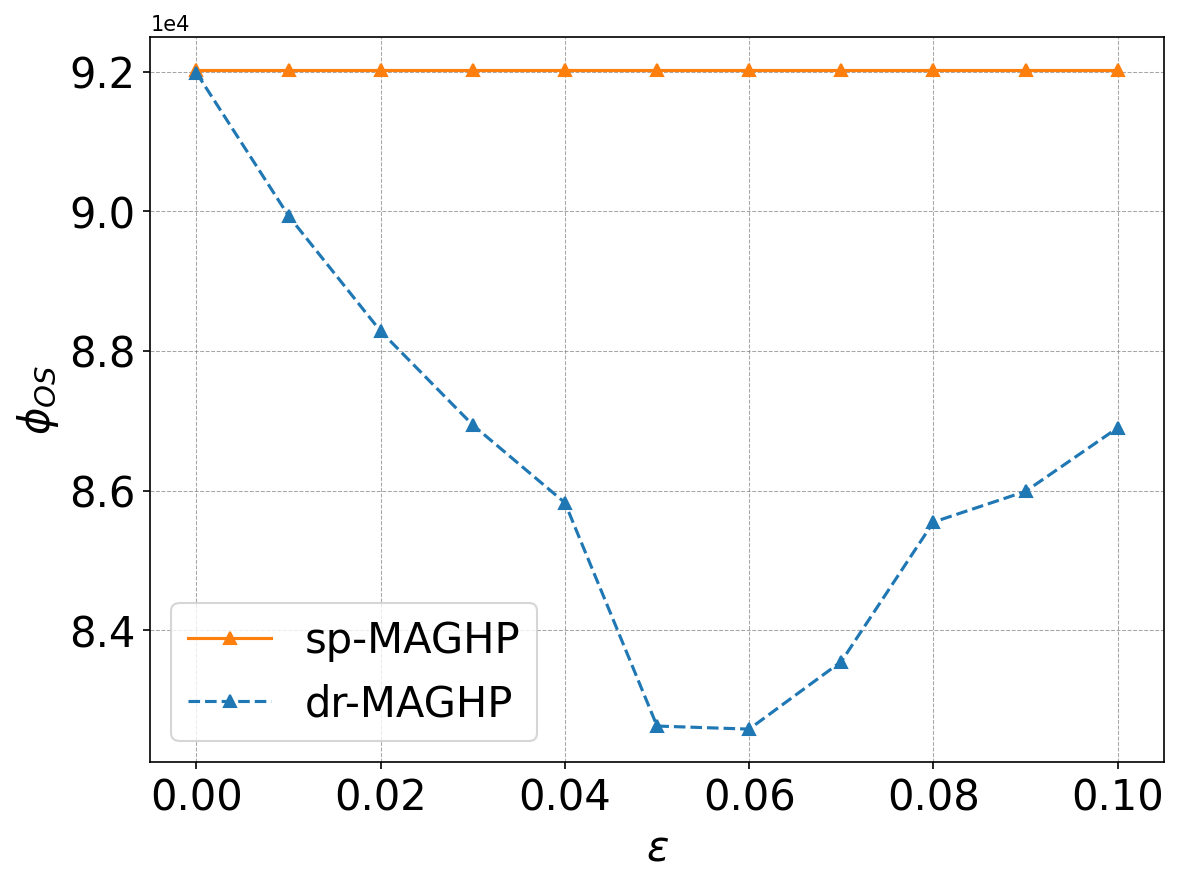}
        \caption{rl = 0.3}
    \end{subfigure}
    \hfill
    \begin{subfigure}[b]{0.49\textwidth}
        \includegraphics[width=\textwidth]{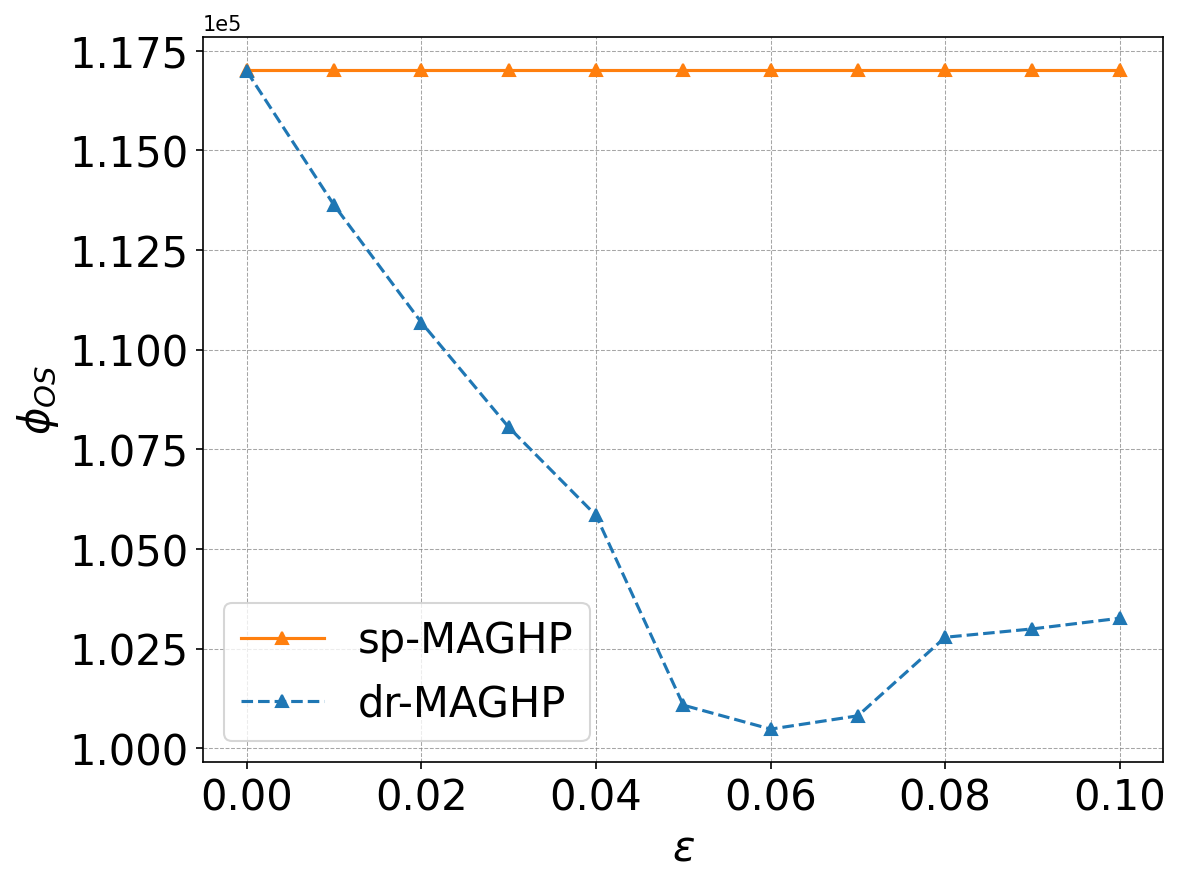}
        \caption{rl = 0.4}
    \end{subfigure}

    \vspace{0.7em}

    \begin{subfigure}[b]{0.5\textwidth}
        \centering
        \includegraphics[width=\textwidth]{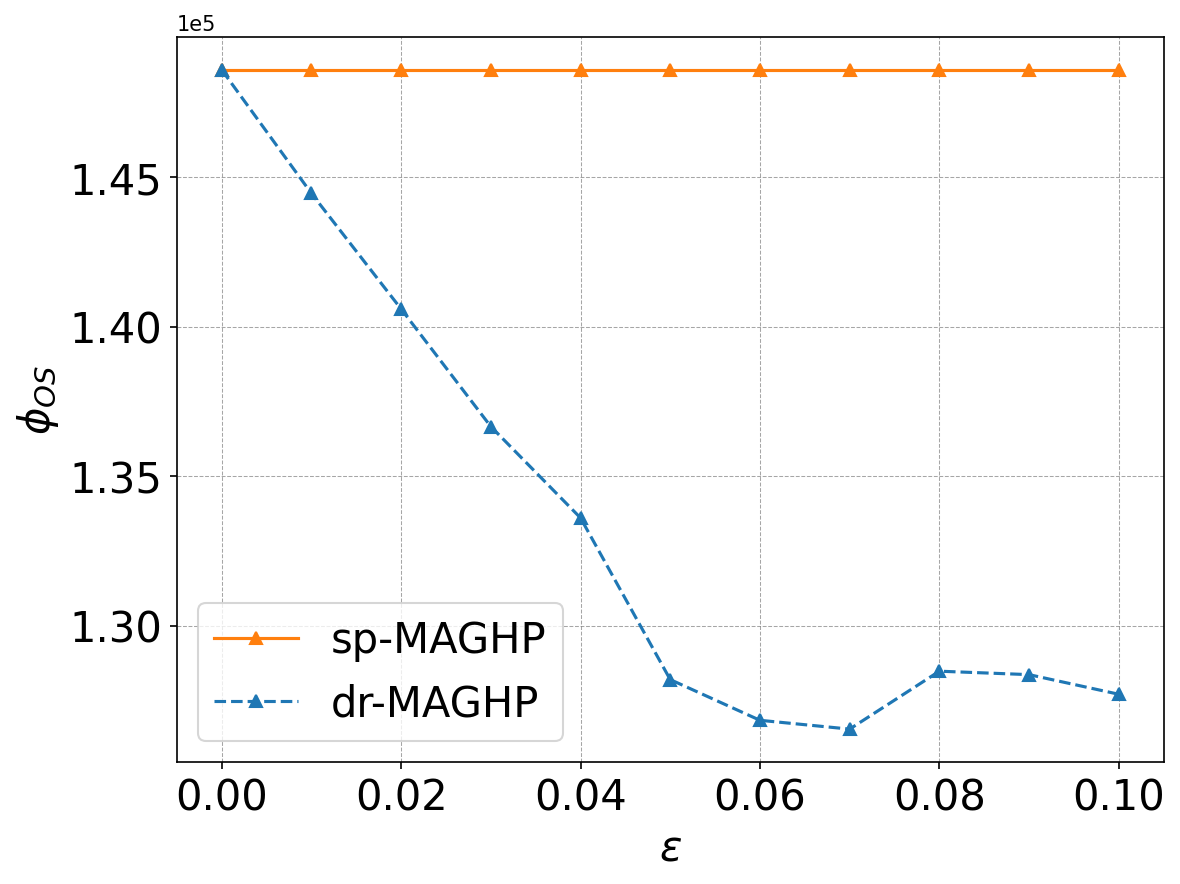}
        \caption{rl = 0.5}
    \end{subfigure}

    \caption{\emph{out-of-sample} performance of \textsc{sp-MAGHP} and \textsc{dr-MAGHP} with different $\epsilon$ values under overestimation scenario on 2019-12-25.}
    \label{fig:dr_cost_eps_overest_20191225}
\end{figure}

\begin{figure}
    \centering
    \begin{subfigure}[b]{0.49\textwidth}
        \includegraphics[width=\textwidth]{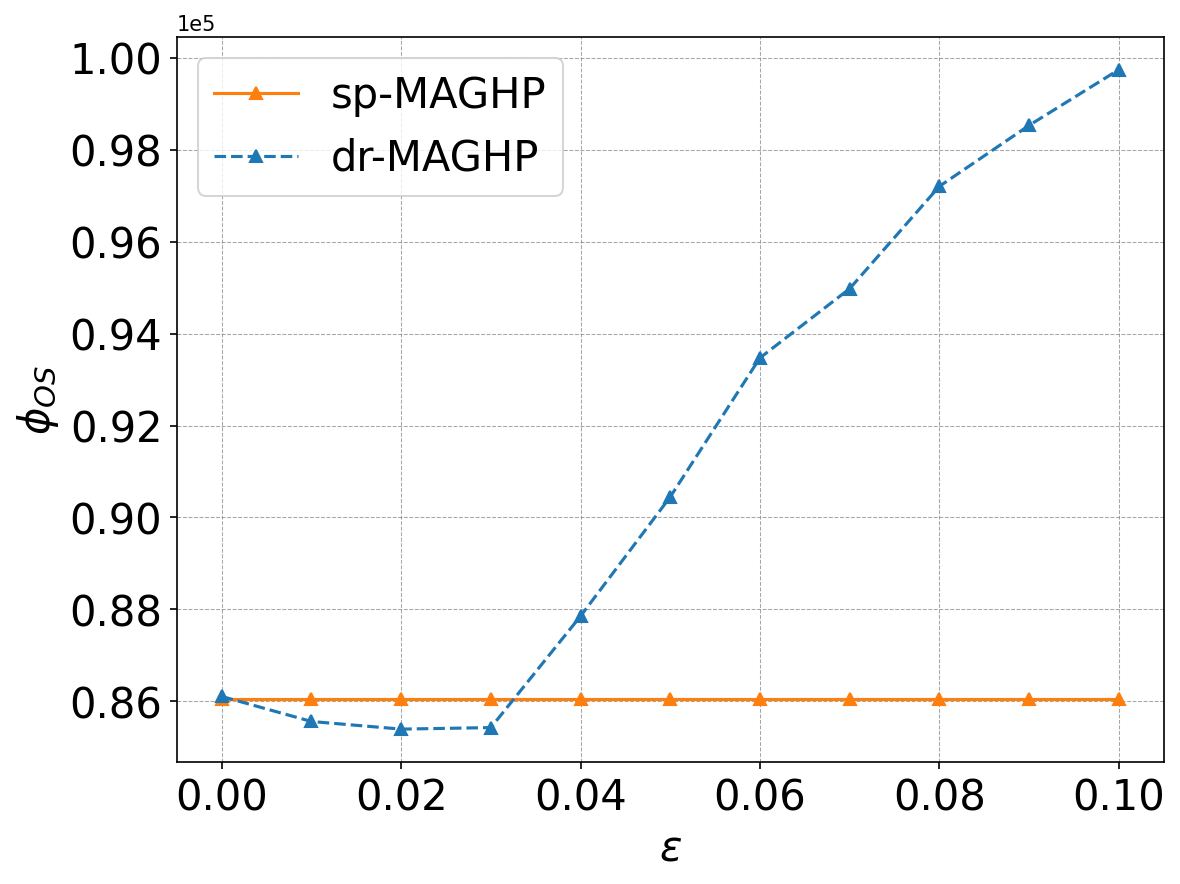}
        \caption{rl = 0.1}
    \end{subfigure}
    \hfill
    \begin{subfigure}[b]{0.49\textwidth}
        \includegraphics[width=\textwidth]{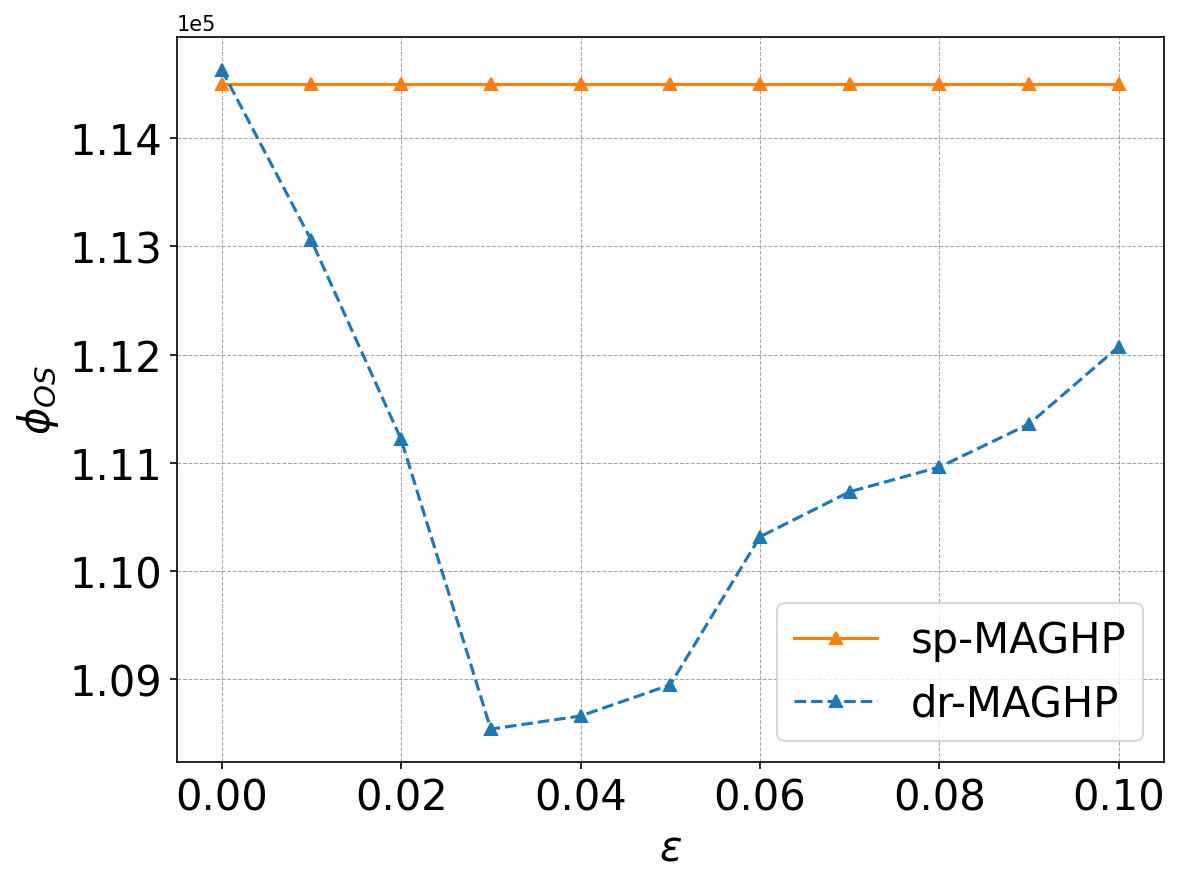}
        \caption{rl = 0.2}
    \end{subfigure}

    \vspace{0.7em}

    \begin{subfigure}[b]{0.49\textwidth}
        \includegraphics[width=\textwidth]{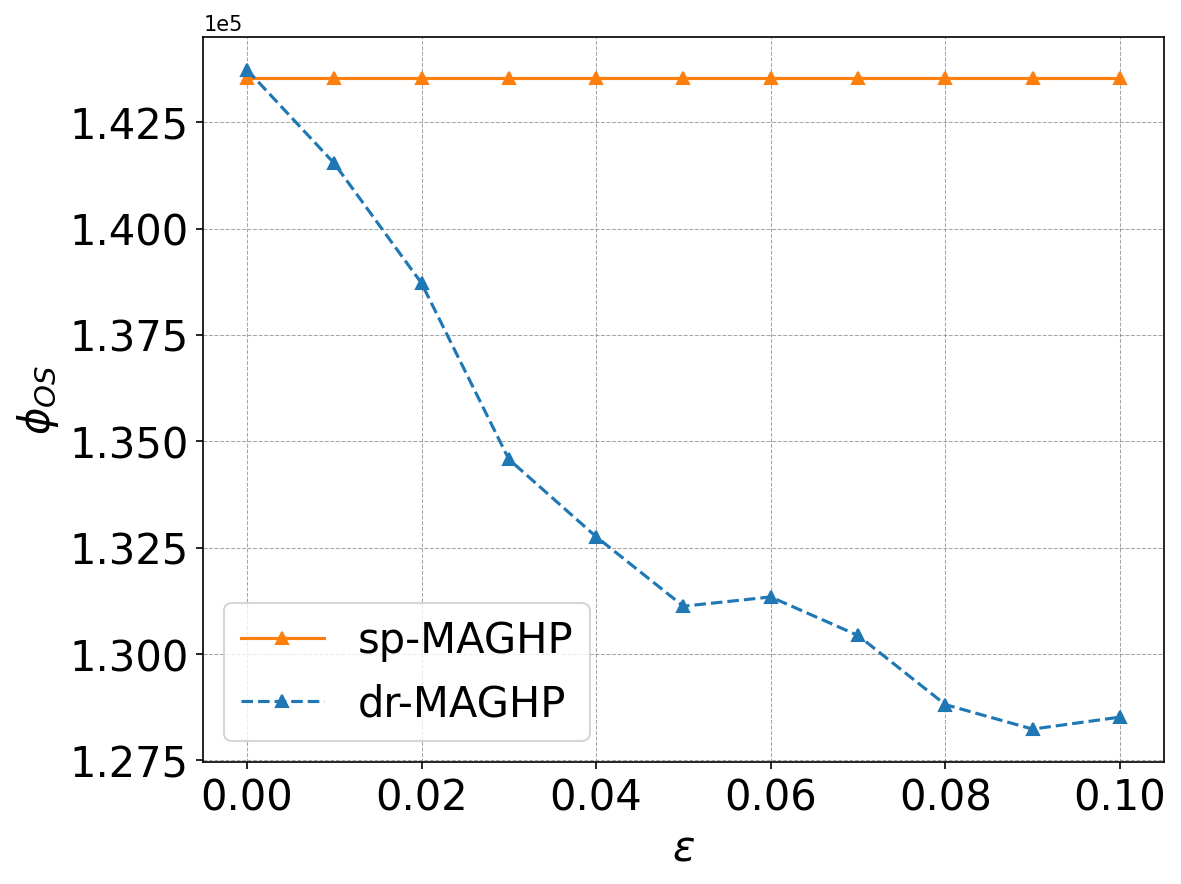}
        \caption{rl = 0.3}
    \end{subfigure}
    \hfill
    \begin{subfigure}[b]{0.49\textwidth}
        \includegraphics[width=\textwidth]{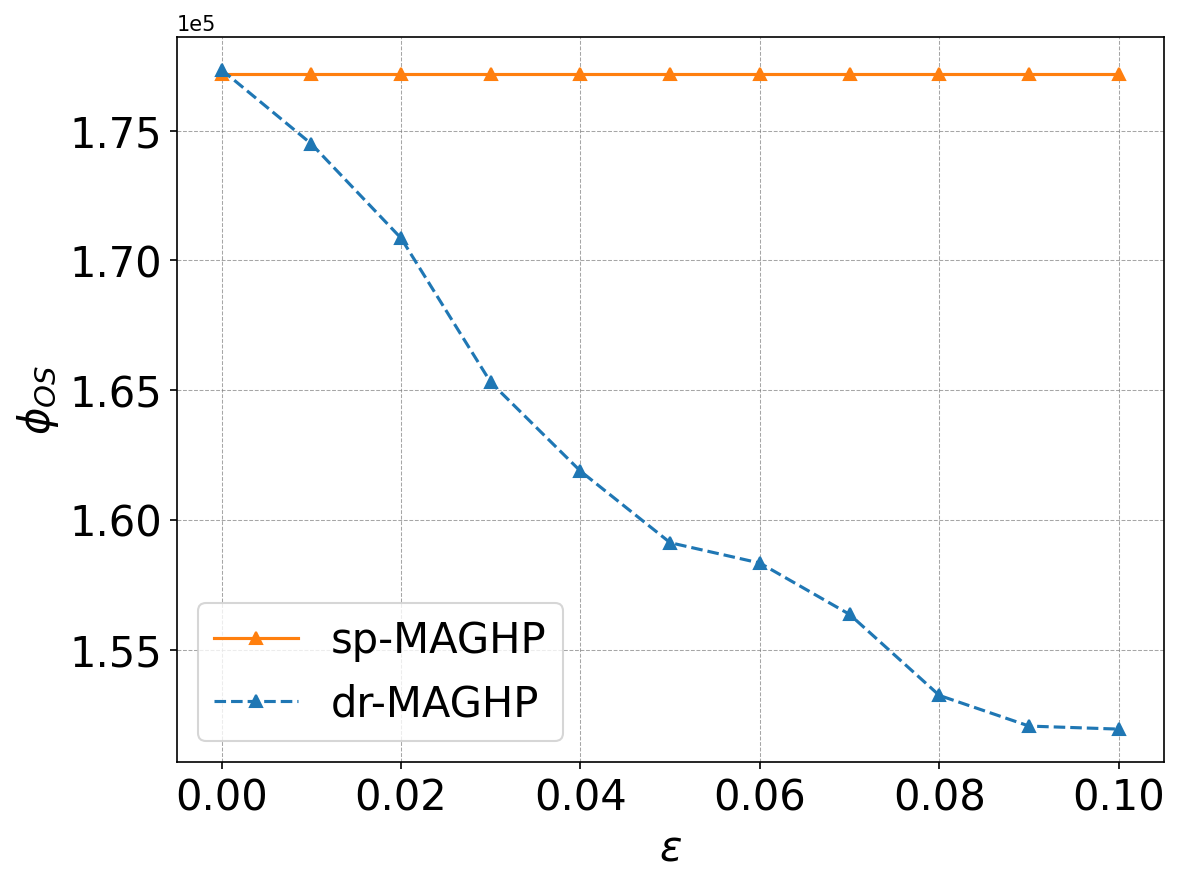}
        \caption{rl = 0.4}
    \end{subfigure}

    \vspace{0.7em}

    \begin{subfigure}[b]{0.5\textwidth}
        \centering
        \includegraphics[width=\textwidth]{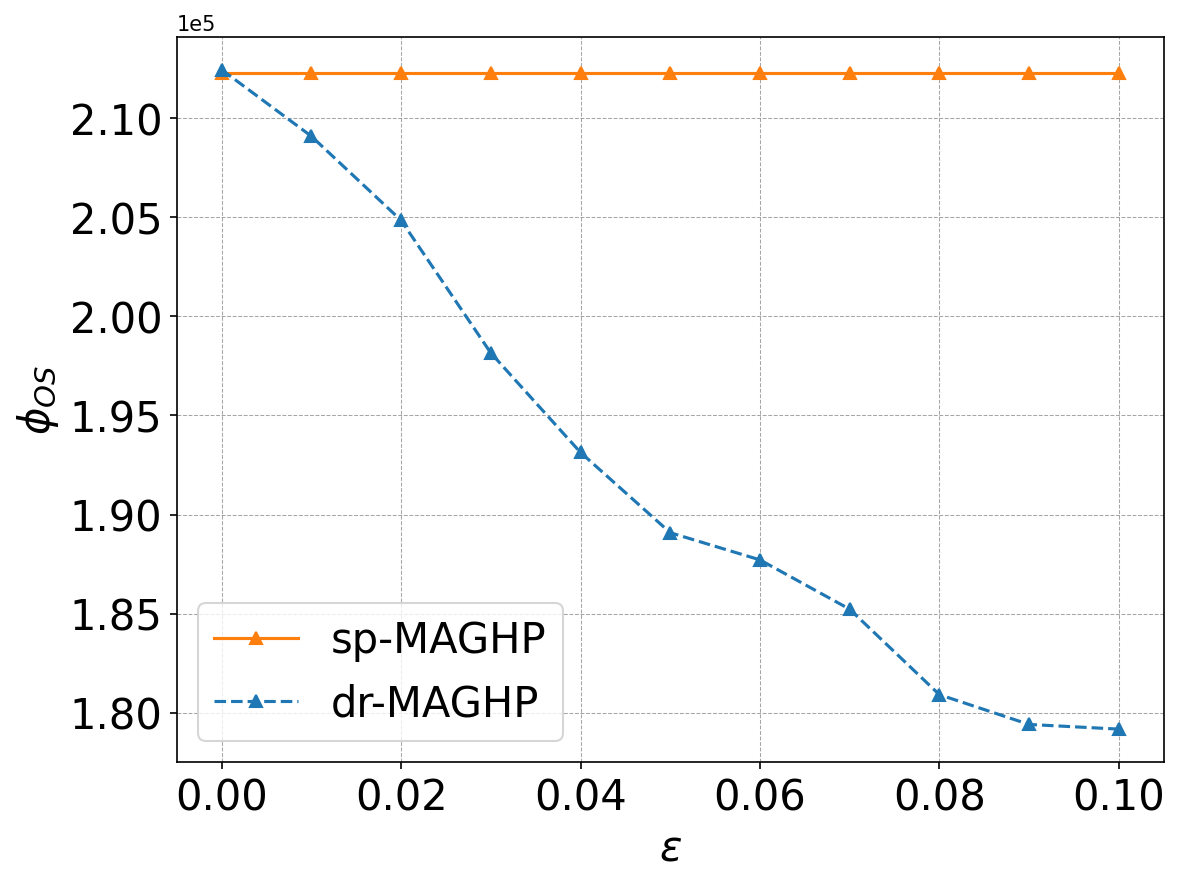}
        \caption{rl = 0.5}
    \end{subfigure}

    \caption{\emph{out-of-sample} performance of \textsc{sp-MAGHP} and \textsc{dr-MAGHP} with different $\epsilon$ values under overestimation scenario on 2019-12-26.}
    \label{fig:dr_cost_eps_overest_20191226}
\end{figure}

\begin{figure}
    \centering
    \begin{subfigure}[b]{0.49\textwidth}
        \includegraphics[width=\textwidth]{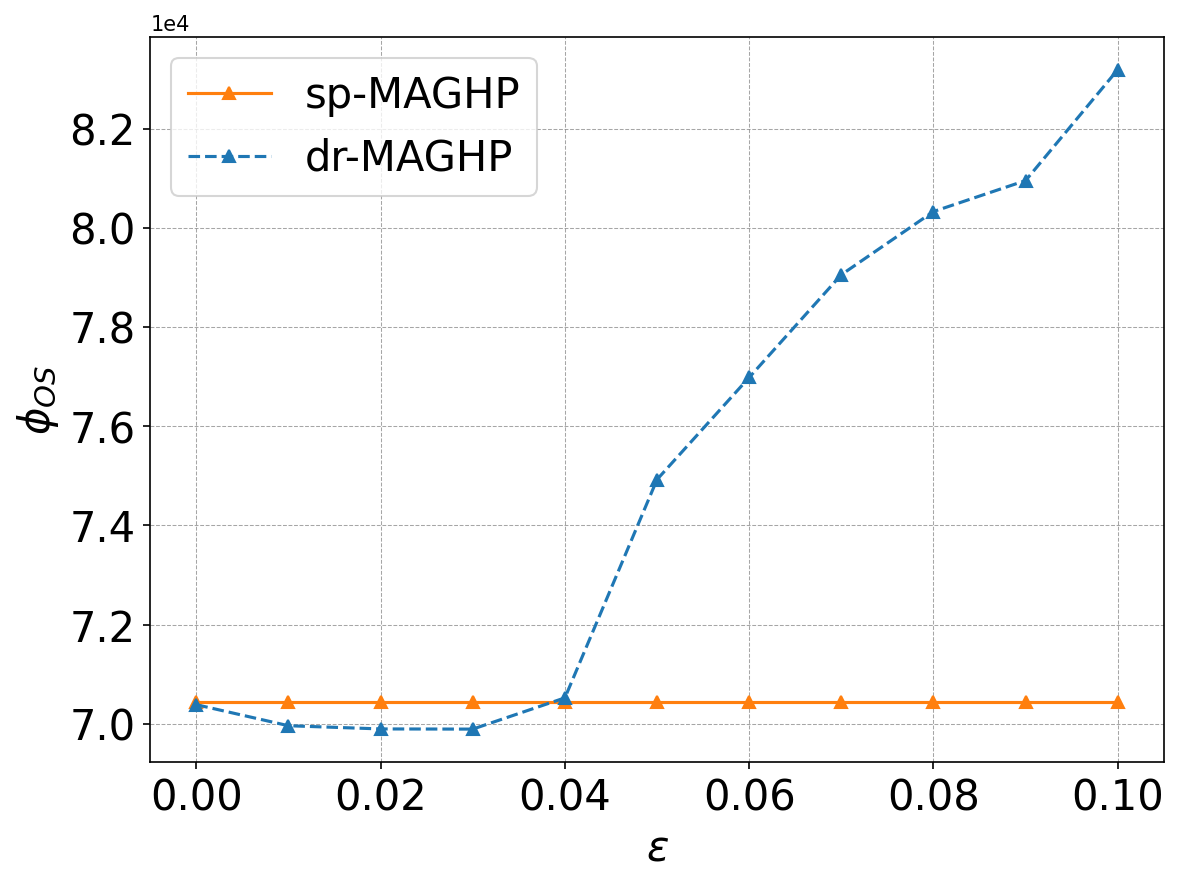}
        \caption{rl = 0.1}
    \end{subfigure}
    \hfill
    \begin{subfigure}[b]{0.49\textwidth}
        \includegraphics[width=\textwidth]{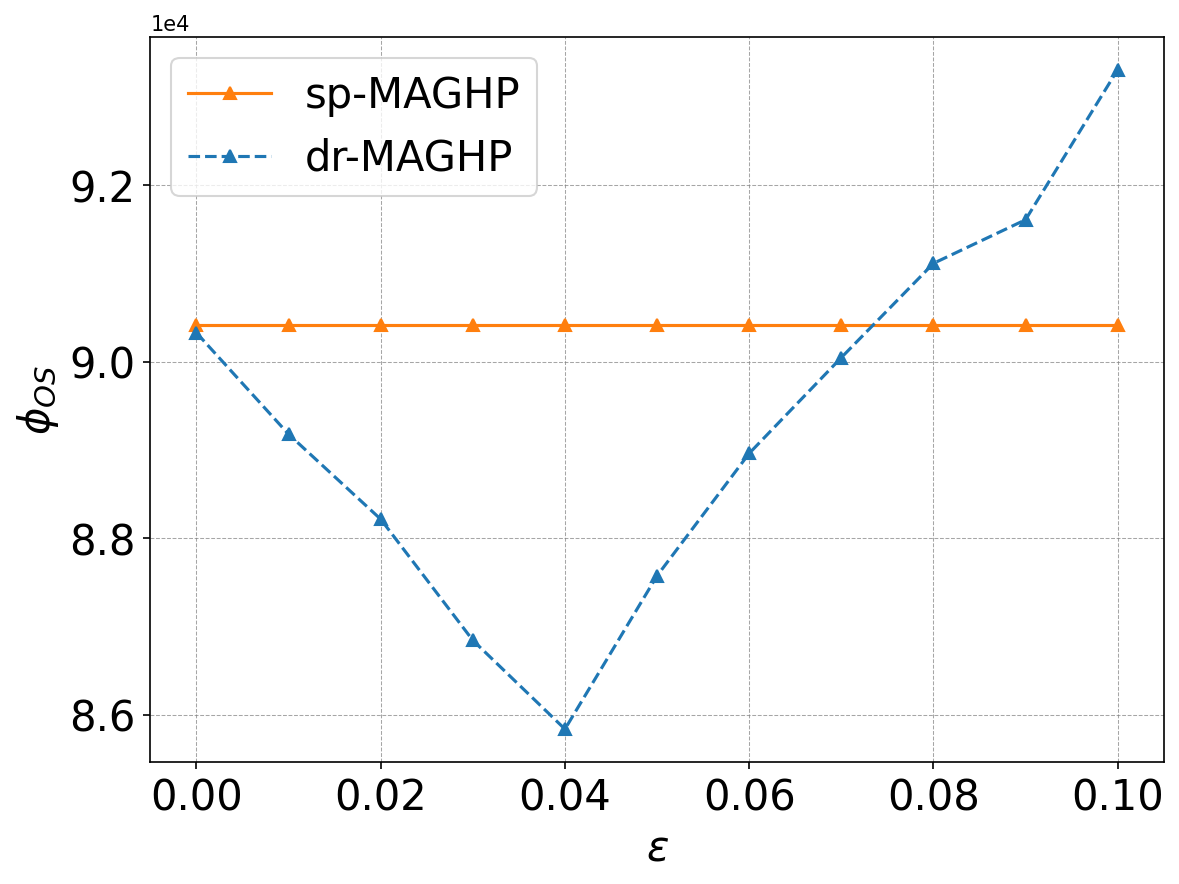}
        \caption{rl = 0.2}
    \end{subfigure}

    \vspace{0.7em}

    \begin{subfigure}[b]{0.49\textwidth}
        \includegraphics[width=\textwidth]{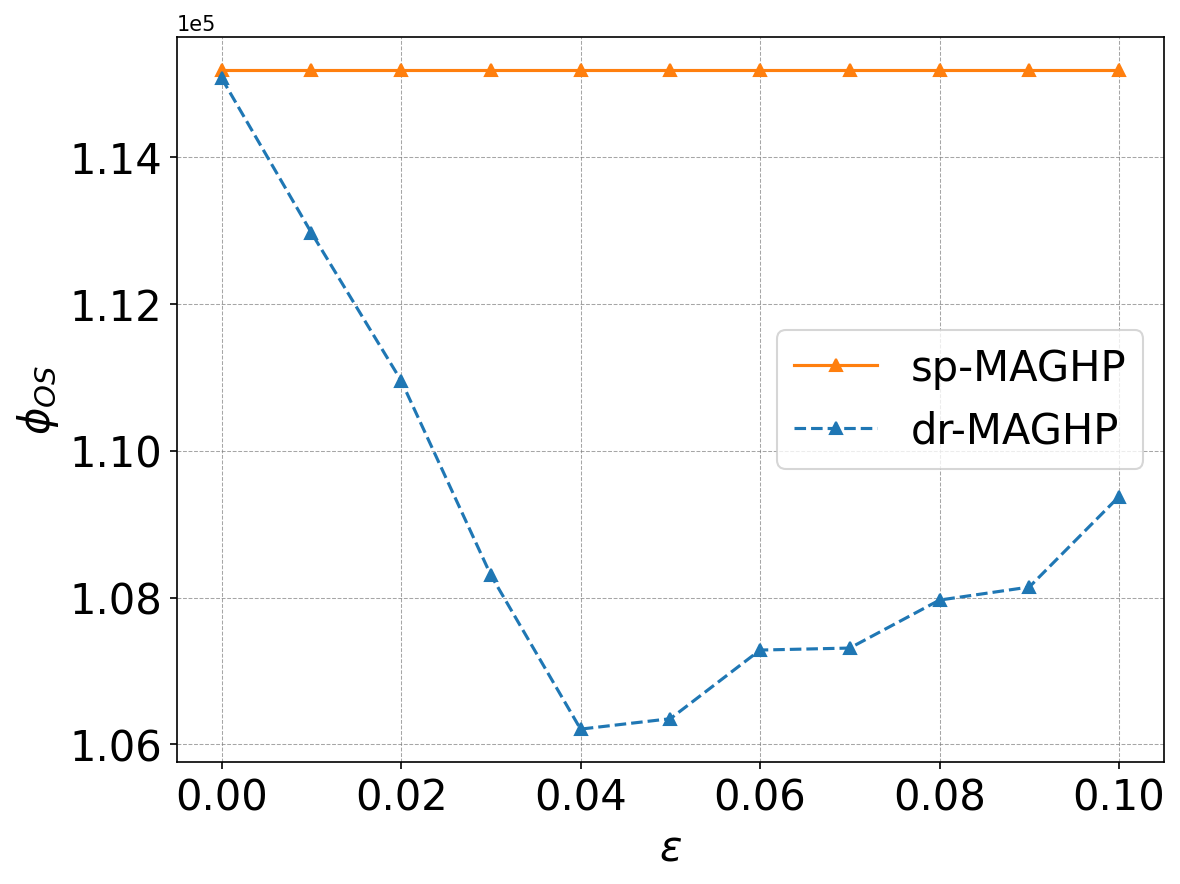}
        \caption{rl = 0.3}
    \end{subfigure}
    \hfill
    \begin{subfigure}[b]{0.49\textwidth}
        \includegraphics[width=\textwidth]{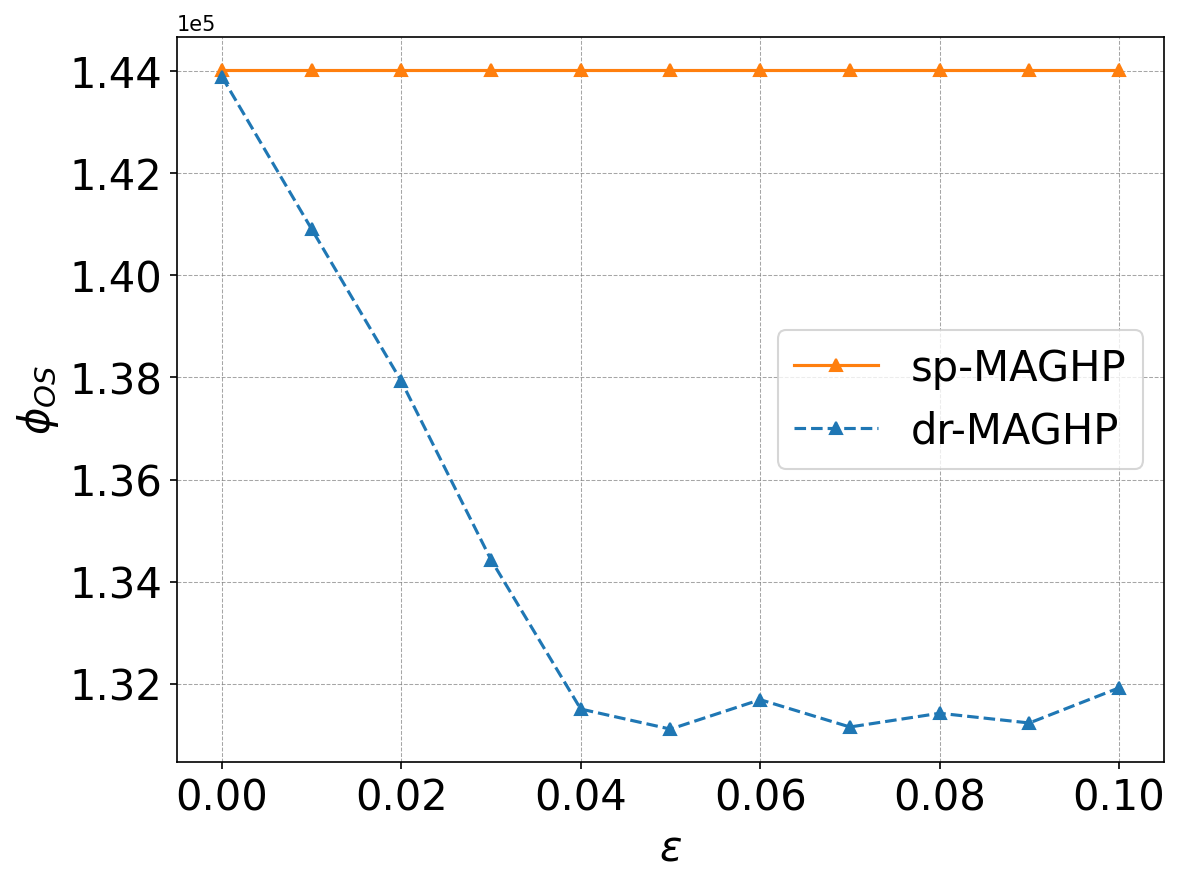}
        \caption{rl = 0.4}
    \end{subfigure}

    \vspace{0.7em}

    \begin{subfigure}[b]{0.5\textwidth}
        \centering
        \includegraphics[width=\textwidth]{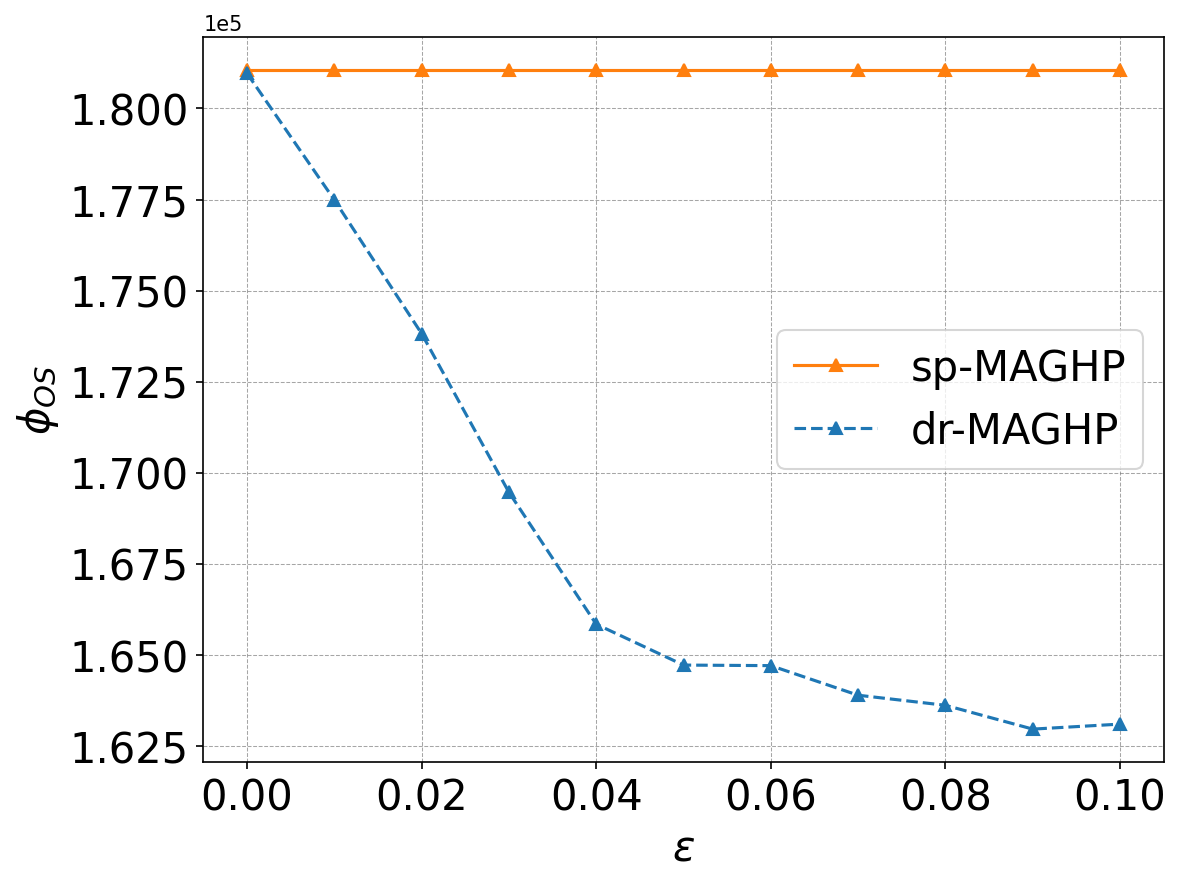}
        \caption{rl = 0.5}
    \end{subfigure}

    \caption{\emph{out-of-sample} performance of \textsc{sp-MAGHP} and \textsc{dr-MAGHP} with different $\epsilon$ values under overestimation scenario on 2019-05-04.}
    \label{fig:dr_cost_eps_overest_20190504}
\end{figure}

\begin{figure}
    \centering
    \begin{subfigure}[b]{0.49\textwidth}
        \includegraphics[width=\textwidth]{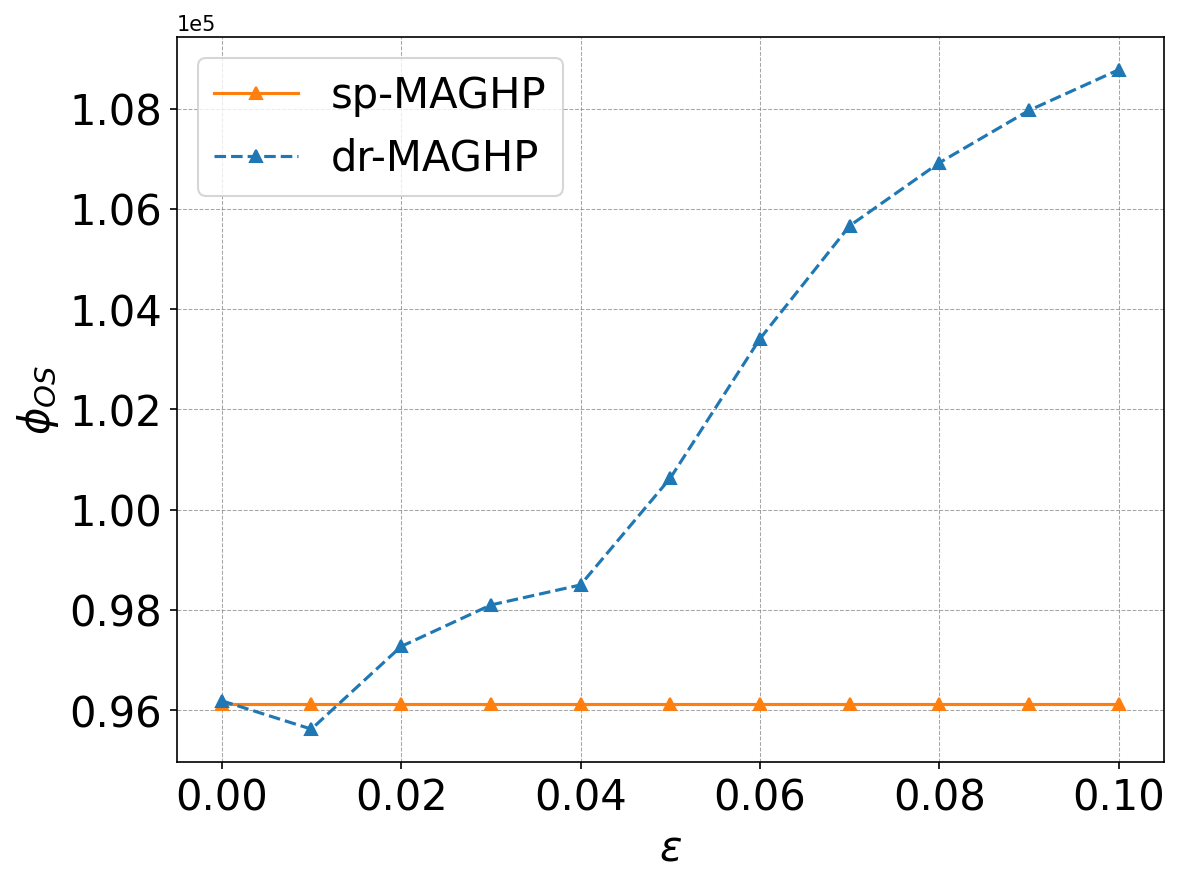}
        \caption{rl = 0.1}
    \end{subfigure}
    \hfill
    \begin{subfigure}[b]{0.49\textwidth}
        \includegraphics[width=\textwidth]{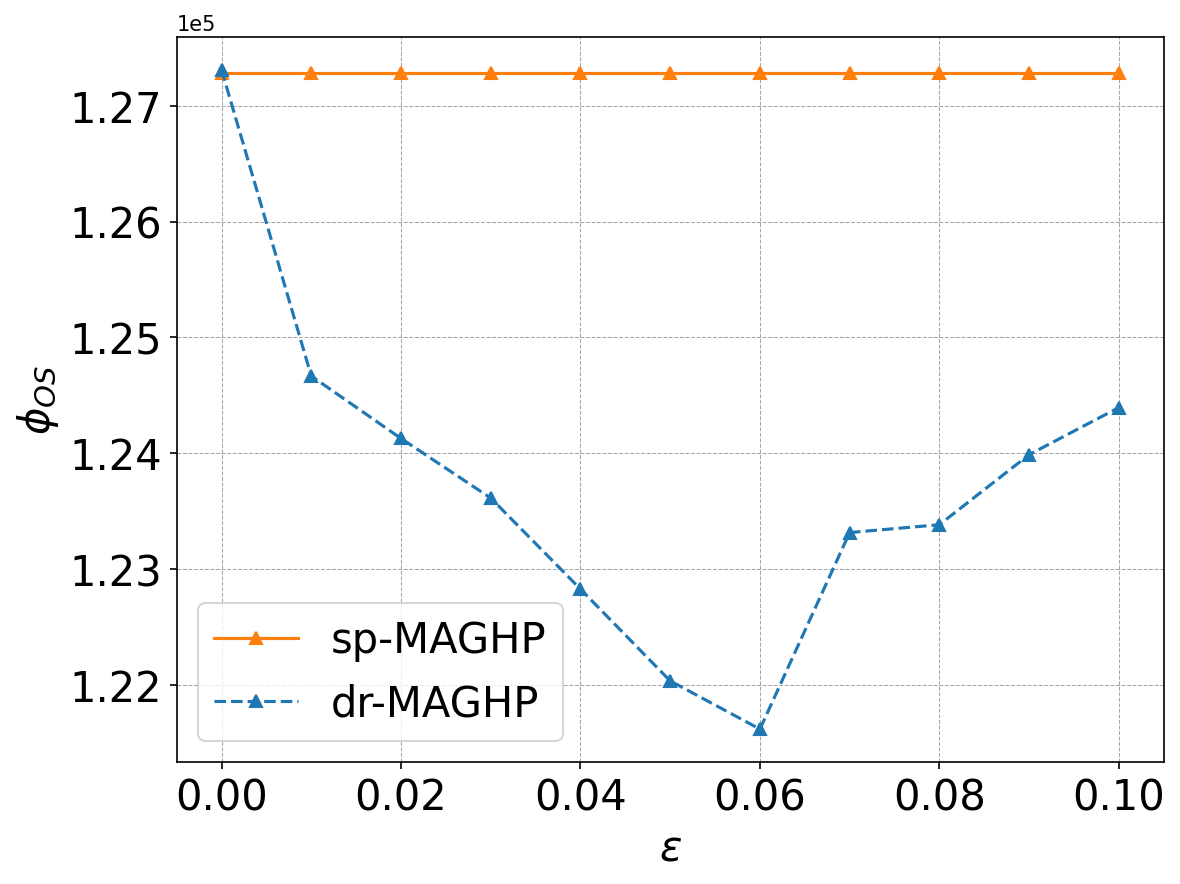}
        \caption{rl = 0.2}
    \end{subfigure}

    \vspace{0.7em}

    \begin{subfigure}[b]{0.49\textwidth}
        \includegraphics[width=\textwidth]{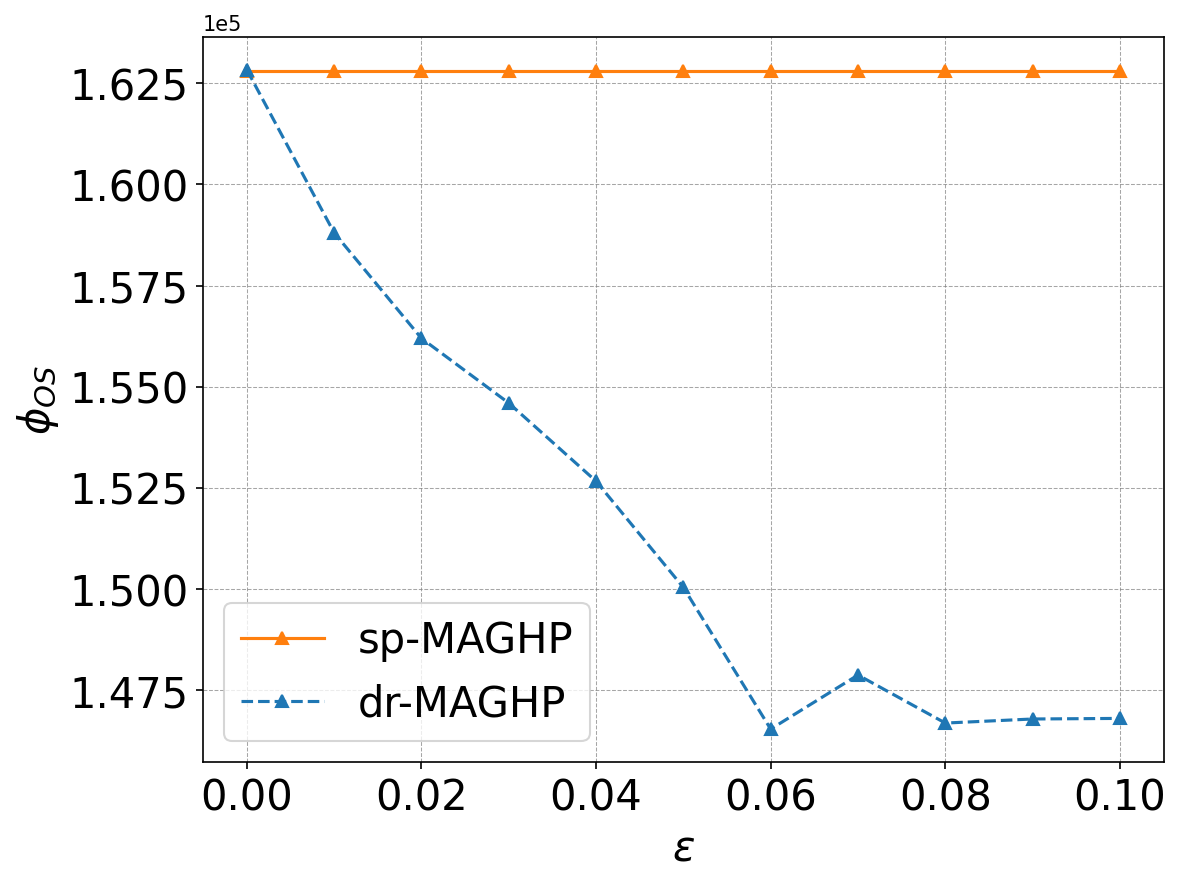}
        \caption{rl = 0.3}
    \end{subfigure}
    \hfill
    \begin{subfigure}[b]{0.49\textwidth}
        \includegraphics[width=\textwidth]{figures/2019_05_04/2019-05-04_overest_rl0.4.png}
        \caption{rl = 0.4}
    \end{subfigure}

    \vspace{0.7em}

    \begin{subfigure}[b]{0.5\textwidth}
        \centering
        \includegraphics[width=\textwidth]{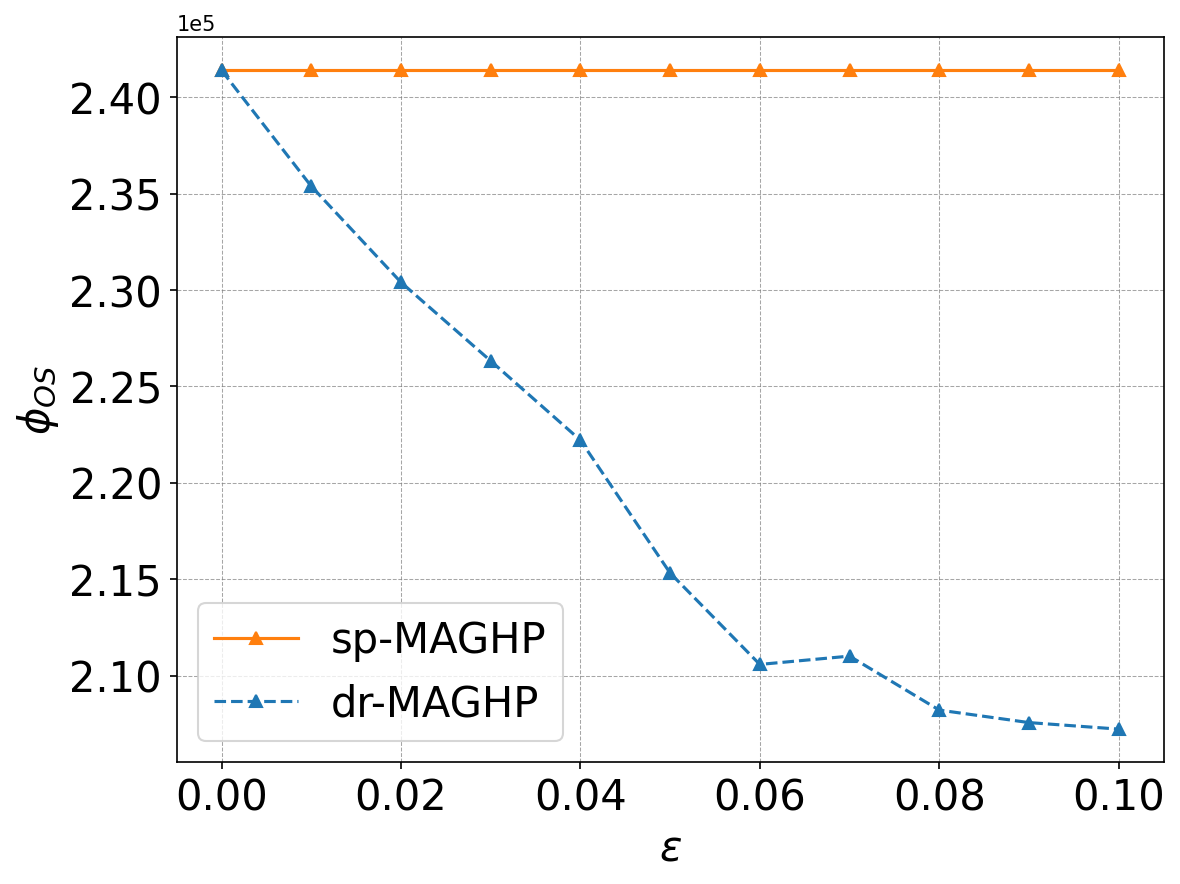}
        \caption{rl = 0.5}
    \end{subfigure}

    \caption{\emph{out-of-sample} performance of \textsc{sp-MAGHP} and \textsc{dr-MAGHP} with different $\epsilon$ values under overestimation scenario on 2019-05-30.}
    \label{fig:dr_cost_eps_overest_20190530}
\end{figure}

\begin{figure}
    \centering
    \begin{subfigure}[b]{0.49\textwidth}
        \includegraphics[width=\textwidth]{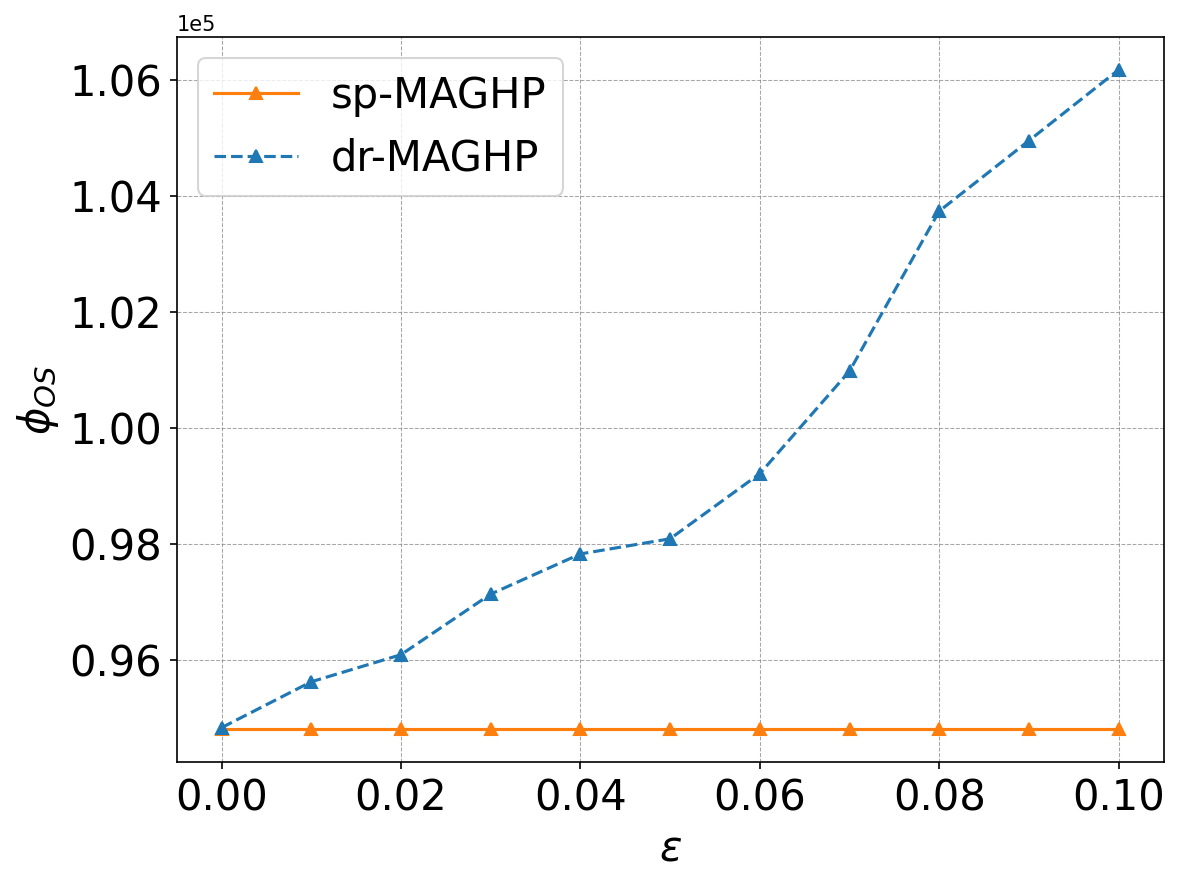}
        \caption{rl = 0.1}
    \end{subfigure}
    \hfill
    \begin{subfigure}[b]{0.49\textwidth}
        \includegraphics[width=\textwidth]{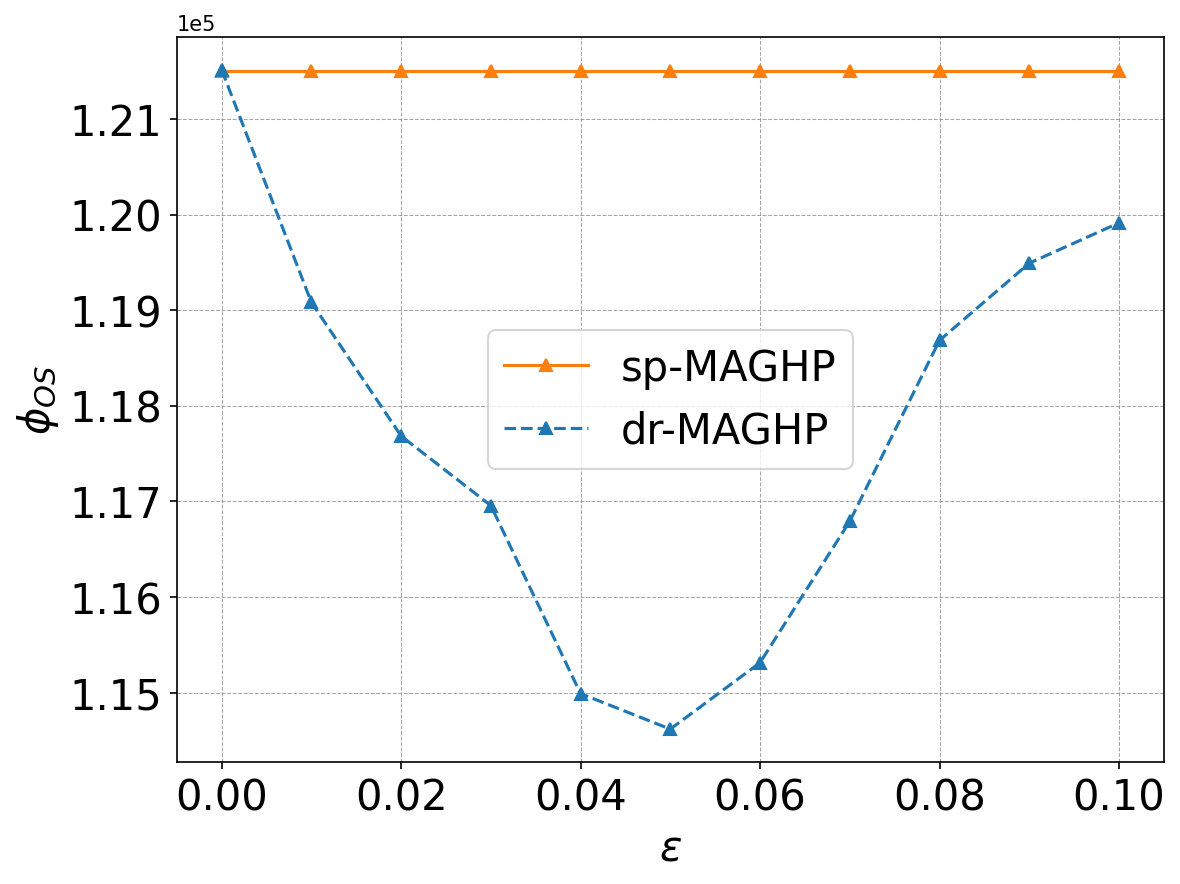}
        \caption{rl = 0.2}
    \end{subfigure}

    \vspace{0.7em}

    \begin{subfigure}[b]{0.49\textwidth}
        \includegraphics[width=\textwidth]{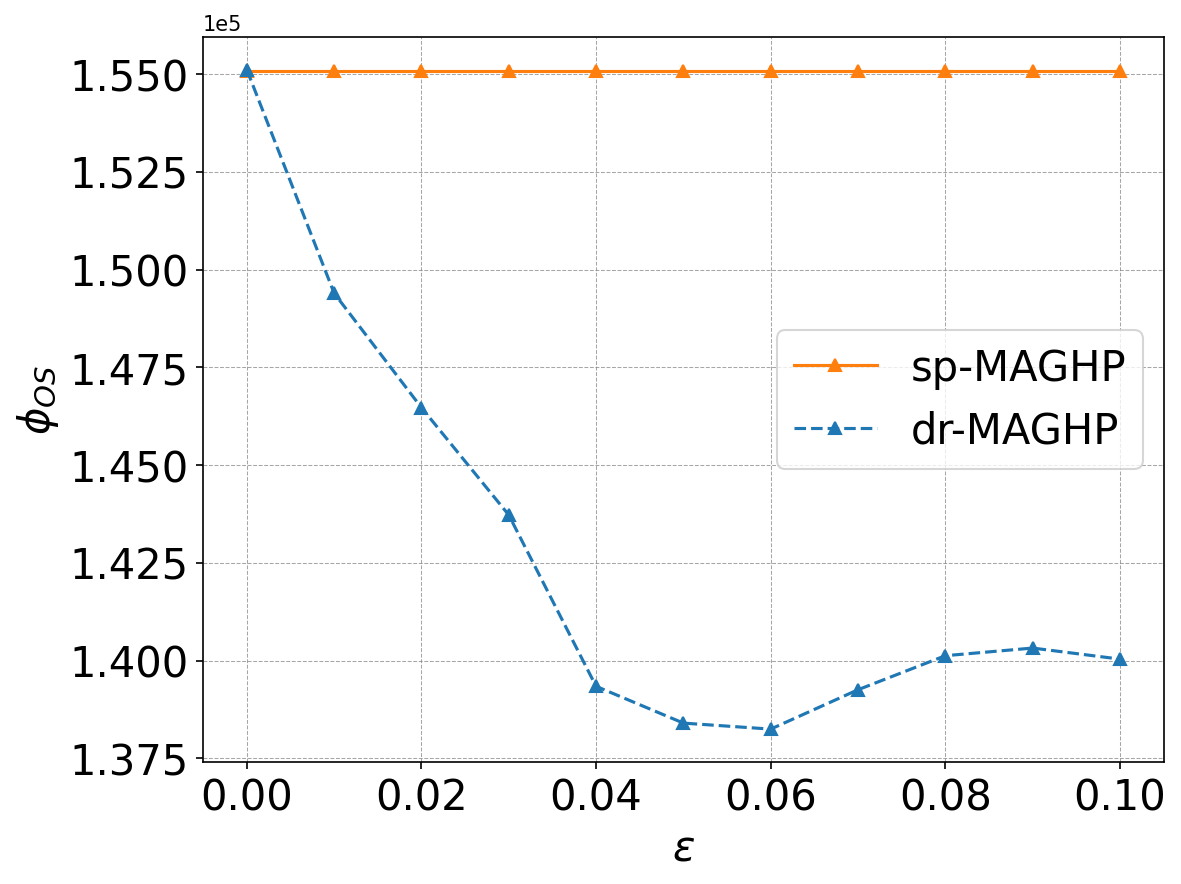}
        \caption{rl = 0.3}
    \end{subfigure}
    \hfill
    \begin{subfigure}[b]{0.49\textwidth}
        \includegraphics[width=\textwidth]{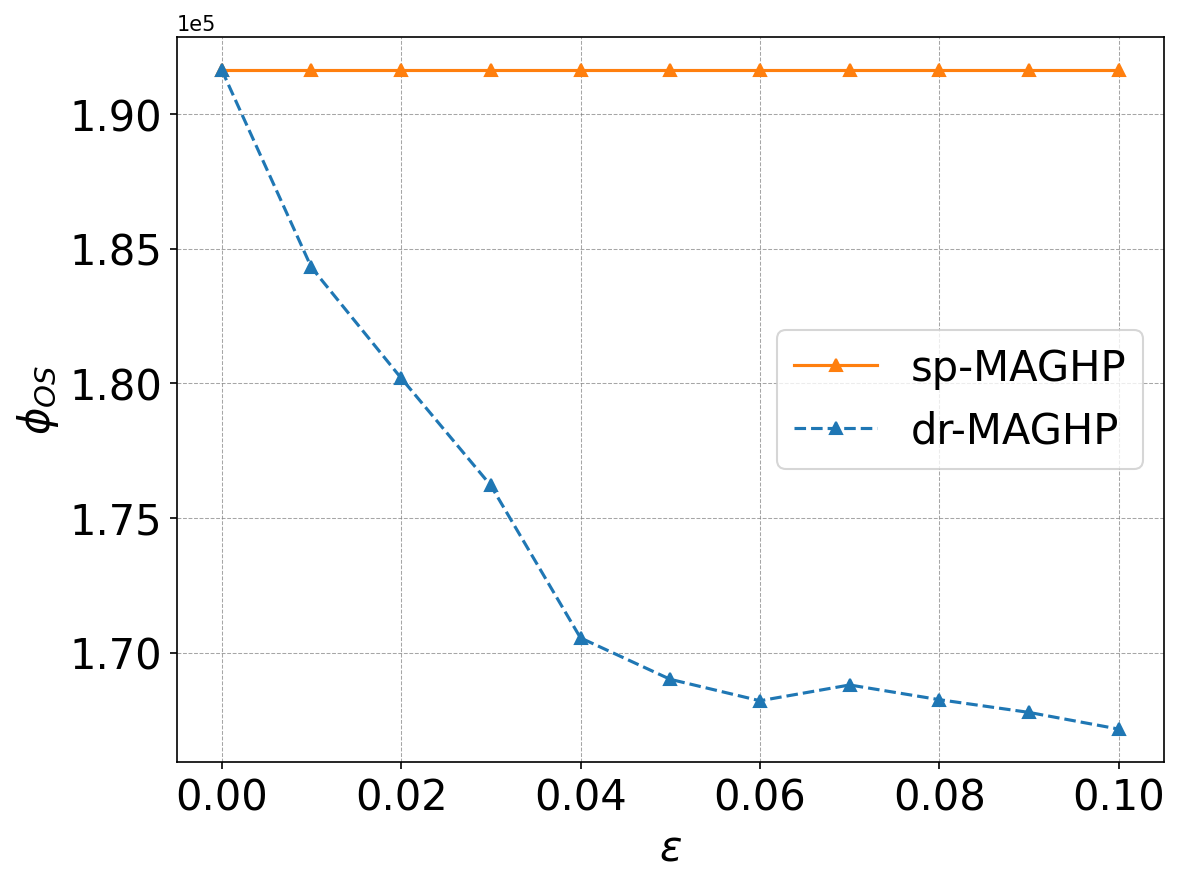}
        \caption{rl = 0.4}
    \end{subfigure}

    \vspace{0.7em}

    \begin{subfigure}[b]{0.5\textwidth}
        \centering
        \includegraphics[width=\textwidth]{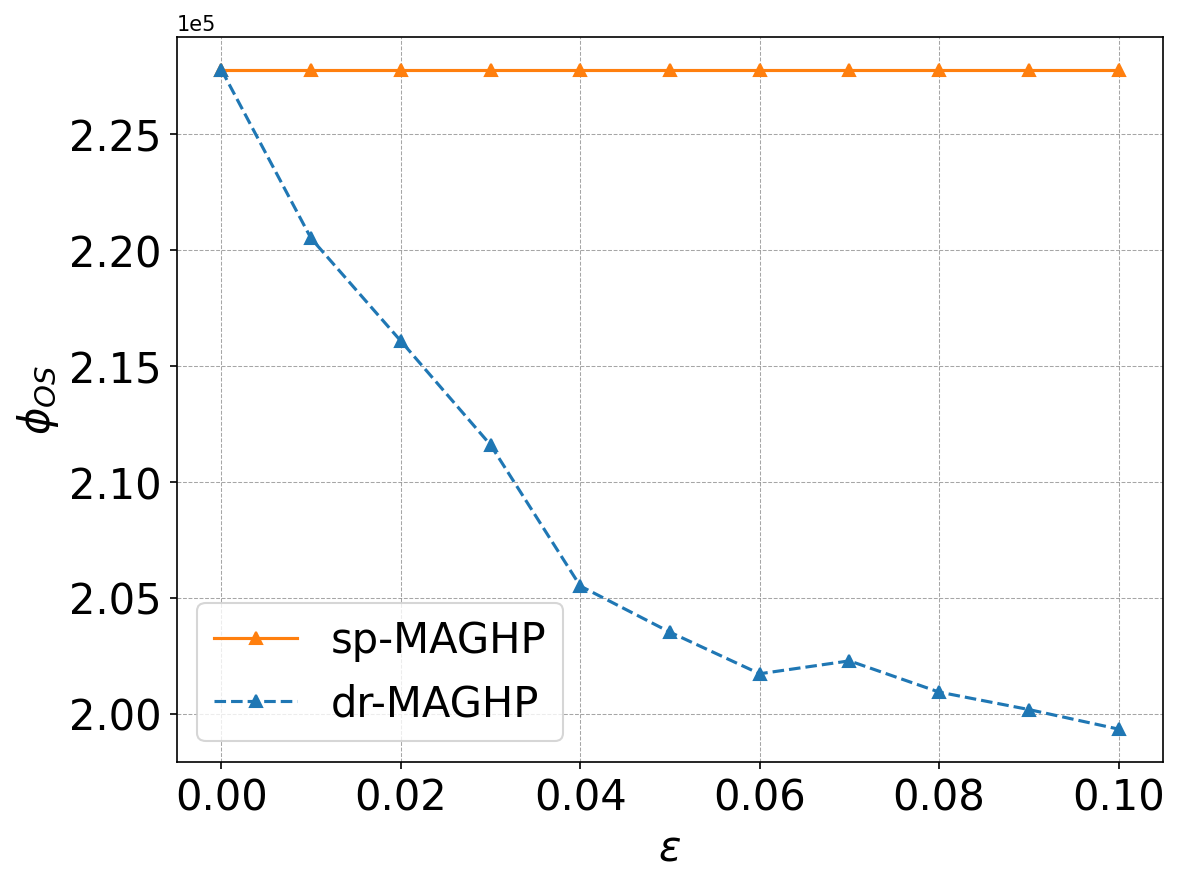}
        \caption{rl = 0.5}
    \end{subfigure}

    \caption{\emph{out-of-sample} performance of \textsc{sp-MAGHP} and \textsc{dr-MAGHP} with different $\epsilon$ values under overestimation scenario on 2019-06-16.}
    \label{fig:dr_cost_eps_overest_20190616}
\end{figure}

\newpage
\bibliographystyle{IEEEtran} 
\small{
\bibliography{main.bib}

\begin{thebibliography}{10}
\providecommand{\url}[1]{#1}
\csname url@samestyle\endcsname
\providecommand{\newblock}{\relax}
\providecommand{\bibinfo}[2]{#2}
\providecommand{\BIBentrySTDinterwordspacing}{\spaceskip=0pt\relax}
\providecommand{\BIBentryALTinterwordstretchfactor}{4}
\providecommand{\BIBentryALTinterwordspacing}{\spaceskip=\fontdimen2\font plus
\BIBentryALTinterwordstretchfactor\fontdimen3\font minus \fontdimen4\font\relax}
\providecommand{\BIBforeignlanguage}[2]{{%
\expandafter\ifx\csname l@#1\endcsname\relax
\typeout{** WARNING: IEEEtran.bst: No hyphenation pattern has been}%
\typeout{** loaded for the language `#1'. Using the pattern for}%
\typeout{** the default language instead.}%
\else
\language=\csname l@#1\endcsname
\fi
#2}}
\providecommand{\BIBdecl}{\relax}
\BIBdecl

\bibitem{GHP1}
P.~B. Vranas, D.~J. Bertsimas, and A.~R. Odoni, ``The multi-airport ground-holding problem in air traffic control,'' \emph{Operations Research}, vol.~42, no.~2, pp. 249--261, 1994.

\bibitem{kicinger2016airport}
R.~Kicinger, J.-T. Chen, M.~Steiner, and J.~Pinto, ``Airport capacity prediction with explicit consideration of weather forecast uncertainty,'' \emph{AIAA JAT}, vol.~24, no.~1, pp. 18--28, 2016.

\bibitem{tittle2013airport}
D.~Tittle, P.~McCarthy, and Y.~Xiao, ``Airport runway capacity and economic development: A panel data analysis of metropolitan statistical areas,'' \emph{Economic Development Quarterly}, vol.~27, no.~3, pp. 230--239, 2013.

\bibitem{avery2015predicting}
\BIBentryALTinterwordspacing
J.~Avery and H.~Balakrishnan, ``Predicting airport runway configuration: A discrete-choice modeling approach,'' in \emph{Eleventh USA/Europe Air Traffic Management Research and Development Seminar}.\hskip 1em plus 0.5em minus 0.4em\relax Lisbon, Portugal: Federal Aviation Administration/EUROCONTROL, June 23-26 2015. [Online]. Available: \url{http://www.atmseminarus.org/seminarContent/seminar11/presentations/509-Balakrishnan_0126150652-PresentationPDF-6-29-15.pdf}
\BIBentrySTDinterwordspacing

\bibitem{vlachou2019simultaneous}
K.~Vlachou, R.~Sharma, and F.~Wieland, ``Simultaneous traffic management initiatives: The double delay problem,'' in \emph{2019 IEEE/AIAA 38th DASC}.\hskip 1em plus 0.5em minus 0.4em\relax IEEE, 2019, pp. 1--6.

\bibitem{Info_centric_NAS}
{FAA}, ``{Initial Concept of Operations for an Info-Centric National Airspace System},'' 2022, accessed: Feb 2024.

\bibitem{SESAR_ATM_MASTER_PLAN}
{SESAR}, ``{European ATM Master Plan 2025 Edition},'' 2025, accessed: Aug 2025.

\bibitem{zhu2021flight}
X.~Zhu and L.~Li, ``Flight time prediction for fuel loading decisions with a deep learning approach,'' \emph{Transportation Research Part C: Emerging Technologies}, vol. 128, p. 103179, 2021.

\bibitem{zhu2022short}
X.~Zhu, Y.~Lin, Y.~He, K.-L. Tsui, P.~W. Chan, and L.~Li, ``Short-term nationwide airport throughput prediction with graph attention recurrent neural network,'' \emph{Frontiers in Artificial Intelligence}, vol.~5, p. 884485, 2022.

\bibitem{Wang2022}
Z.~Wang, C.~Liao, X.~Hang, L.~Li, D.~Delahaye, and M.~Hansen, ``Distribution prediction of strategic flight delays via machine learning methods,'' \emph{Sustainability}, vol.~14, no.~22, 2022.

\bibitem{tamang2025handling}
L.~Tamang, M.~R. Bouadjenek, R.~Dazeley, and S.~Aryal, ``Handling out-of-distribution data: A survey,'' \emph{IEEE Transactions on Knowledge and Data Engineering}, 2025.

\bibitem{filos2020can}
A.~Filos, P.~Tigkas, R.~McAllister, N.~Rhinehart, S.~Levine, and Y.~Gal, ``Can autonomous vehicles identify, recover from, and adapt to distribution shifts?'' in \emph{International Conference on Machine Learning}.\hskip 1em plus 0.5em minus 0.4em\relax PMLR, 2020, pp. 3145--3153.

\bibitem{yang2023generic}
Z.~Yang, X.~He, J.~Zhang, J.~Wu, X.~Xin, J.~Chen, and X.~Wang, ``A generic learning framework for sequential recommendation with distribution shifts,'' in \emph{Proceedings of the 46th International ACM SIGIR Conference on Research and Development in Information Retrieval}, 2023, pp. 331--340.

\bibitem{wang2019drifted}
X.~Wang, Q.~Kang, J.~An, and M.~Zhou, ``Drifted twitter spam classification using multiscale detection test on kl divergence,'' \emph{IEEE Access}, vol.~7, pp. 108\,384--108\,394, 2019.

\bibitem{wu2024distributionally}
H.~Wu, X.~Zhu, S.~Li, Y.~Zhou, L.~Li, and M.~Z. Li, ``Distributionally robust ground delay programs with learning-driven airport capacity predictions,'' \emph{arXiv preprint arXiv:2402.11415}, 2024.

\bibitem{kicinger2012airport}
R.~Kicinger, J.~Krozel, M.~Steiner, and J.~Pinto, ``Airport capacity prediction integrating ensemble weather forecasts,'' in \emph{Infotech Aerospace 2012}, 2012, p. 2493.

\bibitem{gilbo1993airport}
E.~P. Gilbo, ``Airport capacity: Representation, estimation, optimization,'' \emph{IEEE Transactions on Control Systems Technology}, vol.~1, no.~3, pp. 144--154, 1993.

\bibitem{choi2021artificial}
S.~Choi and Y.~J. Kim, ``Artificial neural network models for airport capacity prediction,'' \emph{Journal of Air Transport Management}, vol.~97, p. 102146, 2021.

\bibitem{tien2018using}
S.-L. Tien, C.~Taylor, E.~Vargo, and C.~Wanke, ``Using ensemble weather forecasts for predicting airport arrival capacity,'' \emph{Journal of Air Transportation}, vol.~26, no.~3, pp. 123--132, 2018.

\bibitem{cox2016probabilistic}
J.~Cox and M.~J. Kochenderfer, ``Probabilistic airport acceptance rate prediction,'' in \emph{AIAA Modeling and Simulation Technologies Conference}, 2016, p. 0165.

\bibitem{pang2021data}
Y.~Pang, N.~Xu, and Y.~Liu, ``Data-driven trajectory prediction with weather uncertainties: A bayesian deep learning approach,'' \emph{Transportation Research Part C: Emerging Technologies}, vol. 130, p. 103326, 2021.

\bibitem{zoutendijk2021probabilistic}
M.~Zoutendijk and M.~Mitici, ``Probabilistic flight delay predictions using machine learning and applications to the flight-to-gate assignment problem,'' \emph{Aerospace}, vol.~8, no.~6, p. 152, 2021.

\bibitem{masalonis2004using}
A.~Masalonis, S.~Mulgund, L.~Song, C.~Wanke, and S.~Zobell, ``Using probabilistic demand predictions for traffic flow management decision support,'' in \emph{AIAA Guidance, Navigation, and Control Conference and Exhibit}, 2004, p. 5231.

\bibitem{zhang2023air}
X.~Zhang, J.~Chen, H.~Yang, M.~Li, and Y.~Wang, ``Air traffic density prediction using bayesian ensemble graph attention network (began),'' \emph{Transportation Research Part C: Emerging Technologies}, vol. 152, p. 104140, 2023.

\bibitem{dalmau2024probabilistic}
R.~Dalmau, P.~De~Falco, M.~Spak, and J.~D. Rodriguez-Varela, ``Probabilistic pretactical arrival and departure flight delay prediction with quantile regression,'' \emph{Journal of Air Transportation}, vol.~32, no.~2, pp. 84--96, 2024.

\bibitem{vandal2018prediction}
\BIBentryALTinterwordspacing
T.~Vandal, M.~Livingston, C.~Piho, and S.~Zimmerman, ``Prediction and uncertainty quantification of daily airport flight delays,'' in \emph{Proceedings of The 4th International Conference on Predictive Applications and APIs}, ser. Proceedings of Machine Learning Research, vol.~82.\hskip 1em plus 0.5em minus 0.4em\relax PMLR, 2018, pp. 45--51. [Online]. Available: \url{https://proceedings.mlr.press/v82/vandal18a.html}
\BIBentrySTDinterwordspacing

\bibitem{rodriguez2024air}
Y.~Rodríguez and O.~Díaz~Olariaga, ``Air traffic demand forecasting with a bayesian structural time series approach,'' \emph{Periodica Polytechnica Transportation Engineering}, vol.~52, no.~1, pp. 75--85, 2024.

\bibitem{pang2020probabilistic}
Y.~Pang and Y.~Liu, ``Probabilistic aircraft trajectory prediction considering weather uncertainties using dropout as bayesian approximate variational inference,'' in \emph{AIAA Scitech 2020 Forum}, 2020, p. 1413.

\bibitem{lecchini2006monte}
A.~Lecchini~Visintini, W.~Glover, J.~Lygeros, and J.~Maciejowski, ``Monte carlo optimization for conflict resolution in air traffic control,'' \emph{IEEE Transactions on Intelligent Transportation Systems}, vol.~7, no.~4, pp. 470--482, 2006.

\bibitem{wang2021prediction}
Y.~Wang and Y.~Zhang, ``Prediction of runway configurations and airport acceptance rates for multi-airport system using gridded weather forecast,'' \emph{Transportation Research Part C: Emerging Technologies}, vol. 125, p. 103049, 2021.

\bibitem{odoni1987flow}
A.~R. Odoni, ``The flow management problem in air traffic control,'' in \emph{Flow control of congested networks}.\hskip 1em plus 0.5em minus 0.4em\relax Springer, 1987, pp. 269--288.

\bibitem{bertsimas2008air}
D.~Bertsimas, G.~Lulli, and A.~Odoni, ``The air traffic flow management problem: An integer optimization approach,'' in \emph{International conference on integer programming and combinatorial optimization}.\hskip 1em plus 0.5em minus 0.4em\relax Springer, 2008, pp. 34--46.

\bibitem{bertsimas2000traffic}
D.~Bertsimas and S.~S. Patterson, ``The traffic flow management rerouting problem in air traffic control: A dynamic network flow approach,'' \emph{Transportation Science}, vol.~34, no.~3, pp. 239--255, 2000.

\bibitem{sun2008multicommodity}
D.~Sun and A.~M. Bayen, ``Multicommodity eulerian-lagrangian large-capacity cell transmission model for en route traffic,'' \emph{Journal of guidance, control, and dynamics}, vol.~31, no.~3, pp. 616--628, 2008.

\bibitem{GHP2}
M.~O. Ball, R.~Hoffman, A.~R. Odoni, and R.~Rifkin, ``A stochastic integer program with dual network structure and its application to the ground-holding problem,'' \emph{Operations research}, vol.~51, no.~1, pp. 167--171, 2003.

\bibitem{GHP4}
J.~Chen and D.~Sun, ``Stochastic ground-delay-program planning in a metroplex,'' \emph{Journal of Guidance, Control, and Dynamics}, vol.~41, no.~1, pp. 231--239, 2018.

\bibitem{GHP5}
C.~Chin, M.~Z. Li, K.~Gopalakrishnan, and H.~Balakrishnan, ``Airport ground holding with hierarchical control objectives,'' in \emph{USA-Europe ATM R\&D Seminar}, 2021.

\bibitem{GHP6}
C.~N. Glover and M.~O. Ball, ``Stochastic optimization models for ground delay program planning with equity--efficiency tradeoffs,'' \emph{Transportation Research Part C: Emerging Technologies}, vol.~33, pp. 196--202, 2013.

\bibitem{GHP7}
A.~Jacquillat, ``Predictive and prescriptive analytics toward passenger-centric ground delay programs,'' \emph{Transportation Science}, vol.~56, no.~2, pp. 265--298, 2022.

\bibitem{IATA_WASG_Ed3_2024}
{Airports Council International (ACI) and International Air Transport Association (IATA) and Worldwide Airport Coordinators Group (WWACG)}, \emph{Worldwide Airport Slot Guidelines (WASG), Edition 3, Effective 1 April 2024}, ACI / IATA / WWACG, Montreal / Geneva / Global, Sep. 2023, published jointly by ACI, IATA, and WWACG; available at \url{https://www.iata.org/en/programs/ops-infra/slots/slot-guidelines/} (accessed 09/15/2025).

\bibitem{pellegrini2012metaheuristic}
P.~Pellegrini, L.~Castelli, and R.~Pesenti, ``Metaheuristic algorithms for the simultaneous slot allocation problem,'' \emph{IET Intelligent Transport Systems}, vol.~6, no.~4, pp. 453--462, 2012.

\bibitem{pellegrini2017sosta}
P.~Pellegrini, T.~Boli{\'c}, L.~Castelli, and R.~Pesenti, ``Sosta: An effective model for the simultaneous optimisation of airport slot allocation,'' \emph{Transportation Research Part E: Logistics and Transportation Review}, vol.~99, pp. 34--53, 2017.

\bibitem{ball2020quantity}
M.~O. Ball, A.~S. Estes, M.~Hansen, and Y.~Liu, ``Quantity-contingent auctions and allocation of airport slots,'' \emph{Transportation Science}, vol.~54, no.~4, pp. 858--881, 2020.

\bibitem{zografos2019bi}
K.~G. Zografos and Y.~Jiang, ``A bi-objective efficiency-fairness model for scheduling slots at congested airports,'' \emph{Transportation Research Part C: Emerging Technologies}, vol. 102, pp. 336--350, 2019.

\bibitem{ribeiro2019large}
N.~A. Ribeiro, A.~Jacquillat, and A.~P. Antunes, ``A large-scale neighborhood search approach to airport slot allocation,'' \emph{Transportation Science}, vol.~53, no.~6, pp. 1772--1797, 2019.

\bibitem{DR1}
E.~Delage and Y.~Ye, ``Distributionally robust optimization under moment uncertainty with application to data-driven problems,'' \emph{Operations Research}, vol.~58, no.~3, pp. 595--612, 2010.

\bibitem{DR2}
P.~M. Esfahani and D.~Kuhn, ``Data-driven distributionally robust optimization using the wasserstein metric: Performance guarantees and tractable reformulations,'' \emph{arXiv preprint arXiv:1505.05116}, 2015.

\bibitem{DR3}
K.~Kim, ``Dual decomposition of two-stage distributionally robust mixed-integer programming under the wasserstein ambiguity set,'' \emph{Preprint manuscript}, 2020.

\bibitem{cheramin2022computationally}
M.~Cheramin, J.~Cheng, R.~Jiang, and K.~Pan, ``Computationally efficient approximations for distributionally robust optimization under moment and wasserstein ambiguity,'' \emph{INFORMS Journal on Computing}, vol.~34, no.~3, pp. 1768--1794, 2022.

\bibitem{jiang2019data}
R.~Jiang, M.~Ryu, and G.~Xu, ``Data-driven distributionally robust appointment scheduling over wasserstein balls,'' \emph{arXiv preprint arXiv:1907.03219}, 2019.

\bibitem{hanasusanto2018conic}
G.~A. Hanasusanto and D.~Kuhn, ``Conic programming reformulations of two-stage distributionally robust linear programs over wasserstein balls,'' \emph{Operations Research}, vol.~66, no.~3, pp. 849--869, 2018.

\bibitem{wiesemann2014distributionally}
W.~Wiesemann, D.~Kuhn, and M.~Sim, ``Distributionally robust convex optimization,'' \emph{Operations research}, vol.~62, no.~6, pp. 1358--1376, 2014.

\bibitem{ball2010total}
M.~Ball, C.~Barnhart, M.~Dresner, M.~Hansen, K.~Neels, A.~Odoni, E.~Peterson, L.~Sherry, A.~Trani, and B.~Zou, ``Total delay impact study: a comprehensive assessment of the costs and impacts of flight delay in the united states,'' Federal Aviation Administration, Washington, DC, Tech. Rep., 2010.

\bibitem{mukherjee2012dynamic}
A.~Mukherjee and M.~Hansen, ``Dynamic stochastic optimization model for air traffic flow management with en route and airport capacity constraints,'' \emph{Transportation Research Part C: Emerging Technologies}, vol.~21, no.~1, pp. 61--73, 2012.

\bibitem{wang2017airport}
Y.~Wang, M.~Viscarra, L.~Jin, and M.~V. Bendarkar, ``Airport capacity prediction with explicit consideration of terminal airspace constraints,'' in \emph{17th AIAA Aviation Technology, Integration, and Operations Conference}, 2017, p. 3428.

\bibitem{FAA2016FSM}
\BIBentryALTinterwordspacing
{Federal Aviation Administration}, ``Flight schedule monitor user's guide, version 13.0,'' Technical Report, Feb 2016. [Online]. Available: \url{https://tfmlearning.faa.gov/tfm-training/final_rel_13_tfms_fsm_users_guide_12686.pdf}
\BIBentrySTDinterwordspacing

\bibitem{FAA2017AADC}
\BIBentryALTinterwordspacing
------, ``Airport arrival demand chart information briefing,'' Technical Report, Jul 2017. [Online]. Available: \url{https://tfmlearning.faa.gov/assets/media/tfm-training/TFMS_AADC_Information_Briefing_28Jul2017.pdf}
\BIBentrySTDinterwordspacing

\bibitem{sridhar2011modeling}
B.~Sridhar, K.~S. Sheth, and S.~Grabbe, ``Modeling and optimization in traffic flow management,'' \emph{Proceedings of the IEEE}, vol.~96, no.~12, pp. 2060--2080, 2008.

\bibitem{liu2019using}
Y.~Liu, Y.~Liu, M.~Hansen, A.~Pozdnukhov, and D.~Zhang, ``Using machine learning to analyze air traffic management actions: Ground delay program case study,'' \emph{Transportation Research Part E}, vol. 131, pp. 80--95, 2019.

\bibitem{gilbo1997airport}
E.~P. Gilbo, ``Airport capacity: representation, estimation, optimization,'' \emph{IEEE Transactions on Control Systems Technology}, vol.~5, no.~1, pp. 144--154, 1997.

\bibitem{simaiakis2014balancing}
I.~Simaiakis, M.~Sandberg, and H.~Balakrishnan, ``Balancing runway capacity and delay at airports through pushback rate control: A queuing network approach,'' \emph{Transportation Research Part C: Emerging Technologies}, vol.~44, pp. 378--395, 2014.

\bibitem{ramanujam2009estimation}
V.~Ramanujam and H.~Balakrishnan, ``Estimation of arrival-departure capacity tradeoffs in multi-airport systems,'' in \emph{Proceedings of the 48th IEEE Conference on Decision and Control (CDC)}.\hskip 1em plus 0.5em minus 0.4em\relax IEEE, 2009, pp. 2534--2540.

\bibitem{neufville2013airport}
R.~de~Neufville and A.~R. Odoni, \emph{Airport Systems: Planning, Design, and Management}, 2nd~ed.\hskip 1em plus 0.5em minus 0.4em\relax New York: McGraw-Hill Education, 2013.

\bibitem{allan2001analysis}
S.~Allan, J.~Beesley, J.~Evans, and S.~Gaddy, ``Analysis of delay causality at newark international airport,'' in \emph{4th USA/Europe Air Traffic Management R\&D Seminar}, 2001, pp. 1--11.

\bibitem{renhe2014meteorological}
Z.~Renhe, Q.~Li, and R.~Zhang, ``Meteorological conditions for the persistent severe fog and haze event over eastern china in january 2013,'' \emph{SCES}, vol.~57, pp. 26--35, 2014.

\bibitem{drGHP}
H.~Wu and M.~Z. Li, ``Distributionally robust airport ground holding problem under wasserstein ambiguity sets,'' \emph{arXiv preprint arXiv:2306.09836}, 2023.

\bibitem{delay_cost}
{Airports Council International – North America }, ``{Aircraft Operating and Delay Cost per Enplanement},'' 2014, accessed: Feb 2024.

\bibitem{likas2003global}
A.~Likas, N.~Vlassis, and J.~J. Verbeek, ``The global k-means clustering algorithm,'' \emph{Pattern recognition}, vol.~36, no.~2, pp. 451--461, 2003.

\end{thebibliography}
}
\end{document}